\documentclass[12pt]{amsart}
\usepackage{amsmath, amsthm, amssymb,cite,enumerate}
\usepackage{fullpage}
\usepackage{color}
\usepackage{hyperref}
\usepackage[english]{babel}
\usepackage{framed}

\newtheorem{thm}{Theorem}
\newtheorem{lem}{Lemma}[section]
\newtheorem{prop}[lem]{Proposition}
\newtheorem*{rmk}{Remark}

\newtheorem*{defn}{Definition}
\numberwithin{equation}{section}

\newcommand{\R}{\mathbb{R}}
\newcommand{\D}{\mathcal{D}}

\newcommand{\eps}{\varepsilon}
\newcommand{\norm}[1]{\left\lVert#1\right\rVert}
\newcommand{\fm}[1]{\[\begin{aligned} #1 \end{aligned}\]}
\newcommand{\lra}[1]{\langle #1 \rangle}
\newcommand{\abs}[1]{\left| #1 \right|}
\newcommand{\kh}[1]{\left( #1 \right)}
\newcommand{\wt}[1]{\widetilde{ #1 }}

\begin{document}

\title{Modified wave operators for a scalar quasilinear wave equation satisfying the weak null condition }

\author{Dongxiao Yu}
\address{Department of Mathematics, University of California at Berkeley}
\thanks{}
\email{yudx@math.berkeley.edu}

\begin{abstract}

We prove the existence of the modified wave operators for a scalar quasilinear wave equation satisfying the weak null condition. This is accomplished in three steps. First, we derive a new reduced asymptotic system for the quasilinear wave equation by modifying H\"{o}rmander's method. Next, we construct an approximate solution, by solving our new reduced system given some scattering data at infinite time. Finally, we prove that the quasilinear wave equation has a global solution which agrees with the approximate solution at infinite time.
\end{abstract}

\maketitle
\tableofcontents

\addtocontents{toc}{\protect\setcounter{tocdepth}{1}}

\section{Introduction}

This paper is devoted to the study of the long time dynamics of a scalar quasilinear wave equation in~$\R_{t,x}^{1+3}$, of the form
\begin{equation}\label{qwe}
\wt{g}^{\alpha\beta}(u)\partial_\alpha\partial_\beta u=0.
\end{equation}
Here we use the Einstein summation convention with the sum taken over $\alpha,\beta=0,1,2,3$ with  $\partial_0=\partial_t$, $\partial_i=\partial_{x_i}$, $i=1,2,3$. We assume that  $\wt{g}^{\alpha\beta}(u)$ are smooth functions of $u$, such that $\wt{g}^{\alpha\beta}=\wt{g}^{\beta\alpha}$ and $\wt{g}^{\alpha\beta}(0)\partial_\alpha\partial_\beta=\square=\partial_t^2-\Delta_x$.

This model equation is closely related to General Relativity. The vector-valued version of $\wt{g}^{\alpha\beta}(u)\partial_\alpha\partial_\beta u$ is the principal part of the Einstein equations in wave coordinates. For more physical background for the equation \eqref{qwe}, we refer the readers to \cite{lind,lindrodn,lindrodn2}.

The study of global well-posededness theory of \eqref{qwe} started with Lindblad's paper~\cite{lind2}. Given the initial data \begin{equation}\label{initcond}u(0)=\eps u_0,\ \partial_tu(0)=\eps u_1,\hspace{1cm}\text{where }u_1,u_2\in C_c^\infty(\R^3)\text{ and }\eps >0\text{ small enough},\end{equation} Lindblad conjectured that \eqref{qwe} has a global solution if $\eps$ is sufficiently small. In the same paper, he proved the small data global existence for a special case 
\begin{equation}\label{qwespec}\partial_t^2u-c(u)^2\Delta_xu=0,\hspace{1cm}\text{where }c(0)=1\end{equation}for radially symmetric data.
Later, Alinhac \cite{alin2} generalized  the result to general initial data for \eqref{qwespec}. The small data global existence result to the general case \eqref{qwe} was finally proved by Lindblad in \cite{lind}.

In this paper, we prove the existence of the modified wave operators for \eqref{qwe}, which is closely related to the global well-posedness theory. Precisely, our goal is two-fold. First, we seek to identify a good notion of \emph{asymptotic profile} for this problem, and an associated notion of \emph{scattering data}. Then, for this asymptotic profile, we find a matching solution.

Usually, a global solution to a nonlinear PDE would scatter to a solution to the corresponding linear equation. For example, a  solution to the cubic defocusing three-dimensional NLS scatters to a solution of a linear Schr\"{o}dinger equation; see \cite{tao}. This is however not the case for \eqref{qwe}; its global solution does not decay as  a solution to $\Box u=0$. Thus, a good notion of asymptotic profile is  necessary here. One such candidate is given by a type of asymptotic equations introduced by H\"{o}rmander \cite{horm,horm2,horm3}. From H\"{o}rmander's asymptotic equations, several results on scattering for \eqref{qwe} have been proved; we refer to Lindblad-Schlue \cite{lindschl} and Deng-Pusateri \cite{dengpusa}.

In this paper, we identify a new notion of asymptotic profile by deriving a new reduced system. Our derivation is similar to that of H\"{o}rmander's asymptotic equations, but we choose a different $q$ in the ansatz $u\approx \eps r^{-1}U(s,q,\omega)$, with $s=\eps\ln(t)$, $\omega=x/|x|$ and $r=|x|$. Instead of taking $q=r-t$ as in H\"{o}rmander's derivation, we let $q(t,r,\omega)$ be the solution to the eikonal equation $\wt{g}^{\alpha\beta}(u)q_\alpha q_\beta=0$. By introducing an auxiliary variable $\mu=q_t-q_r$, we are able to derive a first-order ODE system for $\mu$ and $U_q$ in the  coordinate set $(s,q,\omega)$; see \eqref{asypde2} for the reduced system we obtain.

With this new reduced system, we can construct an approximate solution to \eqref{qwe} as follows. First, we solve the new reduced system explicitly with the initial data $(\mu, U_q)|_{s=0}=(-2,A)$. Here $A=A(q,\omega)$ is our \emph{scattering data}. Second, we construct an approximate solution $q(t,r,\omega)$ to the eikonal equation by solving $q_t-q_r=\mu$ and $q(t,0,\omega)=-t$; this equation is an ODE along each characteristic line. Both $s$ and $q$ are now functions of $(t,r,\omega)$, so we obtain a function $U(t,r,\omega)$ from the solution to the reduced system. Here $U(t,r,\omega)$ is our \emph{asymptotic profile}. Third, we define $u_{app}$ which is an approximate solution to \eqref{qwe}. We expect that $u_{app}=\eps r^{-1}U(t,r,\omega)$ in a conic neighborhood of the light cone $\{t=r\}$ and that $u_{app}$ is supported in a slightly larger conic neighborhood of the light cone.

Finally, we show that there is an exact solution to \eqref{qwe} which behaves asymptotically the same as $u_{app}$ as time goes to infinity. Fixing a large time $T>0$, we solve a backward Cauchy problem for $v=u-u_{app}$ with zero data for $t\geq 2T$, such that $v+u_{app}$ solves $\eqref{qwe}$ for $t\leq T$. We then prove that $v=v^T$ converges to some function $v^\infty$ as $T\to\infty$. Now $u^\infty=v^\infty+u_{app}$ is a solution to \eqref{qwe} which matches the asymptotic profile at  infinite time. This shows the existence of the modified wave operators.

\subsection{Background}
The equation \eqref{qwe} is a special case for a general scalar nonlinear wave equation in~$\R^{1+3}_{t,x}$
\begin{equation}\label{nlw}
\square u=F(u,\partial u,\partial^2u).
\end{equation}
Here 
\begin{equation}\label{nonlinearity}F(u,\partial u,\partial^2u)=\sum a_{\alpha\beta}\partial^\alpha u\partial^\beta u+O(|u|^3+|\partial u|^3+|\partial^2 u|^3).\end{equation}
The sum in \eqref{nonlinearity} is taken over    all multiindices $\alpha,\beta$ with $|\alpha|\leq |\beta|\leq 2$, $|\beta|\geq1$ and $|\alpha|+|\beta|\leq 3$.

Since 1980s, there have been many results  on the lifespan of the solutions to the Cauchy problem \eqref{nlw} with initial data \eqref{initcond}. In \cite{john,john2}, John proved that \eqref{nlw} does not necessarily have a global solution for all $t\geq 0$: any nontrivial solution to  $\square u=u_t\Delta u$ or $\square u=u_t^2$ blows up in finite time. In contrast, \eqref{nlw} in $\R^{1+d}$ for $d\geq 4$ has small data global existence, proved by H\"{o}rmander \cite{horm3}. For arbitrary nonlinearities in three space dimensions, the best result on the lifespan is the almost global existence: the solution exists for $t\leq e^{c/\eps}$, for sufficiently small $\eps$ and some constant $c>0$. The almost global existence for \eqref{nlw} was proved by Lindblad \cite{lind3}. We also refer to  John and Klainerman \cite{johnklai},  Klainerman \cite{klai}, and  H\"{o}rmander \cite{horm2,horm} for some earlier work on  almost global existence.

In contrast to the finite-time blowup in John's examples, it was proved by Klainerman~\cite{klai3} and by Christodoulou~\cite{chri} that if  the null condition is satisfied, then  \eqref{nlw} has  small data global existence. The null condition was first introduced by Klainerman~\cite{klai2}. It states that for each $0\leq m\leq n\leq 2$ with $m+n\leq 3$, we have 
\begin{equation}\label{nullcond}A_{mn}(\omega):=\sum_{|\alpha|=m,|\beta|=n}a_{\alpha\beta}\widehat{\omega}^{\alpha}\widehat{\omega}^{\beta}=0,\hspace{1cm}\text{for all }\widehat{\omega}=(-1,\omega)\in\R\times\mathbb{S}^2.\end{equation}
Note that the null condition is sufficient but not necessary for the small data global existence. For example, the null condition fails for \eqref{qwe} in general, but \eqref{qwe} still has  small data global existence.

Later, in \cite{lindrodn,lindrodn2}, Lindblad and Rodnianski introduced the weak null condition. To state the weak null condition, we start with the asymptotic equations first introduced by H\"{o}rmander in \cite{horm2,horm,horm3}. We make the ansatz\begin{equation}\label{ansatz}
u(t,x)\approx \frac{\eps}{r}U(s,q,\omega),\hspace{1cm}r=|x|,\ \omega_i=x_i/r,\ s=\eps\ln(t),\ q=r-t.
\end{equation}
Plug this ansatz into \eqref{nlw} and we can derive the following asymptotic PDE for $U(s,q,\omega)$
\begin{equation}\label{asypde11}2\partial_s\partial_q U+\sum A_{mn}(\omega)\partial_q^mU\partial_q^nU=0.\end{equation}
Here $A_{mn}$ is defined in \eqref{nullcond} and the sum is taken over $0\leq m\leq n\leq 2$ with $m+n\leq 3$. We say that the weak null condition is satisfied if \eqref{asypde11} has a global solution for all $s\geq 0$ and if the solution and all its derivatives grow at most exponentially in $s$, provided that the initial data decay sufficiently fast in $q$. In the same papers, Lindblad and Rodnianski made a conjecture that the weak null condition is sufficient for small data global existence. To the best of the author's knowledge, this conjecture  remains open until today.

There are three remarks about the weak null condition and the corresponding conjecture. First, the weak null condition is weaker than the null condition. In fact, if the null condition is satisfied, then \eqref{asypde11} becomes $\partial_s\partial_q U=0$. Second, though the conjecture   remains open, there are many examples of \eqref{nlw} satisfying the weak null condition and admitting small data global existence at the same time. The equation \eqref{qwe} is one of several such examples: the small data global existence for  \eqref{qwe} has  been proved by Lindblad  \cite{lind}; meanwhile, the asymptotic equation \eqref{asypde11} now becomes
\begin{equation}\label{asypde1}2\partial_s\partial_q U=G(\omega)U\partial_q^2U,\end{equation}
where \fm{G(\omega):=g^{\alpha\beta}\widehat{\omega}_\alpha\widehat{\omega}_\beta,\hspace{1cm}g^{\alpha\beta}=\frac{d}{du}\wt{g}^{\alpha\beta}(u)|_{u=0},\ \widehat{\omega}=(-1,\omega)\in\R\times\mathbb{S}^2,}
whose solutions exist globally in $s$ and satisfy the decay requirements, so \eqref{qwe} satisfies the weak null condition. There are also many examples  violating the weak null condition and admitting finite-time blowup at the same time. Two of such examples are $\square u=u_t\Delta u$ and $\square u=u_t^2$: the corresponding asymptotic equations are  $(2\partial_s-U_q\partial_q)U_q=0$ (Burger's equation) and $\partial_s U_q=U_q^2$, respectively, whose solutions are known to blow up in finite time. Third, in the recent years, Keir has made  important progress. In \cite{keir}, he proved  the small data global existence  for a large class of quasilinear wave equations satisfying the weak null condition, significantly enlarging upon the class of equations for which global existence is known. His proof also applies to \eqref{qwe}. In \cite{keir2},  he proved that if the solutions to the asymptotic system are bounded (given small initial data) and  stable against rapidly decaying perturbations, then the corresponding system of nonlinear wave equations admits  small data global existence.

\subsection{Asymptotic equations}\label{subsec1.2}

Instead of working with H\"{o}rmander's asymptotic system \eqref{asypde1} directly, in this paper we will construct a new system of asymptotic equations. Our analysis starts as in H\"{o}rmander's derivation in \cite{horm2,horm,horm3}, but diverges at a key point: the choice of $q$ is different. One may contend from the paper that this new system is more accurate than \eqref{asypde1}, in that it both describes the long time evolution and contains full information about it. In addition, if we choose the initial data appropriately, our reduced system will reduce to linear first order ODEs on $\mu$ and $U_q$, so it is easier to solve it than to solve \eqref{asypde1}.

To derive the new equations, we still make the ansatz  \eqref{ansatz}, but now we replace $q=r-t$ with a solution $q(t,r,\omega)$ to the eikonal equation related to \eqref{qwe}\begin{equation}\label{eikeqn}\wt{g}^{\alpha\beta}(u)\partial_\alpha q\partial_\beta q=0.\end{equation}
In other words, $q(t,r,\omega)$ is an optical function. There are two reasons why we choose $q$ in this way. First, if we plug $u=\eps r^{-1}U(s,q,\omega)$  in \eqref{qwe} where $q(t,r,\omega)$ is an arbitrary function, we get two terms in the expansion \fm{\eps r^{-1}\wt{g}^{\alpha\beta}(u) q_{\alpha\beta} U_q+\eps r^{-1}\wt{g}^{\alpha\beta}(u) q_{\alpha}q_{\beta} U_{qq}.}All the other terms  either decay faster than $\eps^2 r^{-2}$ for $t\approx r\to\infty$, or do not contain $U$ itself (but may contain $U_{q},U_{qq},U_{sq}$ and etc.). If $q$ satisfies the eikonal equation, then the second term vanishes. From the eikonal equation, we can also prove that the first term is approximately equal to a function depending on $U_q$ but not on $U$. Thus, in contrast to the second-order PDE \eqref{asypde1} for $U$, we expect to get a first-order ODE for $U_q$ which is simpler.

Second, the eikonal equations have been used in the previous works on the small data global existence for \eqref{qwe}. In  \cite{alin2}, Alinhac followed the method used in  Christodoulou and Klainerman \cite{chriklai}, and adapted the vector fields to the characteristic surfaces, i.e.\ the level surfaces of solutions to the eikonal equations. In \cite{lind}, Lindblad considered the radial eikonal equations when he derived the pointwise bounds of solutions to \eqref{qwe}. When they derived the energy estimates, both Alinhac and Lindblad considered a weight $w(q)$ where $q$ is an approximate solution to the eikonal equation. Their works suggest that the  eikonal equation  plays an important role when we study the long time behavior of solutions to \eqref{qwe}.

Since $u$ is  unknown, it is difficult to solve \eqref{eikeqn} directly. Instead,  we introduce a new auxiliary function $\mu=\mu(s,q,\omega)$ such that $q_t-q_r=\mu$. From \eqref{eikeqn}, we can express $q_t+q_r$ in terms of $\mu$ and $U$, and then solve for all partial derivatives of $q$, assuming that all the angular derivatives are negligible. Then from \eqref{qwe}, we can derive the following asymptotic equations for $\mu(s,q,\omega)$ and $U(s,q,\omega)$:
\begin{equation}\label{asypde2}\left\{
\begin{array}{l}\displaystyle
\partial_s\mu=-\frac{1}{4}G(\omega)\mu^2 U_q,\\[.2cm]
\displaystyle\partial_sU_q=\frac{1}{4}G(\omega) \mu U_q^2.
\end{array}\right.
\end{equation}
The derivation of these two equations is given in Section \ref{sae}.

To solve \eqref{asypde2}, we need to assign the initial data at $s=0$. To choose $\mu|_{s=0}$, we use the gauge freedom. Note that if $q_t-q_r=\mu$ and if $\wt{q}=F(q,\omega)$, then we have $\wt{q}_t-\wt{q}_r=(\partial_qF)\mu$. Thus, by choosing the function $F$ appropriately, we can prescribe $\mu|_{s=0}$ freely. We set $\mu|_{s=0}=-2$ since we expect $q\approx r-t$. The initial data of $U_q$ can be chosen arbitrarily, so we  set $U_q(0,q,\omega)=A(q,\omega)$ for an arbitrary $A(q,\omega)\in C_c^\infty(\R\times\mathbb{S}^2)$. Here $A(q,\omega)$ is defined as the \emph{scattering data} for our result. Note that \eqref{asypde2} implies that $\partial_s(\mu U_q)=0$, so $(\mu U_q)(s,q,\omega)=-2A(q,\omega)$. Then,  \eqref{asypde2} is reduced to a linear first order ODE system
\fm{\left\{
\begin{array}{l}\displaystyle
\partial_s\mu=\frac{1}{2}G(\omega)A(q,\omega)\mu,\\[.2cm]
\displaystyle\partial_sU_q=-\frac{1}{2}G(\omega)A(q,\omega)U_q.
\end{array}\right.}
To uniquely solve $U$ from $U_q$, we also assume that $\lim_{q\to-\infty}U(s,q,\omega)=0$. Now we obtain an explicit solution $(\mu,U)$ to our reduced system \eqref{asypde2}.

To construct an approximate solution,  we make a change of coordinates. For a small $\eps>0$,  we set $s=\eps\ln(t)-\delta$, where $\delta>0$ is  a sufficiently small constant to be chosen. We remark that this choice of $s$ is related to the almost global existence, since now $s=0$ if and only if $t= e^{\delta/\eps}$. In fact, when $t\leq e^{\delta/\eps}$, we expect the solution to \eqref{qwe} behaves as a solution to $\Box u=0$, so our asymptotic equations play a role only when $t\geq e^{\delta/\eps}$. Let $q(t,r,\omega)$ be the solution to \fm{q_t-q_r=\mu(\eps\ln(t)-\delta,q(t,r,\omega),\omega),\hspace{1cm}q(t,0,\omega)=-t.} We can use the method of characteristics to solve this equation. Then, any function of $(s,q,\omega)$ induces a new function of $(t,r,\omega)$. With an abuse of notation, we set \fm{U(t,r,\omega)=U(\eps\ln(t)-\delta,q(t,r,\omega),\omega).} Here $U(t,r,\omega)$ is the \emph{asymptotic profile}. We can prove that, near the light cone $\{t=r\}$, $\eps r^{-1}U(t,r,\omega)$ is an approximate solution to \eqref{qwe}, and $q(t,r,\omega)$ is an approximate optical function, i.e.\ an approximate solution to the eikonal equation corresponding with the metric~$\wt{g}^{\alpha\beta}(\eps r^{-1}U)$. See Section \ref{sas} for the explicit formulas and the estimates for $q$ and~$U$.

\subsection{The main result}

Given the asymptotic equations \eqref{asypde2}, we can ask the following two questions. First, given a scattering data $A(q,\omega)$, can we use \eqref{asypde2} to construct an exact solution to \eqref{qwe} which has this scattering data at infinite time? Second, as time goes to infinity, can any small global solution to the Cauchy problem \eqref{qwe} with \eqref{initcond} be well approximated by a solution to our reduced system \eqref{asypde2}? For example, can we recover the scattering data $A(q,\omega)$, approximate optical function $q(t,r,\omega)$ and asymptotic profile $U(t,r,\omega)$ from an exact solution?  In  scattering theory, the first problem is the existence of the (modified) wave operators, and the second one is  asymptotic completeness. We remark that these two questions have also been formulated and studied for many other nonlinear PDEs. For example, in the setting of nonlinear Schr\"{o}dinger equations, the existence of wave operators and asymptotic completeness have been formulated in  many texts such as \cite{ginivelo}. Scattering theory for NLS has also been studied; we refer to \cite{tao,caze} for a collection of such results.

For the equation \eqref{qwe}, there are some previous results on these two questions. In \cite{lindschl}, Lindblad and Schlue  proved the existence of the wave operators for the semilinear models of Einstein's equations. In \cite{dengpusa}, Deng and Pusateri used the original H\"{o}rmander's asymptotic system \eqref{asypde1} to prove a partial scattering result for \eqref{qwe}.  In their proof, they applied the spacetime resonance method; we refer to \cite{pusashat,pusa2} for some earlier applications of this method to the first order systems of wave equation.  To the author's knowledge, there is no previous result on the modified wave operators for \eqref{qwe}. 

In this paper, we will answer the first question, i.e.\ concerning  the existence of the modified wave operators. Let $Z$ be one of the commuting vector fields: translations $\partial_\alpha$, scaling $t\partial_t+r\partial_r$, rotations $x_i\partial_j-x_j\partial_i$ and Lorentz boosts $x_i\partial_t+t\partial_i$.  Our main theorem is the following.

\begin{thm}\label{mthm1}
Consider a scattering data $A(q,\omega)\in C_c^\infty(\R\times\mathbb{S}^2)$ where ${\rm supp}(A)\subset [-R,R]\times \mathbb{S}^2$ for some $R\geq 1$.  Fix an integer $N\geq 2$ and any sufficiently small $\eps>0$ depending on~$A$ and~$N$. Let $q(t,r,\omega)$  and $U(t,r,\omega)$ be the associated approximate optical function and  asymptotic profile. Then, there is a $C^{N}$ solution $u$ to \eqref{qwe} for $t\geq 0$ with the following properties:
\begin{enumerate}[{\normalfont (i)}]
\item The solution vanishes for $|x|=r\leq t-R$.
\item The solution satisfies the energy bounds: for all $|I|\leq N-1$ and all $t\gg_R1$, we have
\fm{\norm{\partial Z^I(u-\eps r^{-1}U)(t)}_{L^2(\{x\in\R^3:\ |x|\leq 5t/4\})}+\norm{\partial Z^Iu(t)}_{L^2(\{x\in\R^3:\ |x|\geq 5t/4\})}\lesssim_I \eps t^{-1/2+C_I\eps}.}
\item The solution satisfies the  pointwise bounds: for all $(t,r,\omega)$ with $t\gg_R 1$, we have
\fm{|(\partial_t-\partial_r)u+2\eps r^{-1}A(q(t,r,\omega),\omega)|\lesssim \eps t^{-3/2+C\eps}.}
Moreover, for all $|I|\leq N-1$ and all $(t,x)$ with $t\gg_R1$, \fm{|\partial Z^I(u-\eps r^{-1}U)(t,x)|\chi_{|x|\leq 5t/4}+|\partial Z^Iu(t,x)|\chi_{|x|\geq 5t/4}\lesssim_I \eps t^{-1/2+C_I\eps}\lra{t+r}^{-1}\lra{t-r}^{-1/2},}
\fm{|Z^I(u-\eps r^{-1}U)(t,x)|\chi_{|x|\leq 5t/4}+| Z^Iu(t,x)|\chi_{|x|\geq 5t/4}\lesssim_I \min\{\eps t^{-1+C_I\eps},\eps t^{-3/2+C_I\eps}\lra{r-t}\}.}
\end{enumerate}

\end{thm}\rm

\rmk{\rm We have several remarks on the main theorem.

(1) The solution in the main theorem is unique in the following sense. Suppose $N\geq 7$. Suppose $u_1,u_2$ are two $C^N$ solutions to \eqref{qwe}, such that both of them satisfy the energy bounds and pointwise bounds in the main theorem. Then, we have $u_1=u_2$, assuming $\eps\ll 1$. We also remark that $u$ does not depend on the value $5/4$ in the estimates: for each fixed $\kappa>1$, if $u_{\kappa}$ is a solution satisfying all the estimates above with $5/4$ replaced by $\kappa$, then $u=u_{\kappa}$ for $\eps\ll_\kappa 1$, where $u$ is the unique solution from the main theorem. We will prove these statements after  the proof of the main theorem.

(2) By the main theorem, we have the following pointwise bound near the light cone (e.g.\ when $|t-r|\lesssim t^{C\eps}$): \begin{equation}\label{ptbdd}|\partial Z^I(u-\eps r^{-1}U)(t,x)|+|Z^I(u-\eps r^{-1}U)(t,x)|\lesssim_I  \eps t^{-3/2+C_I\eps}.\end{equation}
Note that, for the free constant coefficient linear wave equation, we can prove a  stronger pointwise estimate with $t^{-3/2+C_I\eps}$ replaced by $t^{-2}$ on the right hand side. This is suggested by the fact that the solution to the forward Cauchy problem $\Box w=0$ with compactly supported initial data satisfies such a stronger pointwise estimate  (see Theorem 6.2.1 in \cite{horm}). In our construction, we can achieve this stronger estimate if we add an additional assumption $\int_{-\infty}^{\infty}A(q,\omega)\ dq=0$  on the scattering  data. In fact, this assumption implies that $U$, i.e.\ the Friedlander radiation field, is compactly supported for fixed time; see \eqref{ueqn2} for the definition of $U$. Then, our approximate solution $u_{app}$, defined by \eqref{uappdefn2}, is supported in $\{|r-t|\lesssim 1\}$ (compared with $\{C^{-1}t\leq r\leq t+C\}$ in the general case). This leads to  \fm{\norm{\wt{g}^{\alpha\beta}(u_{app})\partial_{\alpha}\partial_\beta u_{app}(t)}_{L^2(\R^3)}=O(\eps t^{-2})} while in general we only have $O(\eps t^{-3/2+C_I\eps})$. This new estimate would lead to a better energy estimate and thus a better pointwise bound; see the proofs in Section~\ref{smp}. In the general case \eqref{qwe}, the author tends to believe that the exponent $-3/2$ cannot be improved. The dependence on $s=\eps\ln(t)-\delta$ of $U$ prevents us from making $U(t,\cdot)$ compactly supported for all $t$ by only putting restrictions on $A(q,\omega)$. To resolve this problem, we introduce a cutoff function $\psi(r/t)$ (see \eqref{uappdefn2}) which unavoidably causes the loss of power of $t$.

(3) In the main theorem, we assume that the scattering data $A(q,\omega)$ is in $C_c^\infty(\R\times\mathbb{S}^2)$. This assumption can be relaxed.  In fact, instead of  $A\in C^\infty$, we only need $A\in C^{N'}$ where $N'\gg_N1$; instead of $A$ having a compact support in $q$, we can assume  $A\equiv 0$ for $q\leq -R$ and $\partial_q^a\partial_\omega^cA=O(\lra{q}^{-a-1-\delta_0})$, for some fixed $R\geq 1$ and $\delta_0>0$ depending on $A$. We remark that the main theorem remains valid under these weaker assumptions, but proving it would require a more delicate analysis and substantial changes of the arguments in the present paper. For example, in Section \ref{sas},  we would need to abandon the assumption $|q(t,r,\omega)|\leq R$ and extend all the estimates  to a larger region $t\sim r$. However, our goal of this paper is to show the power of the new reduced system \eqref{asypde2}, and we believe that this is  already achieved under the strong assumption $A\in C_c^\infty$. Thus, the author prefers to keep such a strong assumption in this paper for simplicity. The author plans to give a detailed proof under the weaker assumptions stated above in his future dissertation.

Note that the assumption $A\equiv 0$ for $q\leq -R$ is necessary in our proof. It guarantees that both the asymptotic profile $U(t,r,\omega)$ and the exact solution $u(t,x)$ in the main theorem vanish for $r-t\leq -R$. Such a property is essential for the Poincar$\acute{\rm e}$'s lemmas; see Section~\ref{sec5.2}.

Here our decay assumption  $\partial_q^a\partial_\omega^cA=O(\lra{q}^{-a-1-\delta_0})$ is motivated by Lindblad and Schlue~\cite{lindschl}. In  \cite{lindschl}, it is assumed that $(\lra{q}\partial_q)^a\partial_\omega^cF_0=O(\lra{q}^{-\gamma})$ for some $\gamma\in(1/2,1)$, where $F_0$ is their radiation field. For a linear wave equation,  in our setting in this paper, we have $q=r-t$ and $-2A=U_q=\partial_qF_0$, so we expect $-2\partial_q^a\partial_\omega^cA=\partial_q^{a+1}\partial_\omega^cF_0=O(\lra{q}^{-a-1-\delta_0})$.

(4) This paper can be viewed as a preparation for the study on  asymptotic completeness and scattering for the forward Cauchy problem \eqref{qwe} and \eqref{initcond}. To achieve this goal, we should consider whether the setting in this paper fits in a forward Cauchy problem. For example, when we choose the initial data of $\mu$ at $s=0$, we do not have any restriction if we only consider the modified wave operator problem. Thus, we may set $\mu|_{s=0}\equiv -2$ merely for simplicity. However, in our future work, we hope to derive an asymptotic equation on $\mu=q_t-q_r$ for some optical function $q(t,x)$. In general, we do not have $q_t-q_r\equiv -2$ at a fixed time corresponding to $s=0$. Why can we set $\mu|_{s=0}\equiv -2$ in the modified wave operator problem while it does not hold in a forward Cauchy problem? We need the gauge freedom to explain our setting; see the discussions below \eqref{asypde2}. Another example is our choice of $s$. In our construction of asymptotic profile, we take $s=\eps\ln(t)-\delta$. In the modified wave operator problem, this choice is no different from $s=\eps\ln(t)$ or $s=\eps\ln(t+1)$. In the forward scattering problem, however, the choice $s=\eps\ln(t)-\delta$ is better. Recall that we solve our reduced system for $s\geq 0$, and that we expect our reduced system to play a role only when $t\geq e^{\delta/\eps}$. Under our choice of $s$, these two inequalities are equivalent to each other.
}

\subsection{Idea of the proof}\rm
Here we outline the main idea of the construction of $u$ in Theorem~\ref{mthm1}. Roughly speaking,  our starting point is the ideas from both Lindblad \cite{lind} and Lindblad-Schlue \cite{lindschl}. To construct a matching global solution, we follow the idea in Lindblad-Schlue \cite{lindschl}: we solve a backward Cauchy problem with some initial data at $t=T$ and then send $T$ to infinity. However, the backward Cauchy problems in \cite{lindschl} are of simpler form, and their solutions can be constructed by Duhamel's formula explicitly. Here, our backward Cauchy problem is quasilinear, and it is necessary to prove that the solution does exist for all $0\leq t\leq T$. We follow the proof of the small data global existence in \cite{lind}: we use a continuity argument with the  help of the adapted energy estimates and  Poincar$\acute{\rm e}$'s lemma.

We now provide more detailed descriptions of the proof. First, we construct an approximate solution to \eqref{qwe}. Let $q(t,r,\omega)$  and $U(t,r,\omega)$ be the  approximate optical function and  asymptotic profile associated to some scattering data $A(q,\omega)$. We set
\begin{equation}\label{uappdefn2}u_{app}(t,x)=\eps r^{-1}\eta(t)\psi(r/t)U(\eps\ln(t)-\delta,q(t,r,\omega),\omega)\end{equation} for all $t\geq 0$ and $x\in\R^3$. Here $\psi\equiv 1$ when $|r-t|\leq t/4$ and $\psi\equiv 0$ when $|r-t|\geq t/2$, which is used to localize $\eps r^{-1}U$ near the light cone $\{r=t\}$; $\eta$ is a cutoff function such that $\eta\equiv 0$ for $t\leq 2R$, which is used to remove the singularity at $|x|=0$ and $t=0$. We can check that $u_{app}$ is a good approximate solution to \eqref{qwe} in the sense that 
\fm{\wt{g}^{\alpha\beta}(u_{app})\partial_\alpha\partial_\beta u_{app}=O(\eps t^{-3+C\eps}),\hspace{1cm}t\gg_R1.}

Next we seek to construct an exact solution matching $u_{app}$ at infinite time. Fixing a large time $T$, we consider the following equation \begin{equation}\label{introeqn}\wt{g}^{\alpha\beta}(u_{app}+v)\partial_\alpha\partial_\beta v=-\chi(t/T)\wt{g}^{\alpha\beta}(u_{app}+v)\partial_\alpha\partial_\beta u_{app},\ t>0;\hspace{1cm} v\equiv 0,\ t\geq 2T.\end{equation}Here $\chi\in C^\infty(\R)$ satisfies $\chi(t/T)=1$ for $t\leq T$ and $\chi(t/T)=0$ for $t\geq 2T$. Note that  $u_{app}+v$ is now an exact solution to \eqref{qwe} for $t\leq T$. In Section \ref{smp} we prove that, if $\eps$ is sufficiently small, then \eqref{introeqn} has a solution $v=v^T$ for all $t\geq 0$ which satisfies some decay in energy as $t\to\infty$. To prove this, we use a continuity argument. The proof relies on the energy estimates and Poincar$\acute{\rm e}$'s lemma, which are established in Section \ref{sep}. Note that   the small constant $\delta>0$ is not chosen until  the proof of the  Poincar$\acute{\rm e}$'s lemma, and we remark that $\delta$ depends only on the scattering data $A(q,\omega)$. We also remark that the energy estimates and Poincar$\acute{\rm e}$'s lemma in our paper are closely related to those in \cite{lind,alin2}.

Finally we prove in Section \ref{slim} that $v^T$ does converge to some $v^\infty$ in suitable function spaces, as $T\to\infty$. Thus we obtain a global solution $u_{app}+v^\infty$ to \eqref{qwe} for $t\geq 0$, such that it ``agrees with'' $u_{app}$ at infinite time, in the sense that the energy of $v^\infty$ tends to $0$ as $t\to\infty$. By the Klainerman-Sobolev inequality, we can derive the pointwise bounds in the main theorem from the estimates for the energy of $v^\infty$.

Note that to obtain a candidate for $v^\infty$, we have a more natural choice of PDE than \eqref{introeqn}. We may consider the Cauchy problem \eqref{qwe} for $t\leq T$ with initial data $(u_{app}(T),\partial_tu_{app}(T))$. The problem with such a choice is that for $u_{app}$ constructed above, $Z^I(u-u_{app})(T)$ does not seem to have a good decay in $T$ if  $Z^I$ only contains the scaling $S=t\partial_t+r\partial_r$ and Lorentz boosts $\Omega_{0i}=t\partial_i+x_i\partial_t$.  For example, we can consider the linear wave equation $\square u=0$. We set $v=u-u_{app}$, then $v=v_t=0$ at $t=T$. Then, at $t=T$ we have $S^2v=t^2v_{tt}=-t^2\square u_{app}$. However, in the linear case, $u_{app}=\eps r^{-1}F_0(r-t,\omega)$ for $t\approx r$ and thus $\square u_{app}=O(\eps r^{-3})$. The power $-3$ cannot be improved, so we can only get $S^2v=O(\eps r^{-1})$ for $t\approx r$, while we expect $S^2v=O(\eps r^{-3/2+C\eps})$ for $t\approx r$ from Theorem \ref{mthm1}. Similarly, the same applies for $S^kv$ if $k\geq 3$.  In the linear case, one possible way to deal with this difficulty is to consider more terms in the asymptotic expansion of the solutions, say take
\fm{u_{app}=\sum_{n=0}^N\frac{\eps}{r^{n+1}}F_n(r-t,\omega)}where $F_0$ is the usual Friedlander radiation field, and $F_n$ satisfies some PDE based on $F_{n-1}$. This method was used by Lindblad and Schlue in their construction. However,  it does not seem to work in the quasilinear case, since we do not have such a good asymptotic expansion for a solution to  \eqref{qwe}. In this paper, we avoid such a difficulty by considering a variant \eqref{introeqn} of \eqref{qwe}. Such a difficulty does not appear in \eqref{introeqn}, since  $v\equiv 0$ for all $t\geq 2T$.

\subsection{Acknowledgement}
The author would like to thank his advisor, Daniel Tataru, for
suggesting this problem and for many helpful discussions. The author would also like to thank the anonymous reviewers for their valuable comments and suggestions on this paper. This research was partially supported by a James H.\ Simons Fellowship and by the NSF grant DMS-1800294.

\section{Preliminaries}\label{spre}

\subsection{Notations}We use $C$ to denote universal positive constants. We write $A\lesssim B$ or $A=O(B)$ if $|A|\leq CB$ for some  $C>0$. We write $A\sim B$ if $A\lesssim B$ and $B\lesssim A$. We  use $C_{v}$ or $\lesssim_v$ if we want to emphasize that the constant depends on a parameter $v$.  The values of all constants in this paper may vary from line to line.

In this paper, $R$ is reserved for the radius of the scattering data in $q$, i.e. $A(q,\omega)=0$ unless $|q|\leq R$. Unless specified otherwise, we always assume that $t\geq T_R>1$ for some sufficiently large constant $T_R$ depending on $R$ (denoted by $T_R\gg 1$, or $t\gg_R1$). We also assume $0<\eps<1$ is sufficiently small (denoted by $\eps\ll 1$). $T_R$ and $\eps$ are allowed to depend on all other constants, and $\eps$ can also depend on $T_R$. 

We always assume that the latin indices $i,j,l$ take values in $\{1,2,3\}$ and the greek indices $\alpha,\beta$ take values in $\{0,1,2,3\}$. We use subscript to denote partial derivatives, unless specified otherwise. For example, $u_{\alpha\beta}=\partial_\alpha\partial_\beta u$, $q_r=\partial_rq=\sum_i \omega_i\partial_iq$, $A_q=\partial_qA$ and etc. For a fixed integer $k\geq 0$, we  use $\partial^k$ to denote either a specific $k$-th partial derivative, or the collection of all $k$-th partial derivatives.

To prevent confusion, we will only use $\partial_\omega$ to denote the angular derivatives  under the coordinate $(s,q,\omega)$, and will never use it under the coordinate $(t,r,\omega)$. We use $\partial_\omega^c$ to denote $\partial_{\omega_1}^{c_1}\partial_{\omega_2}^{c_2}\partial_{\omega_3}^{c_3}$ for a multiindex $c=(c_1,c_2,c_3)$.

\subsection{Commuting vector fields}

Let $Z$ be any of the following vector fields:
\begin{equation}\label{vf} \partial_\alpha,\ S=t\partial t+r\partial_r,\ \Omega_{ij}=x_j\partial_i-x_i\partial_j,\ \Omega_{0i}=x_i\partial_t+t\partial_i.\end{equation}For any multiindex $I$ with length $|I|$, let $Z^I$ denote the product of $|I|$ such vector fields. Then we have Leibniz's rule
\begin{equation}Z^I(fg)=\sum_{|J|+|K|=|I|}C^I_{JK}Z^JfZ^Kg,\hspace{1cm}\text{where $C_{JK}^I$ are constants.}\end{equation}

The vector fields $Z$ have many good properties. First, we have the  commutation properties.
\begin{equation}[S,\square]=-2\square,\hspace{1cm}[Z,\square]=0\text{ for other $Z$};\end{equation}
\begin{equation}[Z_1,Z_2]=\sum_{|I|=1} C_{Z_1,Z_2,I}Z^I,\hspace{1cm}\text{where $C_{Z_1,Z_2,I}$ are constants};\end{equation}
\begin{equation}[Z,\partial_\alpha]=\sum_\beta C_{Z,\alpha\beta}\partial_\beta,\hspace{1cm}\text{where $C_{Z,\alpha\beta}$ are constants}.\end{equation}

\subsection{Several pointwise bounds}

We have the pointwise estimates for partial derivatives.
\begin{lem}\label{l21}
For any function $\phi$, we have 
\begin{equation}|\partial\phi|\leq C(1+|t-r|)^{-1}\sum_{|I|=1}|Z^I\phi|\end{equation}
and
\begin{equation}|(\partial_t+\partial_r)\phi|+|(\partial_i-\omega_i\partial_r)\phi|\leq C(1+t+r)^{-1}\sum_{|I|=1}|Z^I\phi|.\end{equation}
\end{lem}
Finally, we have the Klainerman-Sobolev inequality.
\begin{prop} For $\phi\in C^\infty(\R^{1+3})$ which vanishes for large $|x|$, we have
\begin{equation}(1+t+|x|)(1+|t-|x||)^{1/2}|\phi(t,x)|\leq C\sum_{|I|\leq 2}\norm{Z^I\phi(t,\cdot)}_{L^2(\R^3)}.\end{equation}
\end{prop}

We also state the Gronwall's inequality.
\prop{Suppose $A,E,r$ are bounded functions from $[a,b]$ to $[0,\infty)$. Suppose that $E$ is increasing. If \fm{A(t)\leq E(t)+\int_a^br(s)A(s)\ ds,\hspace{1cm} \forall t\in[a,b],}then\fm{A(t)\leq E(t)\exp(\int_a^t r(s)\ ds),\hspace{1cm} \forall t\in[a,b].}}\rm
The proofs of these results are standard. See, for example, \cite{lind,sogg,horm} for the proofs.

We also need the following lemma, which can be viewed as the estimates for Taylor's series adapted to $Z$ vector fields.

\lem{\label{l24}Fix $\eps>0$, an integer $k\geq 0$ and a multiindex $I$. Suppose there are two functions $u,v$ on $(t,x)$ such that $|u|+|v|\leq 1$  for all $(t,x)$. Suppose $f\in C^\infty(\R)$ with $f(0)=f^\prime(0)=0$. Then, for all $(t,x)$, we have 
\begin{equation}\label{l242}\begin{aligned}&\hspace{2em}|\partial^kZ^I(f(u+v)-f(u))|\\&\lesssim_{k,I} \sum_{k_1+k_2\leq k,\ |I_1|+|I_2|\leq |I|}p_{k,I}|\partial^{k_1}Z^{I_1}v(t,x)|(|\partial^{k_2}Z^{I_2}v(t,x)|+|\partial^{k_2}Z^{I_2}u(t,x)|).\end{aligned}\end{equation}
where \fm{p_{k,I}(t,x)=1+\max_{k_1+|J|\leq (k+|I|)/2}(|\partial^{k_1}Z^Ju(t,x)|+|\partial^{k_1}Z^Jv(t,x)|)^{k+|I|}.}
}

\begin{proof} By the chain rule and Leibniz's rule, $\partial^kZ^I(f(u))$ can be written as a sum of terms of the form
\fm{f^{(l)}(u) \partial^{k_1}Z^{I_1}u \partial^{k_2}Z^{I_2}u\cdots \partial^{k_l}Z^{I_l}u}where $l\leq k+|I|$, $k_i+|I_i|>0$ for each $i$ and $\sum_ik_i= k$, $\sum_iI_i=I$.  Thus, $\partial^kZ^I(f(u+v)-f(u))$ can be written as a sum of terms of the form
\fm{&\hspace{1em}f^{(l)}(u+v) \partial^{k_1}Z^{I_1}(u+v) \partial^{k_2}Z^{I_2}(u+v)\cdots \partial^{k_l}Z^{I_l}(u+v)-f^{(l)}(u) \partial^{k_1}Z^{I_1}u \partial^{k_2}Z^{I_2}u\cdots \partial^{k_l}Z^{I_l}u
\\&=(f^{(l)}(u+v)-f^{(l)}(u))\partial^{k_1}Z^{I_1}(u+v) \cdots \partial^{k_l}Z^{I_l}(u+v)\\&\hspace{1em}+\sum_{j=1}^l f^{(l)}(u)\partial^{k_1}Z^{I_1}u\cdots \partial^{k_{j-1}}Z^{I_{j-1}}u\cdot\partial^{k_j}Z^{I_j}v\cdot \partial^{k_{j+1}}Z^{I_{j+1}}(u+v) \cdots \partial^{k_l}Z^{I_{l}}(u+v)}
where $k_i+|I_i|>0$ for each $i$ and $\sum_ik_i=k$, $\sum_iI_i=I$.  When $l=0$, we must have $k=|I|=0$, so \eqref{l242} follows from 
\fm{|f(u+v)-f(u)|&\leq \sup_{\beta\in[0,1]}|f^{\prime}(u+\beta v)||v|\leq \sup_{|z|\leq 1}|f^{\prime\prime}(z)|\cdot\sup_{\beta\in[0,1]}|u+\beta v|\cdot|v|\leq C(|u|+|v|)|v|.} Note that now $p_{0,0}=2$. When $l\geq 1$, since  $k_i+|I_i|>(k+|I|)/2>0$ for at most one $i$ and since the product of all other terms of the form $\partial_{k_i}Z^{I_i}(u+v)$ can be controlled by $p_{k,I}$, we have 
\fm{&\hspace{1em}|(f^{(l)}(u+v)-f^{(l)}(u))\partial^{k_1}Z^{I_1}(u+v) \cdots \partial^{k_l}Z^{I_l}(u+v)|\\&\leq \sup_{\beta\in[0,1]}|f^{(l+1)}(u+\beta v)||v\cdot \partial^{k_1}Z^{I_1}(u+v) \cdots \partial^{k_l}Z^{I_l}(u+v)|\\&\leq C_{k,I}p_{k,I}|v|\sum_{k_1\leq k,|J|\leq |I|}(|\partial^{k_1}Z^Ju|+|\partial^{k_1}Z^Jv|).}
When $l=1$, we have
\fm{|f^\prime(u)\partial^kZ^Iv|\leq C|u||\partial^kZ^Iv|.}
When $l\geq 2$, since  $k_i+|I_i|>(k+|I|)/2$ for at most one $i$  and since the product of all other terms of the form $\partial_{k_i}Z^{I_i}(u+v)$ or $\partial_{k_i}Z^{I_i}u$ can be controlled by $p_{k,I}$, we have
\fm{&\hspace{1em}|f^{(l)}(u)\partial^{k_1}Z^{I_1}u\cdots \partial^{k_{j-1}}Z^{I_{j-1}}u \cdot \partial^{k_j}Z^{I_j}v\cdot\partial^{k_{j+1}}Z^{I_{j+1}}(u+v) \cdots \partial^{k_l}Z^{I_{l}}(u+v)|\\&\leq C_{k,I}p_{k,I} \sum_{k_1+k_2\leq k,\ |I_1|+|I_2|\leq |I|}|\partial^{k_1}Z^{I_1}v|(|\partial^{k_2}Z^{I_2}u|+|\partial^{k_2}Z^{I_2}v|).}
\end{proof}

\rm

\subsection{A function space}

Suppose $\eps\ll1$. Let $\D$ be a region in $\R^{1+3}_{t,x}$.  We assume that $\D\subset\{t\geq T_R,\ -R\leq r-t\lesssim  t^{C\eps}\}$ where $T_R\gg 1$. We introduce the following definition based on $\D$, which is useful in Section \ref{sec4.2}.
\defn{\rm For any smooth function $F=F(t,r,\omega)$, we say $F\in S^m=S_\D^m$ for a fixed $m\in\R$ if $|Z^IF|\lesssim_I   t^{m+C_I\eps}$ for any multiindex $I$ and $(t,r,\omega)\in\D$. Here $Z^I$ is a product of $|I|$ vector fields in \eqref{vf}. We also set $\eps^nS^m=\{\eps^nF:\ F\in S^m\}$ for $0<\eps<1$. We  allow $F$ to depend on $\eps$, so $\eps^{n_1}S^m\subset \eps^{n_2}S^m$ if $n_1\geq n_2$.}\rm

For example, we have $t^{m},r^{m}\in S^{m}$,  $\partial^m\omega_i\in S^{-m}$ and $r-t\in S^0$.

We have the following  two lemmas.

\lem{\label{lemons}$S^m$ has the following properties.

{\rm (a)} For any $F_1\in S^{m_1}$ and $F_2\in S^{m_2}$, we have $F_1+ F_2\in S^{\max\{m_1,m_2\}}$ and $F_1F_2\in S^{m_1+m_2}$.

{\rm (b)} For any $F\in S^m$,  we have $ZF\in S^{m}$, $(\partial_t+\partial_r)F\in S^{m-1}$ and $(\partial_i-\omega_i\partial_r)F\in S^{m-1}$.
}
\begin{proof}Note that (a) and $ZF\in S^m$ in (b) are obvious from the definition and the Leibniz's rule. It remains to prove $(\partial_t+\partial_r)F\in S^{m-1}$ and $(\partial_i-\omega_i\partial_r)F\in S^{m-1}$.

Note that \fm{\partial_t+\partial_r&=\frac{\sum_i \omega_i\Omega_{0i}+S}{r+t},\hspace{1cm}\partial_i-\omega_i\partial_r=\frac{\Omega_{0i}+(t-r)\omega_i\partial_t}{t}-\frac{\sum_j \omega_j\omega_i\Omega_{0j}+\omega_iS}{r+t}.
}
Let $f_{-1}$ be any element in $S^{-1}$ and we allow $f_{-1}$ to vary from line to line. Since $r-t,Z^I\omega_i\in S^0$ and $t^{-1},(r+t)^{-1}\in S^{-1}$, by applying (a) of this lemma, we can write $\partial_t+\partial_r$ and $\partial_i-\omega_i\partial_r$ as $\sum_{|J|=1} f_{-1}Z^J$. We claim that for each $I$
\fm{\ Z^I(\partial_t+\partial_r)F&=(\partial_t+\partial_r)Z^IF+\sum_{|J|\leq |I|} f_{-1}Z^JF.}
We can induct on $|I|$. If $|I|=0$, there is nothing to prove. If this equality holds for all $|I|<k$, then for $|I|=k$, by writing $Z^I=ZZ^{I^\prime}$ we have
\fm{&\hspace{1.25em}Z^I(\partial_t+\partial_r)F\\&=Z(\partial_t+\partial_r)Z^{I^\prime}F+\sum_{|J|<k} Z(f_{-1}Z^JF)\\&=(\partial_t+\partial_r)Z^IF+[Z,\partial_t+\partial_r]Z^{I^\prime}F+\sum_{|J|<k} ((Zf_{-1})Z^JF+f_{-1}ZZ^JF)\\&=(\partial_t+\partial_r)Z^IF+\sum_{|K|=1}((Zf_{-1})Z^K+f_{-1}[Z,Z^K])Z^{I^\prime}F+\sum_{|J|<k} (f_{-1}Z^JF+f_{-1}ZZ^JF)\\&=(\partial_t+\partial_r)Z^IF+\sum_{|J|\leq k}f_{-1}Z^JF.} 

Since $F\in S^m$, we have $f_{-1}Z^JF\in S^{m-1}$ for all $|J|\leq |I|$ by (a). Since by Lemma \ref{l21} we have \fm{|(\partial_t+\partial_r)Z^IF|\lesssim\lra{t+r}^{-1}\sum_{|J|\leq |I|+1}|Z^JF|,} we conclude that  for each $I$, in $\D$ we have \fm{|Z^I(\partial_t+\partial_r)F|\lesssim_I t^{m-1+C_I\eps}.}Thus $(\partial_t+\partial_r)F\in S^{m-1}$. Following the same proof, we can also show $(\partial_i-\omega_i\partial_r)F\in S^{m-1}$.
\end{proof}

\lem{\label{lemons2}If $f\in C^\infty(\R)$ with $f(0)=f^\prime(0)=0$, and if $u,v\in \eps^nS^{-m}$ with $n,m\geq 1$, $\eps\ll 1$, and $T_R\gg 1$  in $\D$, then $f(u)-f(v)\in \eps^{2n}S^{-2m}$.}
\begin{proof}
Since $n,m\geq 1$, we have $|u|+|v|\lesssim \eps t^{-1+C\eps}$, so when $\eps\ll 1$ and $t\geq T_R\gg 1$, we have $|u|+|v|\leq 1$. Note that here $\eps$ and $T_R$ do not depend on $I$. Now we can apply Lemma \ref{l24} to $Z^I(f(u)-f(v))$. We have \fm{p_{0,I}(t,x)&= 1+\max_{|J|\leq |I|/2}(|Z^Ju(t,x)|+|Z^Jv(t,x)|)^{|I|}\\&\leq 1+(C_I\eps^n t^{-m+C_I\eps})^{|I|}\\&\leq 1+C_I^{|I|}\eps^{n|I|}t^{-m|I|}t^{C_I|I|\eps}\\&\leq(1+C_I)^{|I|}t^{C_I|I|\eps}.}
The last inequality holds since $n,m\geq 1$. Thus, we have 
\fm{|Z^I(f(u)-f(v))|\lesssim_I t^{C_I\eps}\sum_{|I_1|+|I_2|\leq |I|}|Z^{I_1}v|(|Z^{I_2}u|+|Z^{I_2}v|)\lesssim_I\eps^{2n}t^{-2m+C_I\eps}.}
\end{proof}
\rm

\section{The Derivation of the Asymptotic Equations}\label{sae}

\subsection{The asymptotic equations for \eqref{qwe}}

Let $u$ be a global solution to \eqref{qwe}. Let $q(t,r,\omega)$ be a solution of the eikonal equation \eqref{eikeqn} related to \eqref{qwe}, and let $\mu=q_t-q_r$. Suppose $u$ has the form
\begin{equation}\label{anew}
u(t,x)\approx \frac{\eps}{r}U(s,q,\omega)\end{equation}where $\omega_i=x_i/r$, 
$s=\eps\ln(t)$ and $q=q(t,r,\omega)$. Our goal in this section is to derive the asymptotic equations for $(\mu,U)$.

We make the following assumptions:
\begin{enumerate}
\item Every function is smooth.
\item There is a diffeomorphism between two coordinates $(t,r,\omega)$ and $(s,q,\omega)$, so any function $F$ can be  written as  $F(t,r,\omega)$ and $F(s,q,\omega)$ at the same time.
\item $\eps>0$ is sufficiently small, $t,r>0$ are both sufficiently large with $t\approx r$.
\item All the angular derivatives are negligible. In particular, $\partial_i\approx \omega_i\partial_r$.
\item  $\mu,U\sim 1$ and $\nu\lesssim \eps t^{-1}$, where $\nu:=q_t+q_r$. The same estimates hold if we apply $Z^I$ or $\partial_s^a\partial_q^b\partial_\omega^c$ to the left hand sides.
\end{enumerate}

Here are two useful remarks. First, the solutions $(\mu,U)$ to the reduced system may not exactly satisfy the assumptions listed above. They only satisfy some weaker versions of those assumptions. For example, instead of $\mu\sim 1$, we may only get $ t^{-C\eps}\lesssim|\mu|\lesssim t^{C\eps}$; by solving $q_t-q_r=\mu$,  instead of an exact optical function, i.e. a solution to \eqref{eikeqn}, we may only get an approximate optical function $q$ in the sense that
$\wt{g}^{\alpha\beta}(u)q_\alpha q_\beta=O(t^{-2+C\eps})$.  Such differences are usually negligible, so our assumptions at the beginning  make sense. 

Second,  it may seem strange that we ignore the angular derivatives of $q$ which is $\lesssim t^{-1}$ but keep $\nu\lesssim \eps t^{-1}$. This, however, is reasonable according to the form of \eqref{qwe} and \eqref{eikeqn}. For example, if we expand the eikonal equation, we get \eqref{f33} below. The angular derivatives are either squared or multiplied by $\eps r^{-1}U$, while the major terms in \eqref{f33} are of the order $\eps t^{-1}$. On the other hand, $\nu$ is not negligible since there is a term $\mu\nu$ in the expansion.

Recall that \fm{\square u=r^{-1}((\partial_t-\partial_r)(\partial_t+\partial_r)-r^{-2}\Delta_\omega)ru}where $\Delta_\omega=\sum_{i<j}\Omega_{ij}^2$ is the Laplacian on the sphere $\mathbb{S}^2$. By chain rule we have
\fm{\partial_t=\eps t^{-1}\partial_s+q_t\partial_q,\hspace{1cm}\partial_r=q_r\partial_q.}  By the assumptions, we have \fm{\square u&\approx \eps r^{-1} (\partial_t-\partial_r)(\partial_t+\partial_r)U\approx \eps r^{-1}\mu\partial_q(\eps t^{-1}U_s+\nu U_q)\\&\approx \eps^2 (tr)^{-1}\mu U_{sq}+\eps r^{-1}\mu\nu_qU_q+\eps r^{-1}\mu\nu U_{qq}.}
Since 
\fm{q_t&=\frac{1}{2}(\mu+\nu)\approx\frac{1}{2}\mu,&q_i\approx \omega_iq_r\approx \frac{\omega_i}{2}(\nu-\mu)\approx-\frac{1}{2}\omega_i\mu, \\
q_{tt}&\approx \frac{1}{2}\mu_t\approx \frac{1}{2}\mu_qq_t\approx \frac{1}{4}\mu\mu_q,&q_{it}\approx\frac{1}{2}\mu_i\approx \frac{1}{2}\mu_qq_i\approx-\frac{1}{4}\omega_i\mu\mu_q,\\q_{ij}&\approx-\frac{1}{2}\omega_i\mu_j\approx -\frac{1}{2}\omega_i\mu_qq_j\approx \frac{1}{4}\omega_i\omega_j\mu_q\mu,}
we have \fm{g^{\alpha\beta}q_\alpha q_\beta\approx\frac{1}{4}G(\omega)\mu^2,\hspace{1cm}g^{\alpha\beta}q_{\alpha\beta}\approx \frac{1}{4}G(\omega)\mu\mu_q,}where \[G(\omega)=g^{\alpha\beta}\widehat{\omega}_\alpha\widehat{\omega}_\beta,\hspace{1cm }g^{\alpha\beta}=\frac{d}{du}\wt{g}^{\alpha\beta}(u)|_{u=0},\ \widehat{\omega}=(-1,\omega)\in\R\times\mathbb{S}^2.\]
And since
\fm{U_{tt}\approx U_{qq}q_{tt}+U_q q_t^2,\hspace{1cm}U_{it}\approx U_{qq}q_iq_t+U_{q}q_{it},\hspace{1cm}U_{ij}\approx U_{qq}q_iq_j+U_{q}q_{ij},} we have from \eqref{qwe}
\begin{equation}\label{f2}\begin{aligned}0&=\wt{g}^{\alpha\beta}(u)\partial_\alpha\partial_\beta u\approx \square u+ g^{\alpha\beta}u\partial_\alpha\partial_\beta u\\&\approx \eps^2 (tr)^{-1}\mu U_{sq}+\eps r^{-1}\mu\nu_qU_q+\eps r^{-1}\mu\nu U_{qq}+\eps^2r^{-2}g^{\alpha\beta}U(U_qq_{\alpha\beta}+U_{qq}q_\alpha q_\beta)\\&\approx \eps^2 (tr)^{-1}\mu U_{sq}+\eps r^{-1}\mu\nu_qU_q+\eps r^{-1}\mu\nu U_{qq}+\frac{1}{4}G(\omega)\eps^2r^{-2}(\mu\mu_qUU_q+\mu^2UU_{qq}).\end{aligned}\end{equation}

By the eikonal equation, we have 
\begin{equation}\label{f33}\begin{aligned}0&=\wt{g}^{\alpha\beta}(u)q_\alpha q_\beta\approx q_t^2-\sum_i q_i^2+\eps r^{-1}g^{\alpha\beta}U q_\alpha q_\beta\approx\mu\nu+\frac{1}{4}\eps r^{-1}G(\omega)\mu^2U,\end{aligned}\end{equation}
so we conclude that \fm{\nu\approx-\frac{1}{4}\eps r^{-1}G(\omega)\mu U,\hspace{1cm} \nu_q\approx-\frac{1}{4}\eps r^{-1}G(\omega)(\mu_qU+\mu U_q).}
Plug everything back in \eqref{f2}. We thus have 
\fm{0&\approx \eps^2 (tr)^{-1}\mu U_{sq}-\frac{1}{4}\eps^2 r^{-2}G(\omega)(\mu_qU+\mu U_q)\mu U_q\\&\hspace{1em}-\frac{1}{4}\eps^2 r^{-2}G(\omega)\mu^2 UU_{qq}+\frac{1}{4}G(\omega)\eps^2r^{-2}(\mu\mu_qUU_q+\mu^2UU_{qq})\\&=\eps^2 (tr)^{-1}\mu U_{sq}-\frac{1}{4}\eps^2 r^{-2}G(\omega)\mu^2 U_q^2.}
Assuming that $t=r$, we get the first asymptotic equation \fm{U_{sq}=\frac{1}{4}G(\omega)\mu U_q^2.}

Meanwhile, note that from $(\partial_t-\partial_r)\nu=(\partial_t+\partial_r)\mu$, we have \fm{\nu_q\mu\approx\nu_q\mu+\eps t^{-1}\nu_s=\mu_q\nu+\eps t^{-1}\mu_s} and thus \fm{\mu_s&\approx t\eps^{-1}(\nu_q\mu-\mu_q\nu)\approx t\eps^{-1}(-\frac{1}{4}\eps r^{-1}G(\omega)(\mu_qU+\mu U_q)\mu+\frac{1}{4}\eps r^{-1}G(\omega)\mu U\mu_q)\\&\approx -\frac{t}{4r}G(\omega)\mu^2 U_q.}
Again, assuming that $t=r$, we get the second asymptotic equation \fm{\mu_{s}=-\frac{1}{4}G(\omega)\mu^2 U_q.}
In conclusion, our system of asymptotic equations is 
\begin{equation}\label{f2asypde3}
\left\{
\begin{array}{l}\displaystyle\partial_s\mu=-\frac{1}{4}G(\omega)\mu^2U_q,\\[.2cm]\displaystyle\partial_sU_q=\frac{1}{4}G(\omega)\mu U_q^2.\end{array}\right.\end{equation}

Now we can solve \eqref{f2asypde3} if we assign some reasonable initial data. Since we expect $q\approx r-t$ and since $(\partial_t-\partial_r)(r-t)=-2$, we choose $\mu|_{s=0}=-2$. Since there is no restriction on $U_q$,  we  choose arbitrarily $U|_{s=0}=A(q,\omega)\in C_c^\infty(\R\times\mathbb{S}^2)$.  Note that the two asymptotic equations imply that $(\mu U_q)_s=\mu_s U_q+\mu U_{sq}=0$, so we have $\mu U_q=-2A$. Thus, we get two ODEs
\begin{equation}\label{sasf}
\left\{
\begin{array}{l}\displaystyle\partial_s\mu=\frac{1}{2}G(\omega)A(q,\omega)\mu,\\[.2cm]\displaystyle\partial_sU_q=-\frac{1}{2}G(\omega)A(q,\omega)U_q.\end{array}\right.\end{equation}Then we can solve $\mu$ and $U_q$ easily. See \eqref{mu} and \eqref{ueqn} in Section \ref{sas} for the explicit formulas.

We can compare our system \eqref{f2asypde3}  with H\"{o}rmander's system 
\fm{2\partial_s\partial_q U=G(\omega)U\partial_q^2U.}
One advantage of our system is that it reduces to a linear first-order ODE system \eqref{sasf} after we choose the appropriate initial data. The solution to \eqref{f2asypde3} is thus of the simpler form than the solution to H\"{o}rmander's system.

\subsection{Asymptotic equations for general case}
Though \eqref{sasf} is already enough for this paper, let us also do the computations in a more general case. Instead of the fully nonlinear wave equation \eqref{nlw}, we consider the following quasilinear wave equation
\begin{equation}\label{sasqua}
\wt{g}^{\alpha\beta}(u,\partial u)\partial_\alpha\partial_\beta u=f(u,\partial u).
\end{equation}
Assume that we have Taylor expansions
\fm{\wt{g}^{\alpha\beta}(u,\partial u)&=-m^{\alpha\beta}+g^{\alpha\beta}u+g^{\alpha\beta\lambda}\partial_\lambda u+O(|u|^2+|\partial u|^2),\\
f(u,\partial u)&=f_0u^2+f^{\alpha}u\partial_\alpha u+f^{\alpha\beta}\partial_\alpha u\partial_\beta u+O(|u|^3+|\partial u|^3).}
Here $m^{\alpha\beta},g^{*},f_0,f^{*}$ are all real constants.

We still make the ansatz \eqref{anew} with the same $s,\omega,r$, assuming that $q$ is now the solution to the eikonal equation
\begin{equation}\label{saseik}
\wt{g}^{\alpha\beta}(u,\partial u)\partial_\alpha q\partial_\beta q=0.
\end{equation}
Again, we take $\mu=q_t-q_r$ and $\nu=q_t+q_r$.

From \eqref{anew} and \eqref{sasqua}, we have
\fm{0&\approx\eps^2 (tr)^{-1}\mu U_{sq}+\eps r^{-1}\mu\nu_qU_q+\eps r^{-1}\mu\nu U_{qq}+\eps^2 r^{-2}g^{\alpha\beta}UU_{\alpha\beta}+\eps^2 r^{-2}g^{\alpha\beta\lambda}U_\lambda U_{\alpha\beta}\\&\hspace{1em}-\eps^2 r^{-2}f_0U^2-\eps^2 r^{-2}f^\alpha UU_\alpha-\eps^2 r^{-2}f^{\alpha\beta}U_\alpha U_\beta\\&\approx\eps^2 (tr)^{-1}\mu U_{sq}+\eps r^{-1}\mu\nu_qU_q+\eps r^{-1}\mu\nu U_{qq}+\eps^2 r^{-2}g^{\alpha\beta}U(U_{qq}q_\alpha q_\beta+U_{q}q_{\alpha\beta})\\&\hspace{1em}+\eps^2 r^{-2}g^{\alpha\beta\lambda}U_qq_\lambda (U_{qq}q_\alpha q_\beta+U_{q}q_{\alpha\beta})-\eps^2 r^{-2}f_0U^2-\eps^2 r^{-2}f^\alpha UU_qq_\alpha-\eps^2 r^{-2}f^{\alpha\beta}U_q^2q_\alpha q_\beta.}
By \eqref{saseik} and computation in the previous subsection, we have
\fm{0&\approx\mu\nu+\eps r^{-1}g^{\alpha\beta} Uq_\alpha q_\beta+\eps r^{-1}g^{\alpha\beta\lambda} U_qq_\alpha q_\beta q_\lambda\approx\mu\nu+\frac{\eps }{4r}G_2(\omega)\mu^2U-\frac{\eps}{8 r}G_3(\omega) \mu^3U_q,}
where $G_2(\omega)=g^{\alpha\beta}\widehat{\omega}_\alpha\widehat{\omega}_\beta$ and $G_3(\omega)=g^{\alpha\beta\lambda}\widehat{\omega}_\alpha\widehat{\omega}_\beta\widehat{\omega}_\lambda$ for $\widehat{\omega}=(-1,\omega)$. Similarly we can define $F_1(\omega)=f^\alpha\widehat{\omega}_\alpha$ and $F_2(\omega)=f^{\alpha\beta}\widehat{\omega}_\alpha\widehat{\omega}_\beta$. Thus,
\fm{\nu&\approx-\frac{\eps }{4r}G_2(\omega)\mu U+\frac{\eps}{8 r}G_3(\omega) \mu^2U_q,\\
\nu_q&\approx-\frac{\eps }{4r}G_2(\omega)(\mu_q U+\mu U_q)+\frac{\eps}{8 r}G_3(\omega) (\mu^2U_{qq}+2\mu\mu_qU_q).}
By letting $t\approx r$, we have
\fm{0&\approx\mu U_{sq}+\mu(-\frac{1}{4}G_2(\omega)(\mu_q U+\mu U_q)+\frac{1}{8}G_3(\omega) (\mu^2U_{qq}+2\mu\mu_qU_q))U_q\\&\hspace{1em}+g^{\alpha\beta}UU_{q}q_{\alpha\beta}+g^{\alpha\beta\lambda}U_qq_\lambda U_{q}q_{\alpha\beta}-f_0U^2-f^\alpha UU_qq_\alpha-f^{\alpha\beta}U_q^2q_\alpha q_\beta\\&\approx\mu U_{sq}-\frac{1}{4}G_2(\omega)\mu^2U_q^2+\frac{1}{8}G_3(\omega)\mu^3U_{qq}U_q+\frac{1}{8}G_3(\omega)\mu^2\mu_qU_q^2\\&\hspace{1em}-f_0U^2+\frac{1}{2}F_1(\omega) \mu UU_q-\frac{1}{4}F_2(\omega)\mu^2U_q^2.}
Besides, since  $\nu_q\mu+\eps t^{-1}\nu_s=\mu_q\nu+\eps t^{-1}\mu_s$, we have 
\fm{\mu_s&\approx t\eps^{-1}(\nu_q\mu-\mu_q\nu)\\&\approx \mu(-\frac{1}{4}G_2(\omega)(\mu_q U+\mu U_q)+\frac{1}{8 }G_3(\omega) (\mu^2U_{qq}+2\mu\mu_qU_q))\\&\hspace{1em}-\mu_q(-\frac{1}{4}G_2(\omega)\mu U+\frac{1}{8}G_3(\omega) \mu^2U_q)\\&\approx-\frac{1}{4}G_2(\omega)\mu^2 U_q+\frac{1}{8 }G_3(\omega) \mu^3U_{qq}+\frac{1}{8}G_3(\omega) \mu^2\mu_qU_q.}

Note that now
\fm{(\mu U_q)_s\approx f_0U^2-\frac{1}{2}F_1(\omega) \mu UU_q+\frac{1}{4}F_2(\omega)\mu^2U_q^2.}

\begin{defn}\rm We define the following reduced system of \eqref{sasqua} for $(\mu,U)(s,q,\omega)$
\begin{equation}\begin{aligned}
\left\{
\begin{array}{l}
(\mu U_q)_s= f_0U^2-\frac{1}{2}F_1(\omega) \mu UU_q+\frac{1}{4}F_2(\omega)\mu^2U_q^2\\[.2cm]
\mu_s=-\frac{1}{4}G_2(\omega)\mu^2 U_q+\frac{1}{8 }G_3(\omega) \mu^3U_{qq}+\frac{1}{8}G_3(\omega) \mu^2\mu_qU_q
\end{array}\right..
\end{aligned}\end{equation}
\end{defn}
\rmk{\rm
It is unclear to the author whether the lifespan of this new reduced system is the same as that of \eqref{asypde11}. In the special case $f\equiv 0$ and $\wt{g}^{\alpha\beta}(u,\partial u)=\wt{g}^{\alpha\beta}(\partial u)$, the answer is yes. Now $G_2(\omega)\equiv 0$ and $(\mu U_q)_s\equiv 0$, and our new reduced system  admits a finite-time blowup in~$s$, unless the null condition holds, i.e. $G_3(\omega)\equiv 0$. In fact, since $(\mu U_q)_s=0$, with the same choice of initial data $\mu|_{s=0}=-2$ and $U_q|_{s=0}=A\in C_c^\infty(\R\times\mathbb{S}^2)$, we have $\mu U_q=-2A$ for all~$s$. Thus,
\fm{\mu_s=\frac{1}{8}G_3(\omega)\mu^2(\mu U_q)_q=-\frac{1}{4}G_3(\omega)A_q(q,\omega)\mu^2,\hspace{1cm}\mu(0,q,\omega)=-2,}
whose solution is  \fm{\mu(s,q,\omega)=(\frac{1}{4}G_3(\omega)A_q(q,\omega)s-\frac{1}{2})^{-1}.}
If $G_3(\omega)\not\equiv 0$ and $A(q,\omega)\not\equiv 0$, we can choose $(q,\omega)$ such that $G_3(\omega)A_q(q,\omega)>0$. We are able to do this because $A(q,\omega)$ has a compact support. This would lead to a blowup at $s=2/(G_3(\omega)A_q(q,\omega))$. Such a blowup can be related to the blowup of H\"{o}rmander's approximate equation \eqref{asypde11}, which is now a Burgers' equation. We refer to Lemma 6.5.4 in \cite{horm2}. This result implies that our new reduced system may work in a  more general case than \eqref{qwe}. \rm

\section{The Asymptotic Profile and the Approximate Solution}\label{sas}

Our main goal in this section is to construct an approximate solution $u_{app}$ to \eqref{qwe}. Fix a scattering data $A(q,\omega)\in C_c^\infty(\R\times\mathbb{S}^2)$ with ${\rm supp}(A)\subset[-R,R]\times\mathbb{S}^2$ for some $R\geq 1$. Fix a sufficiently small $\eps>0$ and a sufficiently large $T_R>0$, both depending on $A(q,\omega)$. Let $(\mu,U)(s,q,\omega)$ be the solution to \eqref{f2asypde3} with  $(\mu,U_q)|_{s=0}=(-2,A)$ and $\lim_{q\to-\infty}U(s,q,\omega)=0$. Let $q(t,r,\omega)$ be the solution to the ODE \fm{(\partial_t-\partial_r)q(t,r,\omega)=\mu(\eps\ln(t)-\delta,q(t,r,\omega),\omega),\hspace{1cm}q(t,0,\omega)=-t}and set\fm{U(t,r,\omega)=U(\eps\ln(t)-\delta,q(t,r,\omega),\omega).}Here $\delta>0$ is a sufficiently small constant depending only on the scattering data.
Note that near the light cone $\{t=r+R\}$, $\eps r^{-1}U(t,r,\omega)$  and $q(t,r,\omega)$ are the approximate solution to \eqref{qwe} and the approximate optical function, respectively,  in the sense that for all $(t,r,\omega)$ with $t\geq T_R$ and $|q(t,r,\omega)|\leq R$, we have
\fm{\wt{g}^{\alpha\beta}(\eps r^{-1}U)\partial_\alpha\partial_\beta(\eps r^{-1}U)=O(\eps t^{-3+C\eps}),}
\fm{\wt{g}^{\alpha\beta}(\eps r^{-1}U)q_\alpha q_\beta=O( t^{-2+C\eps}).}
For all $t\geq 0$ and $x\in\R^3$, we set
\fm{u_{app}(t,x)=\eps r^{-1}\eta(t)\psi(r/t)U(\eps\ln(t)-\delta,q(t,r,\omega),\omega).}
  Here $\psi\equiv 1$ when $|r-t|\leq t/4$ and $\psi\equiv 0$ when $|r-t|\geq t/2$, which is used to localize $\eps r^{-1}U$ near the light cone $\{r=t\}$; $\eta$ is a cutoff function such that $\eta\equiv 0$ when $t\leq 2R$. The  definitions of $\psi$ and $\eta$ will be given later.

Our main proposition in this section is the following.
\prop{\label{mainprop4} Fix a scattering data $A(q,\omega)\in C_c^\infty(\R\times\mathbb{S}^2)$ with ${\rm supp}(A)\subset[-R,R]\times\mathbb{S}^2$ for some $R\geq 1$. Fix a sufficiently small $\eps>0$ depending on $A(q,\omega)$. Let $u_{app}$ be the function defined as above.  Then, for all $(t,x)$ with $t\geq T_R$, we have \fm{|\partial u_{app}(t,x)|\lesssim \eps(1+t)^{-1}.} Moreover, for all multiindices $I$ and for all $(t,x)$ with $t\geq 0$, we have
\fm{|Z^I u_{app}(t,x)|\lesssim_I \eps(1+t)^{-1+C_I\eps},}
\fm{|Z^I(\wt{g}^{\alpha\beta}(u_{app})\partial_\alpha\partial_\beta u_{app})(t,x)|\lesssim_I \eps(1+t)^{-3+C_I\eps}.}}
\rmk{\rm If we have $0<\delta<1$, then all the constants involved in this section are uniform in $\delta$. Thus, it would not impact any result in this section if we do not choose the value of $\delta$ until the proof of the Poincar$\acute{\rm e}$'s lemma in the next section.}
\rm

This proposition is proved in three steps. First, in Section \ref{sec4.1}, we construct $q(t,r,\omega)$ and $U(t,r,\omega)$ for all $(t,x)$ with $t> 0$, by solving the reduced system \eqref{f2asypde3} and $q_t-q_r=\mu$ explicitly. Next, in Section \ref{sec4.2}, we prove that $\eps r^{-1}U(t,r,\omega)$ is an approximate solution to \eqref{qwe} near the light cone $\{t=r+R\}$ when $t$ is sufficiently large. To achieve this goal we prove several estimates for $q$ and $U$ when $|q(t,r,\omega)|\leq R$. Finally, in Section \ref{sec4.3}, we define $u_{app}$ and prove the pointwise bounds for large $t$. To define $u_{app}$, we use  cutoff functions to restrict $\eps r^{-1}U$ in a conical neighborhood of $\{t=r\}$ and remove the singularities at $|x|=0$ or $t=0$.

\subsection{Construction of $q$ and $U$}\label{sec4.1}
Fix a sufficiently small $\eps>0$. Fix a scattering data $A(q,\omega)\in C^\infty_c(\R\times\mathbb{S}^2)$ with $A(q,\omega)= 0$ for $|q|> R\geq 1$. Also fix  $0<\delta<1$ depending on $A(q,\omega)$ but not on $\eps$. Its value will be chosen in Section \ref{sep}.

Suppose the Taylor expansion of $\wt{g}^{\alpha\beta}$ at $0$ is
\fm{\wt{g}^{\alpha\beta}(u)=-m^{\alpha\beta}+\gamma^{\alpha\beta}(u)=-m^{\alpha\beta}+g^{\alpha\beta}u+O(|u|^2),\hspace{1cm}\text{with }(m^{
\alpha\beta})=\text{diag}(-1,1,1,1).} 
We define $q(t,r,\omega)$ by solving \begin{equation}\label{qeqn}
(\partial_t-\partial_r)q(t,r,\omega)=\mu(\eps\ln(t)-\delta,q(t,r,\omega),\omega),\hspace{1cm}q(t,0,\omega)=-t,\end{equation}
where
\begin{equation}\label{mu}
\mu(s,q,\omega):=-2\exp(\frac{1}{2}G(\omega)A(q,\omega)s),
\end{equation}
where
\fm{G(\omega):=g^{\alpha\beta}\widehat{\omega}_\alpha\widehat{\omega}_\beta,\hspace{2em} \widehat{\omega}:=(-1,\omega)\in\R\times\mathbb{S}^2.} 
Note that \eqref{qeqn} has a solution $q(t,r,\omega)$  for all $t> 0$. In fact, if we apply method of characteristics, for $z(\tau)=q(\tau,r+t-\tau,\omega)$ and $s(\tau)=\ln(\tau)$ we have an autonomous system of ODEs
\fm{\left\{\begin{array}{l}\dot{z}(\tau)=\mu( \eps s(\tau)-\delta,z(\tau),\omega)\\\dot{s}(\tau)=\exp(-s(\tau))\end{array}\right.}
with initial data $(z,s)(r+t)=(-r-t,\ln(r+t))$. Note that whenever $|z(\tau)|\geq R$, we have $\dot{z}(\tau)=-2$ because of the support of $A(q,\omega)$. Thus, $|z(\tau)|$ cannot blow up when $\tau>0$. Neither can $|s(\tau)|$ since $s(\tau)=\ln(\tau)$. We are thus able to solve this system of ODEs for all $\tau> 0$ by Picard's theorem.

We have \begin{equation}\label{qeqn2}
q(t,r,\omega)=-(r+t)-\int_{t}^{r+t}\mu(\eps\ln(\tau)-\delta,q(\tau,r+t-\tau,\omega),\omega)\ d\tau.
\end{equation}
Note that if $G(\omega)\equiv 0$, we have $\mu\equiv -2$ and thus $q=r-t$, which concides with the choice of~$q$ in H\"{o}rmander's setting.

We also define  $U(s,q,\omega)$ by solving the following equation
\begin{equation}\label{ueqn}
(\partial_qU)(s,q,\omega)=A(q,\omega)\exp(-\frac{1}{2}G(\omega)A(q,\omega)s),\hspace{1cm}\lim_{q\to-\infty} U(s,q,\omega)=0.
\end{equation}
The equation \eqref{ueqn} has a solution $U(s,q,\omega)$ for all $s$, which comes from taking the following integral:\begin{equation}\label{ueqn2}
U(s,q,\omega)=\int_{-\infty}^{q}A(p,\omega)\exp(-\frac{1}{2}G(\omega)A(p,\omega)s)\ dp.
\end{equation}
It is clear that $U(s,q,\omega)=0$ unless $q\geq -R$ and $U(s,q,\omega)=U(s,R,\omega)$ for $q\geq R$. Also note that $U$ and all its derivatives are $\leq Ce^{Cs}$. Here $C$ is uniform for all $(s,q,\omega)\in\R\times\R\times\mathbb{S}^2$.

From now on, we  use $U$ to denote the function on $(t,r,\omega)$: \begin{equation}\label{Udef} U=U(t,r,\omega)=U(\eps\ln(t)-\delta,q(t,r,\omega),\omega).\end{equation}
Such a $U$ is the asymptotic profile used in this paper. Note that 
\fm{(\partial_t-\partial_r)U=\mu U_q+\eps t^{-1}U_s=-2A+O(\eps t^{-1+C\eps}).}
This explains the meaning of $A(q,\omega)$ in our construction.

\subsection{Estimates for $q$ and $U$}\label{sec4.2}
Define \begin{equation}\label{domain}
\D:=\{(t,r,\omega):\ t\geq T_R,\ |q(t,r,\omega)|\leq R\},
\end{equation}for some constant $T_R\geq 1$ to be chosen. Here we always assume that $T_R$ is sufficiently large and depends only on $A(q,\omega)$. Our main goal now is to prove that $\eps r^{-1}U(t,r,\omega)\in \eps S^{-1}$ and $\wt{g}^{\alpha\beta}(\eps r^{-1}U)\partial_\alpha\partial_\beta(\eps r^{-1}U)\in \eps S^{-3}$, where $\eps S^{-1}$ and $\eps S^{-3}$ are defined in Section \ref{spre}. In other word, $\eps r^{-1}U$ has some good pointwise bounds and is an approximate solution to \eqref{qwe} in $\D$.

We start with a  more precise description of  the region $\D$. From Lemma \ref{lq1}, we can see that $\D$  is contained in  a conical neighborhood of the light cone $\{t=r\}$ when $t\gg_R1$.

\lem{\label{lq1}For all $(t,r,\omega)$ with $t\geq T_R$, there exist $0< t_0<t_1=(t+r+R)/2$ such that
\begin{equation}|q(\tau,r+t-\tau,\omega)|\leq R\ \Longleftrightarrow\ \tau\in[t_0,t_1];\end{equation}
\begin{equation}q(\tau,r+t-\tau,\omega)=R+2(t_0-\tau),\hspace{1cm} \forall\tau\leq t_0;\end{equation}
\begin{equation}q(\tau,r+t-\tau,\omega)=r+t-2\tau,\hspace{1cm} \forall\tau\geq t_1.\end{equation}
We also have 
\begin{equation}
t_1-t_0=O((t+r)^{C\eps}),
\end{equation}\begin{equation}
q(t,r,\omega)-(r-t)=O((t+r)^{C\eps}).
\end{equation}
When $r\leq t-R$, we have $q=r-t$.

In addition, for  $(t,r,\omega)\in\D$, we have $t_0\leq t\leq t_1$ and \begin{equation}-R\leq r-t\lesssim t^{C\eps}\end{equation}
which implies that $t\sim r$ in $\D$. 
\begin{proof}
Note that $\mu\equiv-2$ for $|q|\geq R$ and $-2e^{Cs}\leq\mu\leq -2e^{-Cs}$ for all $(s,q,\omega)$. Then the existence of $t_0,t_1$ and the estimates related to $\tau$ directly follow from \eqref{qeqn2}. We also have $q=r-t$ if $r\leq t-R$. Now we can assume $r\geq t-R$ i.e. $t\leq t_1$.
We have
\fm{|q(t,r,\omega)-(r-t)|&=|\int_t^{r+t}(\mu(\eps\ln(\tau)-\delta,q(\tau,r+t-\tau,\omega),\omega)+2)\ d\tau|\\&\leq \int_{[t,r+t]\cap[t_0,t_1]}|\mu(\eps\ln(\tau)-\delta,q(\tau,r+t-\tau,\omega),\omega)|+2\ d\tau\\&= \int_{[t,r+t]\cap[t_0,t_1]}-\mu(\eps\ln(\tau)-\delta,q(\tau,r+t-\tau,\omega),\omega)+2\ d\tau\\&\leq2R+2(t_1-\max\{t_0,t\}).}
Moreover, we have\fm{-2R&=\int_{t_0}^{t_1}\mu(\tau,r+t-\tau,\omega)\ d\tau\leq -2e^{-C\delta}\int_{t_0}^{t_1}\tau^{-C\eps}\ d\tau\leq -2e^{-C}(t_1-t_0)t_1^{-C\eps}. }
It follows that 
\fm{t_1-\max\{t,t_0\}\leq t_1-t_0\lesssim t_1^{C\eps}\lesssim_R (t+r)^{C\eps}}
and thus
\fm{|q(t,r,\omega)-(r-t)|\leq 2R+C_R(t+r)^{C\eps}\lesssim_R(t+r)^{C\eps}.}
Here we use $t+r+R\leq 2(t+r)$ if $t\geq T_R\geq R$.

When $(t,r,\omega)\in\D$, it is clear that $t_0\leq t\leq  t_1$. Thus,
\fm{\frac{r-t+R}{2}=t_1-t\lesssim_R(r+t)^{C\eps}} By choosing $T_R$ in \eqref{domain} sufficiently large (e.g. $T_R\geq 2+2R$) and $\eps$ sufficiently small (e.g. $C\eps\leq 1/4$), we have \fm{\frac{r}{t}-1\leq \frac{r-t+R}{t}\lesssim_R t^{-1+C\eps}(1+\frac{r}{t})^{C\eps}\lesssim_R(1+\frac{r}{t})^{1/2}.}This forces $r/t\leq C^\prime_R$ for some constant $C^\prime_R$, which implies that\fm{|q-(r-t)|\lesssim_R (t+C^\prime_Rt)^{C\eps}\lesssim_R t^{C\eps}}and\fm{|r-t|\leq |q|+O(t^{C\eps})\leq R+Ct^{C\eps}\lesssim t^{C\eps}.}
Finally, note that $r-t<-R$ implies $q<-R$. We are done.
\end{proof}

\rm
We now move on to estimates for $\partial q$.  In Lemma \ref{lq2}, we give the pointwise bounds for $\nu=q_t+q_r$ and $\lambda_i=q_i-\omega_i q_r$. In Lemma \ref{lq4}, we find the first terms in the asymptotic expansion of $\nu$ and $\nu_q$ in $\D$.

\lem{\label{lq2}For $t\geq T_R$,
\begin{equation}\label{lq21}\nu(t,r,\omega):=(\partial_t+\partial_r)q=O(\eps (t+r)^{-1+C\eps}).\end{equation}
\begin{equation}\label{lq22}\lambda_i(t,r,\omega):=(\partial_i-\omega_i\partial_r)q=O((t+r)^{-1+C\eps}).\end{equation}}
\begin{proof}
Fix $(t,r,\omega)$. We have
\begin{equation}\label{lq23}(\partial_t-\partial_r)\nu=(\partial_t+\partial_r)\mu=(\partial_q\mu)\nu+\frac{\eps}{2t}G(\omega)A\mu.\end{equation}
By Lemma \ref{lq1}, for all $t>T_R$, we have
\fm{\int_t^{r+t}|\partial_q\mu|\ d\tau&= \int_{[t_0,t_1]\cap[t,r+t]}\frac{1}{2}\abs{G(\omega)\partial_qA}\cdot|\eps\ln(\tau)-\delta|\cdot|\mu|\ d\tau\\&\lesssim (\eps\ln(t+r)+1)\int_{[t_0,t_1]\cap[t,r+t]}|\mu|\ d\tau\\&\lesssim(\eps\ln(t+r)+1)|q(t_1,r+t-t_1,\omega)-q(t_0,r+t-t_0,\omega)|\\&\lesssim \eps\ln(t+r)+1.}
Here the integral is taken along the characteristic $(\tau,r+t-\tau,\omega)$ for $\tau\geq T_R$,  as in \eqref{qeqn2}. Similarly,  we have
\fm{\int_t^{r+t}|G(\omega)A\mu\frac{\eps}{2\tau}|\ d\tau&\lesssim\frac{\eps}{t_0}\int_{[t_0,t_1]\cap[t,r+t]}|\mu|\ d\tau\lesssim\eps (t+r)^{-1}.}
Here we use the fact that for $\eps\ll1$ and $t\geq T_R\gg_R1$, we have \fm{t_0\geq t_1-|t_0-t_1|\geq \frac{t+r+R}{2}-C(t+r)^{C\eps}\geq \frac{r+t}{4}.}Now, we integrate \eqref{lq23} along the characteristic   and  then apply Gronwall's inequality. Note that the initial value of $(\partial_t+\partial_r)q$ is $0$ as $q=r-t$ for $r\leq t-R$, by Lemma \ref{lq1}. So we conclude \eqref{lq21}. The proof for \eqref{lq22} is similar. We have
\begin{equation}\label{lq24}\begin{aligned}
(\partial_t-\partial_r)\lambda_i&=(\partial_i-\omega_i\partial_r)\mu+r^{-1}\lambda_i\\&=(\mu_q+r^{-1})\lambda_i+\frac{1}{2}(\eps\ln(t)-\delta)\sum_l\partial_{\omega_l}(GA)\cdot\frac{\delta_{il}-\omega_i\omega_l}{r}\mu\\&=(\mu_q+r^{-1})\lambda_i+O( r^{-1}|\eps\ln(t)-\delta|\cdot|\mu|)\chi_{|q|\leq R}.
\end{aligned}\end{equation}
Here $\chi_{|q|\leq R}=1$ if $|q(t,r,\omega)|\leq R$ and $\chi_{|q|\leq R}=0$ if $|q(t,r,\omega)|>R$. This term exists in \eqref{lq24} since $A\equiv 0$ if $|q|>R$.
Note that $\lambda_i\equiv 0$ when $r<t-R$ and that for $0<t-R\leq r$, we have\fm{0\leq \int_t^{t_1} (r+t-\tau)^{-1}\ d\tau=\ln\frac{2r}{r+t-R}\leq \ln 2}
and
\fm{\int_t^{t_1} (r+t-\tau)^{-1}|\eps\ln(\tau)-\delta||\mu|\chi_{|q|\leq R}\ d\tau&\leq\int_{[t,t_1]\cap[t_0,t_1]} (r+t-\tau)^{-1}(\eps\ln(\tau)+1)\cdot|\mu|\ d\tau\\&\leq  \frac{\eps\ln(t+r)+1}{r+t-t_1}\cdot 2R\lesssim(t+r)^{-1+C\eps}.} Apply Gronwall's inequality again  and we are done.
\end{proof}

\rmk{\rm Since $|\mu|\geq 2ct^{-C\eps}$ for some small constant $c>0$, we conclude that $q_t,q_r\neq 0$ for all $t\geq T_R$ if $\eps$ is small enough. In particular, for $\eps\ll 1$ and $t\gg_R1$, 
\fm{q_r=\frac{-\mu+\nu}{2}\geq ct^{-C\eps}-C\eps(t+r)^{-1+C\eps}\geq \frac{c}{2}t^{-C\eps}.}
So for each fixed $t\geq T_R$, $r\mapsto q(t,r,\omega)$ is strictly increasing and continuous for each fixed~$t$, and $\lim_{r\to\infty}q(t,r,\omega)=\infty$.
This implies that for each $t$ and $q_0\geq -t$, there exists a unique $r$ such that $q(t,r,\omega)=q_0$. So $(t,r,\omega)\mapsto (\eps\ln(t)-\delta,q(t,r,\omega),\omega)$ has an inverse map $(s,q,\omega)\mapsto (e^{(s+\delta)/\eps},r(s,q,\omega),\omega)$. By inverse function theorem, the map $(t,r,\omega)\mapsto(s,q,\omega)$ is a diffeomorphism.

From now on, any function $V$ can be written as both $V(t,r,\omega)$ and $V(s,q,\omega)$ at the same time. Thus, for any function $V$ on $(t,r,\omega)$, we can define $\partial_s^a\partial_q^b\partial_{\omega}^{c}V$ using the chain rule and Leibniz's rule. Note that in this paper, $\partial_\omega$ will only be used under the coordinate $(s,q,\omega)$ and will never be used under the coordinate $(t,r,\omega)$.}

\lem{\label{lq4}For $(t,r,\omega)\in\D$,
\begin{equation}\label{lq41}\nu+\frac{\eps G(\omega)}{4 t}\mu U=O(\eps t^{-2+C\eps}),\end{equation}
\begin{equation}\label{lq42}\nu_q+\frac{\eps G(\omega)}{4t} \mu_q U+\frac{\eps G(\omega)}{4t}\mu U_q=O(\eps t^{-2+C\eps}).\end{equation} }
\begin{proof}We have
\begin{equation}\label{lq43}\begin{aligned}(\partial_t-\partial_r)(\nu+\frac{\eps G(\omega)}{4  t}\mu U)&=(\partial_t+\partial_r)\mu+\frac{\eps G(\omega)}{4  t}(\partial_t-\partial_r)(\mu U)-\frac{\eps G(\omega)}{4 t^{2}}\mu U\\&=\mu_q\nu+\mu_s\frac{\eps}{t}+\frac{\eps G(\omega)}{4  t}(\partial_q(\mu U)\mu+\partial_s(\mu U)\frac{\eps}{t})-\frac{\eps G(\omega)}{4 t^{2}}\mu U\\&=\mu_q(\nu+\frac{\eps G(\omega)}{4  t}\mu U)+\frac{\eps^2G(\omega)}{4 t^2}(U_s+\frac{1}{2}G(\omega)AU)\mu-\frac{\eps G(\omega)}{4 t^{2}}\mu U.\end{aligned}\end{equation}
In particular, note that $\mu_s+ G(\omega)\mu^2U_q/4=0$.

Fix $(t,r,\omega)\in\D$, so now we have $t\sim r$. Integrate this equation along the characteristic $(\tau,r+t-\tau,\omega)$. Note that $U$ vanishes if $\tau\geq t_1$ and $U,U_s=O(t^{C\eps})$. We have \fm{\int_t^{r+t}\frac{\eps|G(\omega)|}{4\tau^{2}}|\mu U|\ d\tau&\leq \int_{t}^{t_1}\frac{C\eps(t+r)^{C\eps}}{4 t^{2}}|\mu|\ d\tau\lesssim \eps t^{-2+C\eps}}
and \fm{\int_t^{r+t}\abs{\frac{\eps^2G(\omega)}{4 \tau^2}(U_s+\frac{1}{2}G(\omega)AU)\mu}\ d\tau&\leq C\eps^2\int_{t}^{t_1}\frac{(t+r)^{C\eps}}{t^2}|\mu|\ d\tau\lesssim \eps^2t^{-2+C\eps}.}
Finally, since $\int_t^{r+t}|\mu_q|\ d\tau\lesssim \eps\ln(t+r)+1\lesssim\eps\ln(t)+1$ and since $\nu=U=0$ at $\tau=r+t$, by Gronwall's inequality we conclude \eqref{lq41}.

To prove \eqref{lq42}, we first prove it with $\partial_q$ replaced by $\partial_r$. By \eqref{lq43}, we have 
\begin{equation}\label{lq44}\begin{aligned}(\partial_t-\partial_r)\partial_r(\nu+\frac{\eps G(\omega)}{4  t}\mu U)&=\partial_r(\partial_t-\partial_r)(\nu+\frac{\eps G(\omega)}{4  t}\mu U)\\&=\mu_{q}\partial_r(\nu+\frac{\eps G(\omega)}{4  t}\mu U)+q_r\mu_{qq}(\nu+\frac{\eps G(\omega)}{4  t}\mu U)\\&\hspace{1em}-\frac{\eps G(\omega)}{4 t^{2}}(\mu U)_qq_r+\frac{\eps^2G(\omega)}{4 t^2}\partial_q((U_{s}+\frac{1}{2}G(\omega)AU)\mu)q_r.\end{aligned}\end{equation}
Again, by integrating along the characteristic, we have 
\fm{\int_t^{r+t}|q_r\mu_{qq}(\nu+\frac{\eps G(\omega)}{4 \tau}\mu U)|\ d\tau&\lesssim \int_{[t,r+t]\cap[t_0,t_1]}(|\nu|+|\mu|)|\mu|\eps\tau^{-2+C\eps}\ d\tau\\&\lesssim \eps t_0^{-2+C\eps}\lesssim \eps t^{-2+C\eps},}
\fm{\int_t^{r+t}|\frac{\eps G(\omega)}{4\tau^{2}}(\mu U)_qq_r|\ d\tau&\lesssim\int_{[t,r+t]\cap[t_0,t_1]}\eps\tau^{-2}(|A|+|\mu_q U|)(|\mu|+|\nu|)\ d\tau\\&\lesssim\int_{[t,r+t]\cap[t_0,t_1]}\eps\tau^{-2+C\eps}|\mu|\ d\tau+\int_{[t,r+t]\cap[t_0,t_1]}\eps^2\tau^{-3+C\eps}\ d\tau\\&\lesssim \eps t_0^{-2+C\eps}\lesssim \eps t^{-2+C\eps}}and
\fm{\int_t^{r+t}|\frac{\eps^2G(\omega)}{4\tau^2}\partial_q((U_{s}+\frac{1}{2}G(\omega)AU)\mu)q_r|\ d\tau&\lesssim \int_{[t,r+t]\cap[t_0,t_1]}\eps^2\tau^{-2+C\eps}|\mu|\ d\tau\\&\lesssim \eps^2 t_0^{-2+C\eps}\lesssim \eps^2 t^{-2+C\eps}.}
Again, since $\int_t^{r+t}|\mu_q|\ d\tau\lesssim \eps\ln(t+r)+1$, $\nu=U=0$ at $\tau=r+t$ and $\partial_r=q_r\partial_q$ with \fm{|q_r|=\frac{|\nu-\mu|}{2}\geq \frac{|\mu|}{2}-\frac{|\nu|}{2}\gtrsim t^{-C\eps}-C\eps t^{-1+C\eps}\gtrsim t^{-C\eps},\hspace{2em}\text{if }\eps\ll 1,} by Gronwall's inequality we conclude \eqref{lq42}.

\end{proof}

\rmk{\rm The proof of \eqref{lq42} actually gives an estimate not just in $\D$. For $t\geq t_1$, everything in \eqref{lq42} is 0. For $t_0\leq t\leq t_1$, we have $t\sim r$. For $t\leq t_0$, we have proved in Lemma \ref{lq2} that $t_0\gtrsim t+r$. By applying Gronwall's inequality, we thus have 
\begin{equation}\label{lq422}
\nu_q+\frac{\eps G(\omega)}{4t} \mu_q U+\frac{\eps G(\omega)}{4t}\mu U_q=O(\eps(t+r)^{-2+C\eps} t^{C\eps})=O(\eps(t+r)^{-2+C\eps})
\end{equation}
for all $t\geq T_R\gg 1$.
}
\rm\bigskip

Most of the estimates in the previous three lemmas  will still hold in $\D$, if $Z^I$ is applied to the left hand sides for each multiindex $I$.

\lem{\label{lq3} We have $q\in S^0$. In other words, for all $(t,r,\omega)\in\D$ and  $I$,
\begin{equation}\label{lq31}|Z^Iq(t,r,\omega)|\lesssim_I t^{C_I\eps}.\end{equation} 
Similarly, we have  $\partial_s^a\partial_q^b\partial_\omega^c(\mu,U,A)\in S^0$ for all $a,b,c$. We also have $\nu\in \eps S^{-1}, \ \lambda_i\in S^{-1}$, and\fm{\nu+\frac{\eps G(\omega)}{4 t}\mu U\in \eps S^{-2},\ \nu_q+\frac{\eps G(\omega)}{4 t} \mu_q U+\frac{\eps G(\omega)}{4 t}\mu U_q\in \eps S^{-2}.}}
\begin{proof}
We use $f_0$ to denote any element in $S^0$ which is defined in Section \ref{spre}. For example, $f_0$ can be a finite sum of terms of the form \fm{CZ^{I_1}\omega_{i_1}Z^{I_2}\omega_{i_2}\cdots Z^{I_p}\omega_{i_p}.} We allow $f_0$ to vary from line to line. By Lemma \ref{lemons}, we have $Zf_0=f_0$. In addition, \fm{\ [Z,\partial_t-\partial_r]&=[Z,\partial_t]+\sum_i([Z,\omega_i]\partial_i+\omega_i[Z,\partial_i])=f_0\partial,\\ \partial_\alpha&=f_0(\partial_t-\partial_r)+f_0(\partial_t+\partial_r)+\sum_if_0(\partial_i-\omega_i\partial_r).}  We claim that for all multiindices $I$, we have \begin{equation}\label{lq35}
(\partial_t-\partial_r)Z^Iq=Z^I\mu+\sum_{|J|< |I|} f_0Z^J\mu+\sum_{|J|< |I|}[f_0(\partial_t+\partial_r)Z^Jq+\sum_if_0(\partial_i-\omega_i\partial_r)Z^Jq].
\end{equation}
To see this, we first prove that
\begin{equation}\label{lq34}(\partial_t-\partial_r)Z^Iq=Z^I\mu+\sum_{|J|<|I|} f_0\partial Z^Jq.\end{equation}
We can prove this by induction on $|I|$. For $|I|=0$, \eqref{lq34} is \eqref{qeqn}. In general, by writing $Z^I=ZZ^{I^\prime}$, we have 
\fm{(\partial_t-\partial_r)Z^Iq&=Z(\partial_t-\partial_r)Z^{I^\prime}q+[\partial_t-\partial_r,Z]Z^{I^\prime}q\\&=ZZ^{I^\prime}\mu+\sum_{|J|<|I^\prime|}Z(f_0\partial Z^Jq)+f_0\partial Z^{I^\prime}q\\&=Z^I\mu+\sum_{|J|<|I|}f_0\partial Z^Jq.}
We use the induction hypothesis on the second line. Moreover, for each $I$ we have
\fm{\partial Z^Iq&=f_0(\partial_t-\partial_r)Z^Iq+f_0(\partial_t+\partial_r)Z^Iq+\sum_if_0(\partial_i-\omega_i\partial_r)Z^Iq. }
It follows from \eqref{lq34} that
\fm{(\partial_t-\partial_r)Z^Iq=Z^I\mu+\sum_{|J|< |I|} [f_0(\partial_t-\partial_r)Z^Jq+f_0(\partial_t+\partial_r)Z^Jq+\sum_if_0(\partial_i-\omega_i\partial_r)Z^Jq].}
Then \eqref{lq35}  follows  by induction.

Now we prove \eqref{lq31} by induction on $|I|$. The case $|I|=0$  is obvious since $|q|\leq R$ in $\D$. Suppose \eqref{lq31} holds for all $|I|\leq k$ and fix $|I|=k+1$. By the chain rule and Leibniz's rule, $Z^I\mu$ can be expressed as as a sum of terms of the form
\begin{equation}\label{formofmu}
f_0\cdot \partial_q^{a}\partial_s^b\partial_{\omega}^{c}\mu\cdot Z^{I_1}q\cdots Z^{I_a}qZ^{J_1}(\eps\ln (t)-\delta)\cdots Z^{J_b}(\eps\ln (t)-\delta) \prod_lZ^{K_{l,1}}\omega_l\cdots Z^{K_{l,c_l}}\omega_l\end{equation}where all $|I_*|,|J_*|,|K_{*,*}|$ are nonzero, and the sum of them is $k+1$. The only term with some $|I_*|>k$ is $\mu_qZ^Iq$; all the other terms are controlled by $ t^{C\eps}|\mu|$, by induction hypotheses and \eqref{mu}. For the same reason, if $|J|\leq k$,  we have $Z^J\mu=O( t^{C\eps}|\mu|)$. Moreover, by Lemma \ref{l21}, we have
\fm{|(\partial_t+\partial_r)Z^Jq|+\sum_i|(\partial_i-\omega_i\partial_r)Z^Jq|\lesssim (1+t+r)^{-1}\sum_{|K|=1}|Z^KZ^Jq|.}In conclusion, if we integrate \eqref{lq35} along the characteristic and take the sum over all $|I|=k+1$, we have
\fm{&\hspace{1em}\sum_{|I|=k+1}|Z^Iq(t_1,r+t-t_1,\omega)-Z^Iq(t,r,\omega)|\\&\lesssim\int_t^{t_1}(|\mu_q|+(1+t+r)^{-1})\sum_{|I|=k+1}|Z^Iq|+\tau^{C\eps}|\mu|+\tau^{-1+C\eps}\ d\tau\\&\lesssim\int_t^{t_1}(|\mu_q|+(1+t+r)^{-1})\sum_{|I|=k+1}|Z^Iq|\ d\tau+2R t^{C\eps}+ t^{-1+C\eps}(t_1-t).}
Also note that $Z^Iq(t_1,t+r-t_1,\omega)=O(1)$, since $Z(r-t)=O(1)$ when $r=t-R$ and $t\gg_R 1$. Thus we are done if we apply Gronwall's inequality and Lemma \ref{lq1}.

Now we finish the proof on the remaining statements. By the chain rule and Leibniz's rule, we can expand  $Z^IU$ as a sum of terms of the form \eqref{formofmu} with $\mu$ there replaced by $U$. We then get $U\in S^0$ from $q\in S^0$. Similarly, we can show $\partial_s^a\partial_q^b\partial_\omega^c(\mu,U,A)\in S^0$ for all $(a,b,c)$. 

In addition, since $\lambda_i=(\partial_i-\omega_i\partial_r)q$ with $q\in S^0$, by Lemma \ref{lemons} we have $\lambda_i\in S^{-1}$. To prove \fm{\nu\in \eps S^{-1},\ \nu+\frac{\eps G(\omega)}{4 t}\mu U\in \eps S^{-2},\ \nu_q+\frac{\eps G(\omega)}{4 t} \mu_q U+\frac{\eps G(\omega)}{4 t}\mu U_q\in \eps S^{-2},} we only need to prove the second estimate for $\nu+\frac{\eps G(\omega)}{4 t} \mu U$. If we get the second estimate, the first estimate follows from $\frac{\eps G(\omega)}{4 t} \mu U\in \eps S^{-1}$, and the third estimate follows from $\partial_q=q_r^{-1}\partial_r$ and $q_r^{-1}\in S^0$. Here $q_r^{-1}\in S^0$ since $q_r\in S^0$ and $|q_r|\gtrsim t^{-C\eps}$ as proved in Lemma \ref{lq4}.

To prove the second estimate, we start with \eqref{lq43}. If we let $f_m$ denote any element in $S^m$ where $f_m$ is allowed to vary from line to line,  we have
\fm{
(\partial_t-\partial_r)V&=\mu_qV+\eps f_{-2}\mu,\hspace{1cm}\text{where }V=\nu+\frac{\eps G(\omega)}{4  t}\mu U.}
We can prove by induction that for each $I$,
\fm{(\partial_t-\partial_r)Z^IV&=\mu_qZ^IV+\sum_{|J|<|I|}f_0Z^JV+\eps \sum_{|J|\leq |I|}f_{-2}Z^J\mu\\&\hspace{1em}+\sum_{|J|<|I|}(f_0 (\partial_t+\partial_r) Z^JV+\sum_if_0(\partial_i-\omega_i\partial_r) Z^JV).}
Now we can induct on $|I|$ and conclude $|Z^IV|\lesssim_I\eps t^{-2+C_I\eps}$  by applying Gronwall's inequality. The proof here is very similar to that of $q$. 

\end{proof}
\rm

The next lemma can be viewed as a direct application of Lemma \ref{lq3}.

\lem{\label{lq5} In $\D$ we have
\begin{equation}\label{lq51}g^{\alpha\beta}q_\alpha q_\beta-\frac{1}{4}G(\omega)\mu^2\in S^{-1},\end{equation}
\begin{equation}\label{lq52}q_t^2-\sum_iq_i^2-\mu\nu\in S^{-2}.\end{equation}
Then we have
\begin{equation}\label{lq5ae}\wt{g}^{\alpha\beta}(\eps r^{-1}U)q_\alpha q_\beta\in S^{-2}.\end{equation}
In other words, $q$ is an approximate optical function in $\D$.

Moreover,
\begin{equation}g^{\alpha\beta}q_{\alpha\beta}-\frac{1}{4}G(\omega) \mu_q\mu\in S^{-1},\end{equation}
\begin{equation}\square q-(\frac{1}{2}(\mu_q\nu+\nu_q\mu)+\frac{\eps}{2t}\mu_s+r^{-1}\mu)\in S^{-2}.\end{equation}
}
\begin{proof}
Since $q_t=(\mu+\nu)/2$ and $q_i=\lambda_i+\omega_i(\nu-\mu)/2$, by applying Lemma \ref{lq2} we have 
\fm{&\hspace{1.25em}g^{\alpha\beta}q_\alpha q_\beta\\&=g^{00}(\frac{\mu+\nu}{2})^2+2g^{0i}(\frac{\mu+\nu}{2})(\lambda_i+\frac{\omega_i(\nu-\mu)}{2})+g^{ij}(\lambda_i+\frac{\omega_i(\nu-\mu)}{2})(\lambda_j+\frac{\omega_j(\nu-\mu)}{2})\\&=\frac{1}{4}G(\omega)\mu^2+\frac{1}{2}(g^{00}-g^{ij}\omega_i\omega_j)\mu\nu+(g^{0i}-g^{ij}\omega_j)\mu\lambda_i+\frac{1}{4}g^{00}\nu^2+\frac{1}{2}g^{0i}\nu(2\lambda_i+\omega_i\nu)\\&\hspace{1em}+\frac{1}{4}g^{ij}(2\lambda_i+\omega_i\nu)(2\lambda_j+\omega_j\nu)\\&=\frac{1}{4}G(\omega)\mu^2\mod S^{-1}.}

Replace $g^{\alpha\beta}$ with $-m^{\alpha\beta}$, and note that $\sum_i\omega_i(\partial_i-\omega_i\partial_r)=\sum_i\omega_i\partial_i-\omega_i^2\partial_r=0$. Then 
\fm{q_t^2-\sum_iq_i^2&=\mu\nu-\sum_i\lambda_i^2=\mu\nu\mod S^{-2}.}
Thus we have
\fm{\wt{g}^{\alpha\beta}(\eps r^{-1}U)q_\alpha q_\beta&=q_t^2-\sum_i q_i^2+\gamma^{\alpha\beta}(\eps r^{-1}U)q_\alpha q_\beta\\&=\mu\nu+\frac{\eps}{4r} G(\omega)\mu^2 U -\sum_i\lambda_i^2+(\gamma^{\alpha\beta}(\eps r^{-1}U)-g^{\alpha\beta}\eps r^{-1}U)q_\alpha q_\beta\\&=\mu(\nu+\frac{\eps}{4t} G(\omega)\mu U)+\frac{t-r}{4tr}\mu ^2U\mod S^{-2}\in S^{-2}.}
Note that we have $(\gamma^{\alpha\beta}(\eps r^{-1}U)-g^{\alpha\beta}\eps r^{-1}U)\in \eps^2 S^{-2}$ by  Lemma \ref{lemons2}.

In addition, we have \fm{\eps t^{-1}\nu_s&=\nu_t-q_t\nu_q=(\nu_t+\nu_r)-\nu\nu_q,\\
\sum_l\partial_i\omega_l\nu_{\omega_l}&=\nu_i-q_i\nu_q=(\nu_i-\omega_i\nu_r)-\lambda_i\nu_q.
}
Both of these two terms are in $\eps S^{-2}$ by Lemma \ref{lemons}. Thus,
\fm{q_{tt}&=\partial_t(\frac{\mu+\nu}{2})=\frac{1}{4}\mu_q(\mu+\nu)+\frac{\eps}{2 t}\mu_s+\frac{1}{4}\nu_q(\mu+\nu)+\frac{\eps}{2 t}\nu_s\\&=\frac{1}{4}\mu_q\mu+\frac{1}{4}\mu_q\nu+\frac{\eps}{2 t}\mu_s+\frac{1}{4}\nu_q\mu\mod \eps S^{-2},\\
q_{ti}&=\partial_i(\frac{\mu+\nu}{2})=\frac{1}{2}(\mu_q+\nu_q)(\lambda_i+\frac{\omega_i(\nu-\mu)}{2})+\frac{1}{2}\sum_l(\mu_{\omega_l}+\nu_{\omega_l})\partial_i\omega_l\\&=-\frac{1}{4}\omega_i\mu_q\mu\mod S^{-1},\\
q_{ij}&=\partial_i(\lambda_j+\frac{\omega_j(\nu-\mu)}{2})\\&=\partial_i\lambda_j+\frac{1}{2}\partial_i\omega_j (\nu-\mu)+\frac{1}{2}\omega_j (\nu_q-\mu_q)(\lambda_i+\frac{\omega_i(\nu-\mu)}{2})+\frac{1}{2}\omega_j\sum_l(\mu_{\omega_l}+\nu_{\omega_l})\partial_i\omega_l\\&=\frac{1}{4}\omega_i\omega_j\mu\mu_q+\partial_i\lambda_j-\frac{1}{2}\mu\partial_i\omega_j -\frac{1}{4}\omega_j \mu_q(2\lambda_i+\omega_i\nu)\\&\hspace{1em}-\frac{1}{4}\omega_j \nu_q\omega_i\mu+\frac{1}{2}\omega_j\sum_l\mu_{\omega_l}\partial_i\omega_l\mod \eps S^{-2}.}
Thus, 
\fm{g^{\alpha\beta}q_{\alpha\beta}&=\frac{1}{4}G(\omega) \mu_q\mu\mod S^{-1}.}
Replace $g^{\alpha\beta}$ with $-m^{\alpha\beta}$, and note that $\sum_i\omega_i\partial_i\omega_j=\sum_i\omega_i\lambda_i=0$, $\sum_i\partial_i\omega_i=2/r$, 
$\sum_i\omega_i\partial_r\lambda_i=\partial_r(\sum_i\omega_i\lambda_i)=0$.
Thus, 
\fm{\square q&=-\sum_i(\frac{1}{4}\omega_i^2\mu\mu_q+\partial_i\lambda_i-\frac{1}{2}\mu\partial_i\omega_i -\frac{1}{4}\omega_i \mu_q(2\lambda_i+\omega_i\nu)-\frac{1}{4}\omega_i^2 \nu_q\mu+\frac{1}{2}\omega_i\sum_l\mu_{\omega_l}\partial_i\omega_l)\\&\hspace{1em}+\frac{1}{4}\mu_q\mu+\frac{1}{4}\mu_q\nu+\frac{\eps}{2 t}\mu_s+\frac{1}{4}\nu_q\mu\mod S^{-2}\\&=\frac{1}{2}\mu_q\nu+\frac{\eps}{2 t}\mu_s+\frac{1}{2}\nu_q\mu+r^{-1}\mu-\sum_i(\partial_i-\omega_i\partial_r)\lambda_i \mod S^{-2}\\&=\frac{1}{2}\mu_q\nu+\frac{\eps}{2 t}\mu_s+\frac{1}{2}\nu_q\mu+r^{-1}\mu \mod S^{-2}.}
\end{proof}

\rm Finally we prove that $\eps r^{-1}U$ has good pointwise bounds and is an approximate solution to \eqref{qwe} in $\D$.
\prop{\label{prop4}We have 
\begin{equation}\eps r^{-1}U\in \eps S^{-1},\hspace{1cm}\wt{g}^{\alpha\beta}(\eps r^{-1}U)\partial_\alpha\partial_\beta(\eps r^{-1}U)\in \eps S^{-3}.\end{equation}
In other word, for $(t,r,\omega)\in\D$,
\begin{equation}|Z^I(\eps r^{-1}U)|\lesssim_I \eps t^{-1+C_I\eps}, \end{equation}
\begin{equation}
|Z^I(\wt{g}^{\alpha\beta}(\eps r^{-1}U)\partial_\alpha\partial_\beta(\eps r^{-1}U))|\lesssim_I\eps t^{-3+C_I\eps}.
\end{equation}
Note that we  have a better bound for $\partial (\eps r^{-1}U)$: for all $(t,r,\omega)\in\D$, \begin{equation}|\partial(\eps r^{-1}U)|\lesssim \eps t^{-1}. \end{equation}

}
\begin{proof}
Since  $r^{-1}\in S^{-1}$ and $U\in S^0$ which is proved in Lemma \ref{lq3}, we have $\eps r^{-1}U\in \eps S^{-1}$ by Lemma \ref{lemons}. In addition,
\fm{\partial_t(\eps r^{-1}U)&=\eps r^{-1}( U_qq_t+U_s\eps t^{-1}),
\\ \partial_i(\eps r^{-1}U)&=\eps r^{-2}\omega_i U+\eps r^{-1}( U_q q_i+\sum_jU_{\omega_j}\partial_i\omega_j).}
Note that $2U_qq_t=U_q(\mu+\nu)=-2A+O(\eps t^{-1+C\eps})$ and $2U_qq_i=U_q(\nu-\mu)\omega_i+2U_q\lambda_i=2A\omega_i+O( t^{-1+C\eps})$. Thus, $|\partial(\eps r^{-1}U)|\lesssim \eps t^{-1}$.

Next, we have
\fm{U_{tt}&= -U_s\eps t^{-2}+2U_{sq}q_t\eps t^{-1}+U_{ss}\eps^2 t^{-2}+q_{tt}U_q+q_t^2U_{qq}\\&=U_qq_{tt}+U_{qq}q_t^2\mod S^{-1},\\
U_{it}&= U_{qq}q_tq_i+U_{\omega_lq}q_t\partial_i\omega_l+U_qq_{it}+U_{sq}q_i\eps t^{-1}+\sum_lU_{s\omega_l}\partial_i\omega_l\eps t^{-1}\\&=U_{qq}q_tq_i+U_qq_{it}\mod S^{-1}, \\
U_{ij}&=U_{qq}q_iq_j+\sum_lU_{q\omega_l}(q_i\partial_j\omega_l+U_{q\omega_l}q_j\partial_i\omega_l)+U_qq_{ij}+\sum_{l,l^\prime}U_{\omega_l\omega_{l^\prime}}\partial_i\omega_l\partial_j\omega_{l^\prime}\\&=U_{qq}q_iq_j+U_qq_{ij}\mod S^{-1}.}
Thus, 
\fm{\square U&=-U_s\eps t^{-2}+2U_{sq}q_t\eps t^{-1}+U_{ss}\eps^2 t^{-2}+q_{tt}U_q+q_t^2U_{qq}\\&\hspace{1em}-\sum_i(U_{qq}q_i^2+\sum_l2U_{q\omega_l}q_i\partial_i\omega_l+U_qq_{ii}+\sum_{l,l^\prime}U_{\omega_l\omega_{l^\prime}}\partial_i\omega_l\partial_i\omega_{l^\prime})\\&=U_{qq}(q_t^2-\sum_i q_i^2)+U_q\square q+2U_{sq}q_t\eps t^{-1}-\sum_{i,l}2U_{q\omega_l}q_i\partial_i\omega_l\mod  S^{-2}\\&=U_{qq}\mu\nu+U_q(\frac{1}{2}(\mu_q\nu+\nu_q\mu)+\frac{\eps}{2 t}\mu_s+r^{-1}\mu)+\frac{\eps}{t}\mu U_{sq}\mod  S^{-2}.}
The last line holds by Lemma \ref{lq5} and \fm{\sum_i q_i\partial_i\omega_l&=\sum_i\lambda_i\partial_i\omega_l+\sum_i\omega_iq_r\partial_i\omega_l=\sum_i\lambda_i\partial_i\omega_l\in S^{-2}.}
By Lemma \ref{lq3}, we get
\fm{\square U&=-U_{qq}\mu \cdot \frac{\eps G(\omega)}{4 t}\mu U-\frac{1}{2}U_q\mu_q\cdot \frac{\eps G(\omega)}{4 t}\mu U-\frac{1}{2}U_q\mu\cdot \frac{\eps G(\omega)}{4 t}\mu U_q-\frac{1}{2}U_q\mu\cdot \frac{\eps G(\omega)}{4 t}\mu_q U\\&\hspace{1em}+\frac{\eps}{2 t}\mu_sU_q+r^{-1}\mu U_q+\frac{\eps}{t}\mu U_{sq}\mod S^{-2}\\&=(2A_q+\mu_qU_q) \cdot \frac{\eps G(\omega)}{4 t}\mu U-\frac{1}{2}U_q\mu_q\cdot \frac{\eps G(\omega)}{4 t}\mu U-\frac{\eps G(\omega)}{2 t}A^2+ \frac{\eps G(\omega)}{4 t}A\mu_q U\\&\hspace{1em}+\frac{\eps}{4 t}GA\mu U_q-2Ar^{-1}-\frac{\eps}{2 t}GA\mu U_{q}\mod S^{-2}\\&= \frac{\eps G(\omega)}{2 t}A_q\mu U-2Ar^{-1}\mod S^{-2}.}

Besides, by  Lemma \ref{lemons2} we have
\fm{\gamma^{\alpha\beta}(\eps r^{-1}U)-g^{\alpha\beta}\eps r^{-1}U\in \eps^2S^{-2}} 
Thus,
\fm{&\hspace{1.25em}\wt{g}^{\alpha\beta}(\eps r^{-1}U)\partial_\alpha\partial_\beta(\eps r^{-1}U)\\&=\square(\eps r^{-1}U)+\gamma^{\alpha\beta}(\eps r^{-1}U)(\eps r^{-1}U)_{\alpha\beta}\\&=\eps r^{-1}\square U+2\eps r^{-2}\omega_i U_i+g^{\alpha\beta}\eps^2 r^{-2}UU_{\alpha\beta}-2g^{i\beta}\eps^2 r^{-3}\omega_iUU_\beta\\&\hspace{1em}+g^{ij}\eps^2r^{-4}U^2(3\omega_i\omega_j-\delta_{ij})+(\gamma^{\alpha\beta}(\eps r^{-1}U)-g^{\alpha\beta}\eps r^{-1}U)(\eps r^{-1}U)_{\alpha\beta}\\&=\eps r^{-1}\square U+2\eps r^{-2}U_qq_r+g^{\alpha\beta}\eps^2 r^{-2}UU_{\alpha\beta}\mod \eps^2S^{-3}.}

Note that 
\fm{&\hspace{1em}\eps r^{-1}\square U+2\eps r^{-2} U_qq_r+g^{\alpha\beta}\eps^2 r^{-2}UU_{\alpha\beta}\\&=\eps r^{-1}( \frac{\eps G(\omega)}{2 t}A_q\mu U-2Ar^{-1})-\eps r^{-2}\mu U_q+\eps^2r^{-2}U(U_{qq}g^{\alpha\beta}q_\alpha q_\beta+U_qg^{\alpha\beta}q_{\alpha\beta})\mod \eps S^{-3}
\\&=\eps r^{-1}( \frac{\eps G(\omega)}{2 t}A_q\mu U-2Ar^{-1})+2\eps r^{-2}A\\&\hspace{1em}+\eps^2r^{-2}U(U_{qq}\cdot\frac{1}{4}G(\omega)\mu^2+U_q\cdot\frac{1}{4}G(\omega)\mu\mu_q)\mod \eps S^{-3}
\\&=\eps r^{-1} \frac{\eps G(\omega)}{2 t}A_q\mu U+\eps^2r^{-2}U(-(2A_q+\mu_qU_q)\frac{1}{4}G(\omega)\mu+U_q\cdot\frac{1}{4}G(\omega)\mu\mu_q)\mod \eps S^{-3}
\\&=\frac{\eps^2(r-t)}{2r^2 t}A_q\mu U\mod \eps S^{-3}\in \eps S^{-3}.
}
The last equality holds  since $r-t\in S^0$. Done.

\end{proof}
\rm

\subsection{Approximate solution $u_{app}$}\label{sec4.3} 

Let $T_R$ be the constant used in the definition of $\D$, such that all the estimates in Section \ref{sec4.2} hold for $t\geq T_R$. Choose $\eta\in C^\infty(\R)$ such that $\eta\equiv 1$ on $[2T_R,\infty)$ and $\eta\equiv 0$ on $(-\infty,T_R]$.  In addition, choose $\psi\in C_c^\infty(\R)$ such that $\psi\equiv 1$ on $[3/4,5/4]$ and $\psi\equiv 0$ outsides $[1/2,3/2]$.

We now define the approximate solution $u_{app}$ by
\begin{equation}\label{uappdef}u_{app}(t,x):=\eps r^{-1}\eta(t)\psi(r/t)U(\eps\ln(t)-\delta,q(t,r,\omega),\omega),\hspace{1cm}r=|x|,\ \omega_i=x_i/r.\end{equation}
Note that $u_{app}(t,x)$ is defined for all $(t,x)\in[0,\infty)\times\R^{3}$. If $t\leq T_R$, then we have $u_{app}\equiv 0$. If $t\geq T_R\geq 2R$, since $U\equiv 0$ for $r\leq t-R$, $u_{app}$ has no singularity at $|x|=0$.  Moreover, since $\psi\equiv 1$ when $|t-r|<t/4$, we have $\psi\equiv 1$ in $\D$ if $t\geq T_R\gg 1$; since $\psi\equiv 0$ when $|t-r|>t/2$, we have $u_{app}\equiv 0$ unless $t\sim r$.

We now prove the  estimates for $u_{app}$ in Proposition \ref{mainprop4}.  The estimates are in fact the same as those in Proposition \ref{prop4}. However, note that in Proposition \ref{prop4} we assume that $(t,r,\omega)\in\D$ while here we only assume $t\geq 0$.

\begin{proof}[Proof of Proposition \ref{mainprop4}]

When $t\leq T_R$, we have $u_{app}\equiv 0$. When $T_R\leq t\leq 2T_R$, we have $Z^Iu_{app}=O_R(\eps)$. This is because the support of $u_{app}$ lies in $|x|\sim_R 1$, and because $U,\eta,\psi$ and all their derivatives are  $O(1)$. Also note that $\eps\leq (2T_R)^M\eps t^{-M}$ for each $M$ and all $t\leq 2T_R$.

Suppose $t\geq 2T_R$. Now $\eta$ plays no role since $\eta(t)=1$ for all $t\geq 2T_R\gg 1$. For $(t,r,\omega)\in\D$, all the estimates  follow directly from Proposition \ref{prop4}. If $q(t,r,\omega)\leq -R$ i.e. $r-t\leq -R$, or if $r>3t/2$, then $u_{app}\equiv 0$ so there is nothing to prove. So now we can assume $t\geq 2T_R$, $q(t,r,\omega)\geq R$ and $t-R\leq r\leq 3t/2$. By construction, we now have $U(t,r,\omega)=U(\eps\ln(t)-\delta,R,\omega)$. By applying chain rule and Leibniz's rule, for all $k$ and $I$, we have $\partial^kZ^IU=O(t^{-k+C\eps})$.
Since $\partial^kZ^I(r/t)=O(t^{-k})$ for $t\sim r$ and $r>t-R\gg 1$, we have $\partial^kZ^I(\psi(r/t))=O(t^{-k})$ for all $t\geq T_R$. In particular, we have $\partial (\psi(r/t))=\psi^\prime\partial(r/t)=O(t^{-1})$. In addition, $\partial^kZ^Ir^{-1}=O(t^{-1-k})$ for $t\sim r$. Now all the estimates follow directly from the Leibniz's rule.

\end{proof}

\section{Energy estimates and Poincare's lemma}\label{sep}
\rm
We now derive the energy estimates and Poincar$\acute{\rm e}$'s lemma, which are the main tools in the proof of our main theorem. The results in this section are closely related to those in \cite{lind,alin2}.

\subsection{Setup}\rm Suppose $t\geq T_R\gg 1$ and $\eps\ll 1$. Assume that $u$ is a solution to \eqref{qwe} vanishing for $r\leq t-R$ and satisfying the pointwise estimates: for all $t\geq T_R\gg 1$ we have
\begin{equation}\label{assu}|u|\lesssim \eps t^{-1+C\eps},\hspace{1cm} |\partial u(t,x)|\lesssim \eps t^{-1};\end{equation}
if $q(t,r,\omega)\leq t^{1/4}$ and $t\geq T_R$, we have 
\begin{equation}|u-\eps r^{-1}U|\lesssim \eps t^{-5/4+C\eps}.\end{equation}Recall that $U=U(t,r,\omega)$ is the asymptotic profile defined in \eqref{Udef}. In Section \ref{smp} we will check these  estimates when we apply the energy estimates.

We first prove the following two lemmas. 
\lem{\label{le2}Suppose $t\geq T_R$. Then
\begin{equation}\label{le21}q_t=\frac{\mu+\nu}{2}=\frac{\mu}{2}+O(\eps (t+r)^{-1+C\eps}),\end{equation}

\begin{equation}\label{le22} q_i=\lambda_i+\frac{\omega_i(\nu-\mu)}{2}=-\frac{\omega_i\mu}{2}+O( (t+r)^{-1+C\eps}),\end{equation}

\begin{equation}\label{le23} q_t^2-\sum_iq_i^2=\mu\nu-\sum_i\lambda_i^2=\mu\nu+O( (t+r)^{-2+C\eps}).\end{equation}
}

\begin{proof} Note that \eqref{le21} and \eqref{le22} follow directly from Lemma \ref{lq2}, and  \eqref{le23} follows from   the proof of  Lemma \ref{lq5}. Note that these three inequalities hold for all $(t,r,\omega)$ with $t\geq T_R$, not just in $\D$.
\end{proof}

\lem{\label{le1}Suppose $t\geq T_R$ and $-R\leq q(t,r,\omega)\leq  t^{1/4}$. Then, 
\begin{equation}\label{le11} -R\leq r-t\lesssim  t^{1/4},\end{equation}
\begin{equation}\label{le13} \nu+\frac{\eps G(\omega)}{4 t}\mu U=O(\eps t^{-7/4+C\eps}),\end{equation}
\begin{equation}\label{le12} g^{\alpha\beta} q_\alpha q_\beta=\frac{1}{4}G(\omega)\mu^2+O(t^{-1+C\eps}).\end{equation}}
\begin{proof} By Lemma \ref{lq1}, we have\fm{|r-t|&\leq |q|+C(t+r)^{C\eps}\lesssim t^{1/4}+(t+r)^{C\eps}.} This implies that \fm{0<r/t\lesssim 1+t^{-3/4}+t^{-1+C\eps}(1+r/t)^{C\eps} \lesssim1+(1+r/t)^{C\eps}.} Thus we have $r/t\lesssim 1$ and then we conclude \eqref{le11}.

The proof of \eqref{le13} is essentially the same as that of \eqref{lq41}. The only difference is that we have
\fm{\int_t^{t_1}|\mu(\tau,r+t-\tau,\omega)|\ d\tau=R+q(t,r,\omega)\leq R+t^{1/4},}while we put $2R$ on the right hand side in Lemma \ref{lq4}.

The estimate \eqref{le12} follows from the computations in Lemma \ref{lq5} and Lemma \ref{lq2}.
\end{proof}\rm

\subsection{Energy estimates} Fix a smooth function $\phi(t,x)$ with $\phi(t)\in C_c^\infty(\R^3)$ for each $t\geq T_R$ and $\phi$ is supported in $r\geq t-R$. We define the energy 
\begin{equation}\label{energydef}\begin{aligned}
E_u(\phi)(t)&=\int_{\R^3}w(t,x)(2\widetilde{g}^{0\alpha}(u)\phi_t\phi_\alpha-\widetilde{g}^{\alpha\beta}(u)\phi_\alpha\phi_\beta)(t,x)\ dx\\&=\int_{\R^3}w(t,x)(|\partial \phi|^2+2\gamma^{0\alpha}(u)\phi_t\phi_\alpha-\gamma^{\alpha\beta}(u)\phi_\alpha\phi_\beta)(t,x)\ dx.
\end{aligned}\end{equation}
The weight function $w$ is defined by\begin{equation}w(t,x)=\exp(c_0\eps\ln(t)\cdot \sigma(q(t,r,\omega)))\end{equation}with \fm{\sigma(q)=(R+q+1)^{-1/16}.} Here $q(t,r,\omega)$ is defined in Section \ref{sas}; $c_0>0$ is a large constant to be chosen, which depends only on the scattering data $A$. Note that $\phi\equiv 0$ unless $r\geq t-R$, and $q(t,r,\omega)\geq -R$ when $r\geq t-R$. So $w(t,x)$ is well-defined in the support of $\phi$.

We remark that the exact value of the power in $\sigma$ is not important. We can replace $-1/16$ with any fixed constant $\lambda\in(-1/8,0)$.

We also remark that this type of the weight $w$ was already used in the previous work on small data global existence by Lindblad \cite{lind} and Alinhac \cite{alin2}. It can be viewed as an extended version of the method of ghost weight introduced by Alinhac. See \cite{alin}.

Our goal  is to prove the following energy estimates.

\prop{\label{propene}For $1\ll T_R\leq t\leq T$, we have \begin{equation}\label{energymain}
E_u(\phi)(t)\leq E_u(\phi)(T)+\int_t^T2\norm{\wt{g}^{\alpha\beta}(u)\partial_\alpha\partial_\beta\phi(\tau)}_{L^2(w)}\norm{\partial\phi(\tau)}_{L^2(w)}+C\eps \tau^{-1}\norm{\partial\phi}^2_{L^2(w)}\ d\tau.
\end{equation}
Here $\norm{f}_{L^2(w)}^2:=\int_{\R^3}|f|^2w\ dx$ and  $C>0$ is a constant (may depend on $u,\partial u$).}\rm\bigskip

The proof starts with a computation of $\frac{d}{dt}E_u(\phi)(t)$. For simplicity, we write $\wt{g}^{\alpha\beta}=\wt{g}^{\alpha\beta}(u)$. Then, by applying integration by parts, we have 
\fm{&\hspace{1em}\frac{d}{dt}E_u(\phi)(t)\\&=\int_{\R^3}w_t(2\wt{g}^{0\alpha}\phi_t\phi_\alpha-\wt{g}^{\alpha\beta}\phi_\alpha\phi_\beta)\\&\hspace{3em}+w(2\wt{g}^{0\alpha}\phi_{tt}\phi_\alpha+2\wt{g}^{0\alpha}\phi_{t}\phi_{\alpha t}+2\partial_t\wt{g}^{0\alpha}\phi_t\phi_\alpha-2\wt{g}^{\alpha\beta}\phi_{\alpha t}\phi_\beta-\partial_t\wt{g}^{\alpha\beta}\phi_\alpha\phi_\beta)\ dx\\&=\int_{\R^3}w_t(2\wt{g}^{0\alpha}\phi_t\phi_\alpha-\wt{g}^{\alpha\beta}\phi_\alpha\phi_\beta)+w(2\wt{g}^{0\alpha}\phi_{\alpha t}\phi_t-2\wt{g}^{i\beta}\phi_{it}\phi_{\beta}+2\partial_t\wt{g}^{0\alpha}\phi_t\phi_\alpha-\partial_t\wt{g}^{\alpha\beta}\phi_\alpha\phi_\beta)\ dx\\&=\int_{\R^3}w_t(2\wt{g}^{0\alpha}\phi_t\phi_\alpha-\wt{g}^{\alpha\beta}\phi_\alpha\phi_\beta)+2w_i\wt{g}^{i\beta}\phi_t\phi_\beta\\&\hspace{3em}+w(2\wt{g}^{0\alpha}\phi_{\alpha t}\phi_t+2\wt{g}^{i\beta}\phi_{t}\phi_{i\beta}+2\partial_t\wt{g}^{0\alpha}\phi_t\phi_\alpha+2\partial_i\wt{g}^{i\beta}\phi_t\phi_\beta-\partial_t\wt{g}^{\alpha\beta}\phi_\alpha\phi_\beta)\ dx\\&=\int_{\R^3}-w_t\wt{g}^{\alpha\beta}\phi_\alpha\phi_\beta+w(2\wt{g}^{\alpha\beta}\phi_{\alpha \beta}\phi_t+2\partial_\alpha\wt{g}^{\alpha\beta}\phi_t\phi_\beta-\partial_t\wt{g}^{\alpha\beta}\phi_\alpha\phi_\beta)+2w_\alpha \wt{g}^{\alpha\beta}\phi_t\phi_\beta\ dx.}

By setting $T_\alpha:=q_t\partial_\alpha-q_\alpha\partial_t$, we have $\phi_\alpha=q_t^{-1}(T_\alpha \phi+q_\alpha\phi_t)$. Note that
\fm{w_t&=c_0(\eps t^{-1}\sigma(q) +\eps\ln(t)\sigma^\prime(q)q_t)w,\hspace{1cm} w_i=c_0\eps\ln(t)\sigma^\prime(q)q_iw.} Thus,
\fm{-\widetilde{g}^{\alpha\beta}\phi_\alpha\phi_\beta q_t+2\widetilde{g}^{\alpha\beta}\phi_t\phi_\beta q_\alpha&=-\widetilde{g}^{\alpha\beta}q_t^{-1}(T_\alpha\phi+q_\alpha\phi_t)(T_\beta\phi+q_\beta\phi_t)+2\widetilde{g}^{\alpha\beta}q_\alpha\phi_tq_t^{-1}(T_\beta\phi+q_\beta\phi_t)\\&=-\widetilde{g}^{\alpha\beta}q_t^{-1}T_\alpha\phi T_\beta\phi+\widetilde{g}^{\alpha\beta}q_t^{-1}q_\alpha q_\beta\phi_t^2}
and 
\fm{-w_t\widetilde{g}^{\alpha\beta}\phi_\alpha\phi_\beta+2w_\alpha \widetilde{g}^{\alpha\beta}\phi_t\phi_\beta&=c_0\eps\ln(t)\sigma^\prime(q)w(-\widetilde{g}^{\alpha\beta}q_t^{-1}T_\alpha\phi T_\beta\phi+\widetilde{g}^{\alpha\beta}q_t^{-1}q_\alpha q_\beta\phi_t^2)\\&\hspace{1em}+c_0\eps t^{-1}\sigma(q)w(\widetilde{g}^{00}\phi_t^2-\widetilde{g}^{ij}\phi_i\phi_j).}

Note that $T_0=0$, $(-\wt{g}^{ij})=(\delta_{ij}+O(\eps t^{-1+C\eps}))$ is positive definite for $\eps\ll 1$ and $t\geq T_R\gg 1$; $\sigma^\prime(q)=-\frac{1}{16} (R+q+1)^{-17/16}<0$; by Lemma \ref{lq2} we have \fm{q_t=(\mu+\nu)/2\leq -c t^{-C\eps}+C\eps (t+r)^{-1+C\eps}<0.} We  conclude that \fm{-c_0\eps\ln(t)\sigma^\prime(q)w\widetilde{g}^{\alpha\beta}q_t^{-1}T_\alpha\phi T_\beta\phi\geq 0.}

In addition, we claim that 
\begin{equation}\label{encl} |c_0\eps\ln(t)\sigma^\prime(q)w \widetilde{g}^{\alpha\beta}q_t^{-1}q_\alpha q_\beta\phi_t^2|\leq Cc_0\eps t^{-1}\sigma(q) w\phi_t^2.\end{equation}
\subsubsection*{Case 1.} Suppose $q(t,r,\omega)\leq t^{1/4}$.  By Lemma \ref{l24} and Lemma \ref{le1}, we have
\fm{\wt{g}^{\alpha\beta}(u)q_\alpha q_\beta&=q_t^2-\sum_iq_i^2+\gamma^{\alpha\beta}(\eps r^{-1}U)q_\alpha q_\beta+(\gamma^{\alpha\beta}(u)-\gamma^{\alpha\beta}(\eps r^{-1}U))q_\alpha q_\beta\\&=\mu\nu-\sum_i\lambda_i^2+g^{\alpha\beta}\eps r^{-1}Uq_\alpha q_\beta+g^{\alpha\beta} (u-\eps r^{-1}U)q_\alpha q_\beta+O(\eps^2 t^{-2+C\eps})\\&=-\frac{\eps G(\omega)}{4t}\mu^2U+\frac{\eps G(\omega)}{4r}\mu^2U+\frac{1}{4}G(\omega)\mu^2(u-\eps r^{-1}U)+O(t^{-2+C\eps})\\&=O(t^{-5/4+C\eps}).}Note that we use the assumption $u-\eps r^{-1}U=O(\eps t^{-5/4+C\eps})$ in the last inequality.
Since $\frac{1}{16}\ln(t)\leq t^{1/16}$ and $|q_t|\geq |\mu|/2-|\nu|/2\gtrsim t^{-C\eps}$, we have
\fm{|c_0\eps\ln(t)\sigma^\prime(q)w\widetilde{g}^{\alpha\beta}q_t^{-1}q_\alpha q_\beta\phi_t^2|&= c_0\eps\ln(t)\cdot\frac{1}{16}(q+R+1)^{-17/16} w|q_t^{-1}\wt{g}^{\alpha\beta}(u)q_\alpha q_\beta|\phi_t^2\\&\leq  c_0\eps t^{1/16} \sigma(q)(q+R+1)^{-1}\cdot Ct^{-5/4+C\eps}w\phi_t^2\\&\leq Cc_0\eps \sigma(q)t^{-1}w\phi_t^2.}
\subsubsection*{Case 2.} Suppose $q(t,r,\omega)\geq t^{1/4}$. Then 
\fm{(q+R+1)^{-1}|\wt{g}^{\alpha\beta}(u)q_\alpha q_\beta|&= (q+R+1)^{-1}|q_t^2-\sum q_i^2+\gamma^{\alpha\beta}(u)q_\alpha q_\beta|\\&\lesssim t^{-1/4} (|\mu\nu|+\sum_i\lambda_i^2+O(|u||\partial q|^2))\\&\lesssim t^{-5/4+C\eps}.}Here we use \eqref{assu} and Lemma \ref{le2}. So we also have 
\fm{|c_0\eps\ln(t)\sigma^\prime(q)w\widetilde{g}^{\alpha\beta}q_t^{-1}q_\alpha q_\beta\phi_t^2|&= c_0\eps\ln(t)\cdot\frac{1}{16}(q+R+1)^{-17/16} w|q_t^{-1}\wt{g}^{\alpha\beta}(u)q_\alpha q_\beta|\phi_t^2\\&\leq c_0\eps t^{1/16}\sigma(q) (q+R+1)^{-1}|q_t^{-1}\wt{g}^{\alpha\beta}(u)q_\alpha q_\beta|w\phi_t^2\\&\leq Cc_0\eps t^{-1}\sigma(q)w\phi_t^2.}
Now we finish the proof of \eqref{encl}.

Since \fm{\widetilde{g}^{00}\phi_t^2-\widetilde{g}^{ij}\phi_i\phi_j&=|\partial\phi|^2+O(|u||\partial\phi|^2)\sim |\partial\phi|^2,}
 we have 
\fm{-\widetilde{g}^{\alpha\beta}\phi_\alpha\phi_\beta w_t+2\widetilde{g}^{\alpha\beta}\phi_t\phi_\beta w_\alpha\geq -Cc_0\eps t^{-1}\sigma(q)w|\partial\phi|^2.}
In conclusion, 
\fm{\frac{d}{dt}E_u(\phi)(t)&\geq \int_{\R^3}w(2\wt{g}^{\alpha\beta}\phi_{\alpha \beta}\phi_t+2\partial_\alpha\wt{g}^{\alpha\beta}\phi_t\phi_\beta-\partial_t\wt{g}^{\alpha\beta}\phi_\alpha\phi_\beta)-Cc_0\eps t^{-1}\sigma(q)w|\partial\phi|^2\ dx\\&\geq \int_{\R^3}-2w|\wt{g}^{\alpha\beta}\phi_{\alpha \beta}||\phi_t|-C\eps t^{-1}w|\partial\phi|^2\ dx\\&\geq -2\norm{\wt{g}^{\alpha\beta}\phi_{\alpha \beta}}_{L^2(w)}\norm{\phi_t}_{L^2(w)}-C\eps t^{-1}\norm{\partial\phi}_{L^2(w)}^2. }
Integrate this inequality with respect to $t$ on $[t,T]$ and we conclude \eqref{energymain}.

\subsection{Poincar$\acute{\bf e}$'s lemma} \label{sec5.2}
Fix a smooth function $\phi(t,x)$ with $\phi(t)\in C_c^\infty(\R^3)$ for each $t\geq T_R$ and $\phi$ is supported in $r\geq t-R$. As in the previous sections, we shall assume that $t\geq T_R\gg 1$ and $\eps\ll 1$.

\lem{\label{lp1}For $\phi$ as above, we have\begin{equation}\label{lp11}
\int_{\R^3}\lra{t-r}^{-2}|\phi|^2\ dx\lesssim \int_{\R^3}|\partial\phi|^2\ dx.
\end{equation}}
\begin{proof}We have
\fm{\int \lra{t-r}^{-2}|\phi|^2\ dx&\lesssim_R\int_{\mathbb{S}^2}\int_0^\infty(r-t+R+1)^{-2}|\phi|^2\ r^2drdS_\omega\\&=\int_{\mathbb{S}^2}\int_0^\infty|\phi|^2\ r^2\partial_r(-(r-t+R+1)^{-1})\ drdS_\omega\\&=\int_{\mathbb{S}^2}\int_0^\infty\partial_r(|\phi|^2 r^2)(r-t+R+1)^{-1}\ drdS_\omega\\&=\int_{\mathbb{S}^2}\int_0^\infty(2|\phi|^2r+2\phi\phi_rr^2)(r-t+R+1)^{-1}\ drdS_\omega\\&\lesssim_R\int_{\mathbb{S}^2}\int_{0}^\infty 2|\phi r^{-1}+\phi_r|\cdot |\phi|\lra{t-r}^{-1}\ r^2drdS_\omega\\&\lesssim \left(\int \lra{t-r}^{-2}|\phi|^2\ dx\right)^{1/2}\left(\int |\phi r^{-1}+\phi_r|^2\ dx\right)^{1/2}.}Since 
\fm{\int 2\phi \phi_rr^{-1}\ dx&=\int_{\mathbb{S}^2}\int_0^\infty \partial_r(\phi^2) r\ drdS_\omega=\int_{\mathbb{S}^2}\int_0^\infty -\phi^2 \ drdS_\omega=-\int \phi^2r^{-2}\ dx,}
we have 
\fm{\int |\phi r^{-1}+\phi_r|^2\ dx&=\int\phi_r^2 \ dx.}We then conclude \eqref{lp11}.
\end{proof}
\rm
We can also prove a weighted version of the Poincar$\acute{\rm e}$'s lemma. Note that the value of $\delta$ in $s=\eps\ln(t)-\delta$ is chosen in the proof of this lemma.

\lem{\label{lp2}For $\phi$ as above, we have\begin{equation}\label{lp21}
\int\phi^2q_r^2\lra{q}^{-2}w\ dx\lesssim\int|\partial\phi|^2w\ dx.
\end{equation}}

\begin{proof}
Note that $\lra{q}\sim(q+R+1)$ since $\phi$ is supported in $q\geq -R$.

If $|q|\leq R$, we have\fm{&\hspace{1.5em}\partial_q(q_r)w+q_rw_q\\&=w(\partial_q(q_r)-q_r c_0\eps\ln(t)\cdot\frac{1}{16} (q+R+1)^{-17/16})\\&\leq \frac{1}{2}w(\nu_q-\mu_q-c_0\eps\ln(t)\cdot\frac{1}{16} (q+R+1)^{-17/16}(\nu-\mu))\\&\leq -\frac{w}{2}(-\frac{1}{2}G(\omega) A_q(\eps\ln(t)-\delta)+\frac{c_0\eps}{16}\ln(t)\cdot(q+R+1)^{-17/16})|\mu|+O(\eps t^{-1+C\eps}w)\\&\leq  (C(\eps\ln(t)+\delta)-\frac{c_0\eps}{16}(2R+1)^{-17/16}\ln(t))w|\mu|+C\eps t^{-1+C\eps}w\\&\leq C\delta wq_r+\eps(C\ln(t)|\mu|-c_0C^{-1}\ln(t)|\mu|+Ct^{-1+C\eps})w\\&\leq C\delta wq_r.}
The last inequality holds because if $c_0\gg_R 1$, the second term is negative. Also note that $|\mu|\gg t^{-1+C\eps}$.

If $q>R$, then $\mu\equiv -2$, so $q_r=(\nu-\mu)/2=1+O(\eps (t+r)^{-1+C\eps})\in(1/2,2)$ by \eqref{lq21}, and $\partial_q(q_r)=\nu_q/2=O(\eps(t+r)^{-2+C\eps})$ by \eqref{lq422}. Besides, by Lemma \ref{lq1}, we have \fm{2R+1<q+R+1=r-t+R+1+O((r+t)^{C\eps})\leq 2(r+t).}  Then,
\fm{&\hspace{1.5em}\partial_q(q_r)w+q_rw_q\\&=w(\partial_q(q_r)-q_r c_0\eps\ln(t)\cdot\frac{1}{16} (q+R+1)^{-17/16})\\&\leq w(C\eps (t+r)^{-2+C\eps}-\frac{c_0}{32}\eps\ln(t) \cdot(q+R+1)^{-17/16})\\&\leq w(C\eps (t+r)^{-2+C\eps}-\frac{c_0}{64}\eps\ln(t)\cdot (t+r)^{-17/16})\\&\leq 0.}The last inequality holds if we have $\eps\ll 1$, $t\geq T_R\gg 1$ and $c_0\gg_R 1$.

Now we have
\fm{&\hspace{1.25em}\int |\phi|^2q_r^2\lra{q}^{-2}w\ dx\\&\lesssim_R\int_{\mathbb{S}^2}\int_0^\infty |\phi(t,r\omega)|^2r^2q_r^2(q+R+1)^{-2}w\ drdS_\omega\\&=\int_{\mathbb{S}^2}\int_0^\infty (q+R+1)^{-1}\partial_r(\phi^2r^2q_rw)\ drdS_{\omega}\\&=\int_{\mathbb{S}^2}\int_{0}^\infty (q+R+1)^{-1}[2\phi\phi_rr^2w+2\phi^2rw+\phi^2r^2\partial_q(q_r)w+\phi^2r^2q_rw_q]q_r\ drdS_{\omega}\\&\leq\int_{\mathbb{S}^2}\int_{0}^\infty (q+R+1)^{-1}(2\phi\phi_r+2\phi^2r^{-1})r^2q_rw\ drdS_{\omega}\\&\hspace{1em}+\int_{\mathbb{S}^2}\int_{0}^\infty (q+R+1)^{-1}\phi^2r^2\cdot C\delta q_rw\chi_{|q|\leq R}\cdot q_r\ drdS_{\omega}\\&\leq 4
\kh{\int( |\phi_r|^2+r^{-2}|\phi|^2)w\ dx}^{1/2}\kh{\int \phi^2\lra{q}^{-2}q_r^2w\ dx}^{1/2}\\&\hspace{1em}+C_R\delta\int \lra{q}^{-2}\phi^2 q_r^2w\ dx. }
Here we use the fact that $\lra{q}\sim_R1$ when $|q|\leq R$. Here $C_R$ in the second term only depends on $R$ and on the scattering data, and in particular it does not depend on $\eps$, $t$ or $T_R$. Thus, by choosing $\delta:=\frac{1}{4C_R}$, we conclude that 
\fm{\int |\phi|^2q_r^2\lra{q}^{-2}w\ dx&\lesssim_R\int( |\phi_r|^2+r^{-2}|\phi|^2)w\ dx. }

Now recall that $r\geq t-R$ when $\phi\neq 0$. If $q\leq t^{1/2}$ we have $\lra{q}^2\leq Ct$ and $q_r\geq C^{-1}t^{-C\eps}$, as proved before. Thus, if $t\geq T_R\gg 1$,\fm{\int_{q\leq t^{1/2}} r^{-2}\phi^2w\ dx&\lesssim (t-R)^{-2}\cdot Ct^{C\eps}\cdot Ct\int \phi^2 q_r^2 \lra{q}^{-2} w\ dx\\&\lesssim_R t^{-1+C\eps}\int \phi^2q_r^2\lra{q}^{-2} w\ dx.}

If $q\geq t^{1/2}$, we have $w(q)\leq\exp(Cc_0\eps\ln(t)\cdot t^{-1/32})\leq C$ for $t\gg_R1$ and $\eps\ll 1$. Besides, we also have $w\geq 1$.  Thus, by Hardy's inequality,
\fm{\int_{q\geq t^{1/2}} r^{-2}\phi^2 w\ dx\lesssim \int r^{-2}\phi^2\ dx\lesssim\int |\partial \phi|^2\ dx\lesssim \int |\partial \phi|^2w\ dx.}

By choosing $T_R\geq 1$ and $\eps\ll1$, we have \fm{\int |\phi|^2q_r^2\lra{q}^{-2}w\ dx&\leq C_R\int |\phi_r|^2w\ dx+C_R\int_{q\geq t^{1/2}} r^{-2}|\phi|^2w\ dx+C_R\int_{q< t^{1/2}} r^{-2}|\phi|^2w\ dx\\&\leq C_R\int |\partial\phi|^2w\ dx+C_Rt^{-1+C\eps}\int \phi^2q_r^2\lra{q}^{-2} w\ dx\\&\leq C_R\int |\partial\phi|^2w\ dx+\frac{1}{2}\int \phi^2q_r^2\lra{q}^{-2} w\ dx.}
This finishes the proof.
\end{proof}

\lem{\label{lp3}Suppose $\phi$ is supported in $|x|-t\geq -R$ and $\phi(t)\in C_c^\infty(\R^3)$ for each $t$. Let $F:=g^{\alpha\beta}\partial_\alpha\partial_\beta u_{app}$ where $u_{app}$ is defined in \eqref{uappdef}. Then for $t\geq T_R\gg1$, we have \fm{\norm{\phi F}_{L^2(w)}\lesssim \eps t^{-1}\norm{\partial\phi}_{L^2(w)}.}}
\begin{proof}Write $F=\frac{\eps}{4r}G(\omega)q_r A_q+F_2$. By the weighted Poincar$\acute{\rm e}$'s lemma, i.e. Lemma \ref{lp2}, we have
\fm{\norm{\eps r^{-1}G(\omega)q_rA_{q}\phi}_{L^2(w)}^2&\lesssim \eps^2(t-R)^{-2}\int q_r^2\lra{q}^{-2}\phi^2w\ dx\lesssim \eps^2t^{-2}\norm{\partial\phi}_{L^2(w)}^2.}

For $|q|\leq R$, $u_{app}=\eps r^{-1}U$. Following the proof of Proposition \ref{prop4}, we can show
\fm{F&=\frac{\eps}{4r}G(\omega)(\mu^2U_{qq}+\mu\mu_qU_q)+O(\eps t^{-2+C\eps})=\frac{\eps}{4r}G(\omega)q_r A_q+O(\eps t^{-2+C\eps}),} so we have $F_2=O(\eps t^{-2+C\eps})$ in $\D$. For $|q|\geq R$, note that $u_{app}$ is only supported for $t\sim r$, we have  $F_2=F=O(\eps t^{-3+C\eps})$ outsides $\D$ by Proposition \ref{mainprop4}.  Since $1\leq w\leq Ct^{C\eps}$, 
\fm{\norm{\phi F_2}_{L^2(w)}^2&=\norm{\phi F_2\chi_{|q|\leq R}}_{L^2(w)}^2+\norm{\phi F_2\chi_{|q|\geq R}}_{L^2(w)}^2\\&\lesssim \int_{|q|\leq R} \eps^2t^{-4+C\eps}|\phi|^2\ dx+\int_{q> R, r\lesssim t} \eps^2t^{-6+C\eps}|\phi|^2\ dx\\&\lesssim \int \eps^2t^{-2}\lra{t-r}^{-2}\phi^2\ dx\\&\lesssim\eps^2t^{-2}\norm{\partial\phi}_{L^2(\R^3)}^2\lesssim\eps^2t^{-2}\norm{\partial\phi}_{L^2(w)}^2.}
Here we use the Poincar$\acute{\rm e}$'s lemma, i.e. Lemma \ref{lp1}. We are done.
\end{proof}

\section{Continuity Argument}\label{smp}\rm

\subsection{Setup}
Fix $\chi(s)\in C_c^\infty(\R)$ such that $\chi\in[0,1]$ for all $s$, $\chi\equiv 1$ for $|s|\leq 1$ and $\chi\equiv 0$ for $|s|\geq 2$. Also fix a large time $T>0$. Consider the equation of $v=v^T(t,x)$
\begin{equation}\label{eqn}
\wt{g}^{\alpha\beta}(u_{app}+v)\partial_\alpha\partial_\beta v=-\chi(t/T)\wt{g}^{\alpha\beta}(u_{app}+v)\partial_\alpha\partial_\beta u_{app},\ t>0;\hspace{1cm} v\equiv 0,\ t\geq 2T.
\end{equation}
We have the following results.
\begin{enumerate}[(a)]
\item By the local existence theory of quasilinear wave equations, we can find a local smooth solution to \eqref{eqn} near $t=2T$.
\item The solution on $[T_1,\infty)$ can be extended to $[T_1-\epsilon,\infty)$ for some small $\epsilon>0$ if  \[\norm{\partial^kv}_{L^\infty([T_1,\infty)\times\R^3)}<\infty,\hspace{1cm}\text{ for all }k\leq 4.\]
\item The solution to \eqref{eqn} has a finite speed of propagation: $v^T(t,x)=0$ if $r+t>6T$ or $r<t-R$, so $Z^I(t/T)=O(1)$ when $T/2\leq t\leq 2T$.
\item If the solution exists for $t\leq T$,  we have  $\wt{g}^{\alpha\beta}(u)\partial_\alpha\partial_\beta u=0$ for $t\leq T$ where  $u=u_{app}+v$.
\end{enumerate}

The proofs of these statements are standard. We refer to \cite{sogg} for the proofs of (a) and (b). In this section, our goal is to prove the following proposition.

\prop{\label{prop6} Fix an integer $N\geq 6$. Then there exist constants $\eps_{N}>0$ which depend on $N$ and $R$, such that for any $0<\eps<\eps_{N}$, \eqref{eqn} has a solution $v=v^T(t,x)$ for all $t\geq0$.  In addition,  $v\equiv 0$ if $r<t-R$; for all $|I|\leq N$, we have
\begin{equation}\label{prop6f1}\norm{\partial Z^Iv(t)}_{L^2(\R^3)}\lesssim_I\eps(1+t)^{-1/2+C_{I}\eps},\hspace{1cm} \forall t\geq 0.\end{equation}
Recall that we choose $R$ based on the support of our scattering data $A(q,\omega)$.
}
\rm

It should be pointed out that the $N$ in this proposition is different from the $N$ in the main theorem.

We  use a continuity argument to prove this proposition. From now on we assume $\eps\ll1$, which means $\eps$ is arbitrary  in $(0,\eps_N)$ for some fixed small constant $\eps_N$ depending on $N$. First we  prove the result for  $t\geq T_{N,R}$, where $T_{N,R}\gg_{N,R}1$ is a sufficiently large constant depending on $N$. We start with a solution $v(t,x)$ for $t\geq T_1$ such that for all $t\geq T_1\geq T_{N,R}$ and $k+i\leq N$,  \begin{equation}\label{ca1}
E_{k,i}(t):=\sum_{l\leq k,|I|\leq i}E_u(\partial^lZ^Iv)(t)\leq B_{k,i}
\eps^2 t^{-1+C_{k,i}\eps},
\end{equation} 
\begin{equation}\label{ca2}
|u|\leq B_0\eps t^{-1+C_{0,2}\eps/2},\ |\partial u|\leq B_1\eps t^{-1}.
\end{equation} 
Here $u:=v+u_{app}$ and $E_u$ is defined in \eqref{energydef}. We remark that $C_{k,i},B_{k,i}$ depend on $k,i$ but not on $N$. Our goal is to prove that \eqref{ca1} and \eqref{ca2} hold with $B_{k,i},B_0,B_1$ replaced by  smaller constants $B_{k,i}^\prime,B_0^\prime,B_1^\prime$, and with $C_{k,i}$ unchanged, assuming that $\eps\ll 1$ and $T_{N,R}\gg 1$. To achieve this goal, we first induct on $i$, and then we induct on $k$ for each fixed $i$. For each $(k,i)$, we want to prove the following inequality
\begin{equation}\label{fff}\begin{aligned}
\sum_{l\leq k,|I|\leq i}\norm{\wt{g}^{\alpha\beta}(u)\partial_\alpha\partial_\beta \partial^lZ^Iv}_{L^2(w)}&\leq C_N\eps t^{-1}E_{k,i}(t)^{1/2}\\&\hspace{1em}+C_N\eps t^{-1+C\eps}(E_{k-1,i}(t)^{1/2}+E_{k+1,i-1}(t)^{1/2})\\&\hspace{1em}+C\eps t^{-3/2+C\eps}.
\end{aligned}\end{equation}
Here $E_{-1,\cdot}=E_{\cdot,-1}=0$, and $C,C_N$ are constants whose meanings will be explained later. We then combine \eqref{fff} with the energy estimates \eqref{energymain} to derive an inequality on $E_{k,i}(t)$. 

We remark that the proof in this section is closely related to that of the energy estimates in Section 9 of Lindblad \cite{lind}.

In the following computation, let $C$ denote a universal constant  or a constant from the previous estimates on $q$ and $u_{app}$ (e.g. from Proposition \ref{mainprop4}). Here $C$ is allowed to depend on $(k,i)$ or $N$, but we will never write it as $C_{k,i}$ or $C_N$. We will choose the constants in the following order:
\fm{&\hspace{1em}C\to C_{0,0},B_{0,0}\to C_{1,0},B_{1,0}\to\dots\to C_{N,0},B_{N,0}\\&\to C_{0,1},B_{0,1}\to\dots\to C_{N-1,1},B_{N-1,1}\\&\to C_{0,2},B_{0,2}\to\dots\to C_{N-2,2},B_{N-2,2}\\&\dots\\&\to C_{0,N},B_{0,N}\\&\to B_0,B_1\to C_N\to T_{N,R}\to\eps.}
In particular, if a constant $A$ appears before a constant $B$, then $A$ cannot depend on $B$.

In addition, since $\eps \ll 1$ and $T_{N,R}\gg 1$ are chosen at the end, we can control terms like $C_N\eps$ and $C_N T_{N,R}^{-\gamma+C_N\eps}$ for $\gamma>0$ for any $k,i$ by a universal constant, e.g. $1$.

To end the setup, we derive a differential equation for $Z^Iv$ from \eqref{eqn}. If we commute \eqref{eqn} with $Z^I$, we have
\begin{equation}\label{eqn2}
\begin{aligned}
&\hspace{1.25em}\wt{g}^{\alpha\beta}(u)\partial_\alpha\partial_\beta Z^Iv\\&=[\square,Z^I]v+[\gamma^{\alpha\beta}(u),Z^I]\partial_\alpha\partial_\beta v+\gamma^{\alpha\beta}(u)[\partial_\alpha\partial_\beta,Z^I]v\\&\hspace{1em}-Z^I(\chi(t/T)(\gamma^{\alpha\beta}(u)-\gamma^{\alpha\beta}(u_{app}))\partial_\alpha\partial_\beta u_{app})-Z^I(\chi(t/T)\wt{g}^{\alpha\beta}(u_{app})\partial_\alpha\partial_\beta u_{app})\\&=:R_1+R_2+R_3+R_4+R_5
\end{aligned}
\end{equation} with $Z^Iv\equiv 0$ for $t\geq 2T$.

\subsection{Pointwise bounds \eqref{ca2}}

In the next few subsections, we always assume $t\geq T_{N,R}\gg 1$. Since $1\leq w\leq Ct^{C\eps}$, by \eqref{ca2} and \eqref{energydef} we have
\begin{equation}\label{equivnorm}C^{-1}\norm{\partial\phi}_{L^2(\R^3)}\leq\norm{\partial\phi}_{L^2(w)}\sim E_u(\phi)^{1/2}\leq Ct^{C\eps}\norm{\partial\phi}_{L^2(\R^3)}.\end{equation}  Here we can choose $\eps\ll 1$ and $T_{N,R}\gg 1$ so that all constants in this inequality are universal.
If we combine this inequality with \eqref{ca1}, we have
\fm{\norm{\partial Z^Iv(t)}^2_{L^2(\R^3)}\leq C E_u(Z^Iv)(t)\leq CB_{0,i} \eps^2t^{-1+C_{0,i}\eps},\hspace{1cm}|I|=i\leq N,}
so by the Klainerman-Sobolev inequality, we have 
\begin{equation}\label{ptb1}|\partial Z^Iv(t)|\leq CB_{0,i+2}^{1/2}\eps t^{-1/2+C_{0,i+2}\eps/2}(1+t+r)^{-1}\lra{t-r}^{-1/2},\hspace{1cm}|I|=i\leq N-2.\end{equation}
Note that 
\fm{\int_{0}^{2t} (1+t+\rho)^{-1}\lra{t-\rho}^{-1/2}\ d\rho&\leq(1+t)^{-1}\int_{0}^{2t} \lra{t-\rho}^{-1/2}\ d\rho\\&\leq 2(1+t)^{-1}\int_{0}^{t} (1+\rho)^{-1/2}\ d\rho\\&\lesssim (1+t)^{-1/2},\\\int_{2t}^\infty(1+t+\rho)^{-1}\lra{t-\rho}^{-1/2}\ d\rho&\lesssim\int_{2t}^\infty(1+\rho)^{-3/2}\ d\rho\lesssim(1+t)^{-1/2}.}
Thus, by integrating $\partial_r Z^Iv(t,\rho\omega)$ from $\rho=t-R$ to $\rho=r$, we conclude that 
\begin{equation}\label{ptb2}|Z^Iv(t)|\leq CB_{0,i+2}^{1/2}\eps t^{-1+C_{0,i+2}\eps/2},\hspace{1cm}|I|=i\leq N-2.\end{equation}
If we let $I=0$ in \eqref{ptb1} and \eqref{ptb2}, we have \fm{|\partial v|\leq CB_{0,2}^{1/2}\eps t^{-3/2+C_{0,2}\eps/2},\hspace{1cm}|v|\leq CB_{0,2}^{1/2}\eps t^{-1+C_{0,2}\eps/2}.} Note that $|u_{app}|\leq C\eps t^{-1+C\eps}$ and $|\partial u_{app}|\leq C\eps t^{-1}$. This allows us to replace $B_0,B_1$ with $B_0/2,B_1/2$ in \eqref{ca2} as long as we choose $T_{N,R},B_0,B_1$ sufficiently large and $\eps$ sufficiently small (e.g. $CB_{0,2}^{1/2}<B_0/4$, $C<B_0/4$; same for $B_1$; $T_{N,R}>10$; $C_{0,2}\eps<1/4$).

In the following computation, we will use \eqref{ptb1} and \eqref{ptb2} directly instead of \eqref{ca2} for the pointwise bounds, so the choice of $C_{k,i},B_{k,i}$ will be independent of $B_0,B_1$.

We remark that if $N\geq 6$, \eqref{ptb1} and \eqref{ptb2} allow us to extend the solution $v(t,x)$ of \eqref{eqn} below $t=T_1$, by the local existence theory of quasilinear wave equations. Moreover, these two pointwise bounds, together with $Z^Iu_{app}=O(\eps(1+t)^{-1+C\eps})$, allow us to use Lemma \ref{l24} freely, as long as $\eps\ll 1$ and $T_{N,R}\gg 1$.

\subsection{Energy estimate \eqref{ca1} with $k=i=0$}
Let $k=i=0$ and fix $T_1\leq t\leq 2T$. Now $R_1=R_2=R_3=0$ in \eqref{eqn2}. 

For $R_4$, since $|\chi(t/T)|\leq 1$, we have \fm{\norm{R_4}_{L^2(w)}&\leq \norm{g^{\alpha\beta}v\partial_\alpha\partial_\beta u_{app}}_{L^2(w)}+C\norm{|v|(|u_{app}|+|v|)|\partial^2u_{app}|}_{L^2(w)}\\&\leq C\eps t^{-1}\norm{\partial v}_{L^2(w)}+C_N\eps^2 t^{-2+C_N\eps}\norm{|v|\lra{t-r}^{-1}}_{L^2(w)}\\&\leq C\eps t^{-1}E_u(v)(t)^{1/2}.}Here we apply Lemma \ref{l24} in the first inequality, Lemma \ref{lp3} in the second inequality,  Lemma \ref{lp1} and \eqref{equivnorm} in the third inequality.

For $R_5$, since $u_{app}$ is supported in the ball centered at origin with radius $2t$, by Proposition \ref{mainprop4} we have \fm{\norm{R_5}_{L^2(w)}\leq C \eps t^{-3/2+C\eps}.}

Thus, by \eqref{energymain}, we conclude that
\fm{E_u(v)(t)&\leq\int_t^{2T}C_N\eps \tau^{-1}E_u(v)(\tau)+C\eps\tau^{-3/2+C\eps}E_u(v)(\tau)^{1/2}\ d\tau\\&\leq \int_t^{2T}C_NB_{0,0}\eps^3\tau^{-2+C_{0,0}\eps}+CB_{0,0}^{1/2}\eps^2\tau^{-2+(C+C_{0,0}/2)\eps}\ d\tau\\&\leq CC_NB_{0,0}\eps^3t^{-1+C_{0,0}\eps}+CB_{0,0}^{1/2}\eps^2t^{-1+(C+C_{0,0}/2)\eps}.}
In particular, the constants $C$ do not depend on $C_N$ or $C_{k,i},B_{k,i}$ in \eqref{ca1}. If $\eps\ll1$ (say $CC_N\eps\leq 1/4$) and $C_{0,0},B_{0,0}$ are large enough (say $C_{0,0}/2+C<C_{0,0},C\sqrt{B_{0,0}}< B_{0,0}/4$), we obtain \eqref{ca1} with $B_{0,0}$ replaced by $B_{0,0}/2$.\par

\subsection{Energy estimate \eqref{ca1} with $i=0$ and $k>0$}
Let $i=0$ and $k>0$ and fix $T_1\leq t\leq 2T$. Now $R_1=R_3=0$.

For $R_2$, we have
\fm{\norm{R_2}_{L^2(w)}&\leq \norm{[g^{\alpha\beta}u,\partial^k]\partial_\alpha\partial_\beta v}_{L^2(w)}+\norm{[\gamma^{\alpha\beta}(u)-g^{\alpha\beta}u,\partial^k]\partial_\alpha\partial_\beta v}_{L^2(w)}
\\&\leq C\sum_{k_1+k_2\leq k, k_1>0}\norm{|\partial^{k_1}u||\partial^{k_2+2}v|}_{L^2(w)}\\&\hspace{1em}+C\sum_{k_1+k_2+k_3\leq k, k_3<k}\norm{|\partial^{k_1}u||\partial^{k_2}u||\partial^{k_3+2}v|}_{L^2(w)}.}
The second sum comes from Lemma \ref{l24}. By writing $u=v+u_{app}$, we have the following terms in the sums:
\fm{&\norm{|\partial u_{app}||\partial^{k_2+2}v|}_{L^2(w)}\leq C \eps t^{-1}E_{k,0}(t)^{1/2},\hspace{1cm} k_2<k;\\
&\norm{|\partial^{k_1} u_{app}||\partial^{k_2+2}v|}_{L^2(w)}\leq C \eps t^{-1+C\eps}E_{k-1,0}(t)^{1/2},\hspace{1cm}k_1+k_2\leq k,k_1>1;\\
&\norm{|\partial^{k_1} v||\partial^{k_2+2}v|}_{L^2(w)}\leq C_N \eps t^{-3/2+C_N\eps}E_{k,0}(t)^{1/2},\hspace{1cm}k_1+k_2\leq k,k_1>0;\\
&\norm{|\partial^{k_1} u_{app}||\partial^{k_2} u_{app}||\partial^{k_3+2}v|}_{L^2(w)}\leq C \eps^2 t^{-2+C\eps}E_{k,0}(t)^{1/2},\hspace{1cm}k_1+k_2+k_3\leq k,k_3<k;\\
&\norm{|\partial^{k_1} u_{app}||\partial^{k_2} v||\partial^{k_3+2}v|}_{L^2(w)}\leq C_N \eps^2 t^{-2+C_N\eps}E_{k,0}(t)^{1/2},\hspace{1cm}k_1+k_2+k_3\leq k,k_3<k;\\
&\norm{|\partial^{k_1} v||\partial^{k_2} v||\partial^{k_3+2}v|}_{L^2(w)}\leq C_N \eps^2 t^{-2+C_N\eps}E_{k,0}(t)^{1/2},\hspace{1cm}k_1+k_2+k_3\leq k,k_3<k.} 
Here we use Proposition \ref{mainprop4}, \eqref{ptb1} and \eqref{ptb2}. We take $L^2(w)$ norm on the derivative of $v$ with the highest order, and apply pointwise bounds on the derivatives of $u_{app}$ or derivatives of $v$ with lower orders. Here we need $N/2+1\leq N-2$, i.e. $N\geq 6$, to apply the pointwise bounds. Thus, we have 
\fm{\norm{R_2}_{L^2(w)}\leq C\eps t^{-1}E_{k,0}(t)^{1/2}+C\eps t^{-1+C\eps} E_{k-1,0}(t)^{1/2}.}The constants here are universal, as long as we choose $\eps\ll1$ (say $C_N\eps< 1$) and $T_{N,R}$ sufficiently large (say $C_N/\sqrt{T_{N,R}}\leq 1$). 

For $R_4$, since $\partial^l(\chi(t/T))=O(1)$ for all $l$, by Lemma \ref{l24} we have 
\fm{\norm{R_4}_{L^2(w)}&\leq C\sum_{k_1\leq k}\norm{g^{\alpha\beta}\partial^{k_1}v\partial_\alpha\partial_\beta u_{app}}_{L^2(w)}+C\sum_{k_1+k_2\leq k,\ k_2> 0}\norm{|\partial^{k_1}v||\partial^{k_2+2}u_{app}|}_{L^2(w)}\\&\hspace{1em}+C\sum_{k_1+k_2+k_3\leq k}\norm{|\partial^{k_1}v|(|\partial^{k_2}u_{app}|+|\partial^{k_2}v|)|\partial^{k_3+2}u_{app}|}_{L^2(w)}.}
By Lemma \ref{lp3}, the first sum has an upper bound\fm{C\eps t^{-1}\sum_{k_1\leq k}\norm{\partial\partial^{k_1}v}_{L^2(w)}\leq C\eps t^{-1} E_{k,0}(t)^{1/2}.}
By Lemma \ref{lp1}, the second sum has an upper bound
\fm{C\eps t^{-1+C\eps}\sum_{k_1<k}\norm{\partial^{k_1}v\lra{t-r}^{-2}}_{L^2(w)}&\leq C\eps t^{-1+C\eps}\sum_{k_1<k}\norm{\partial\partial^{k_1}v}_{L^2(\R^3)}\\&\leq C\eps t^{-1+C\eps} E_{k-1,0}(t)^{1/2}. }
The third sum is controlled by the second one, because $|\partial^{k_2}u_{app}|\leq C\eps t^{-1+C\eps}\leq 1$, and at least one of $|\partial^{k_1}v|$ and $|\partial^{k_2}v|$ is $\leq C_N\eps t^{-1+C_N\eps}\leq 1$ (since $\min\{k_1,k_2\}\leq k/2\leq N-2$). In conclusion, \fm{\norm{R_4}_{L^2(w)}\leq C\eps t^{-1}E_{k,0}(t)^{1/2}+C\eps t^{-1+C\eps} E_{k-1,0}(t)^{1/2}.}The constants here are again universal.

For $R_5$, we have 
\fm{\norm{R_5}_{L^2(w)}\leq C\eps t^{-3/2+C\eps}.}

Thus, by \eqref{energymain}, we have
\fm{E_{k,0}(t)&\leq\int_t^{2T} C_N\eps(1+\tau)^{-1}E_{k,0}(\tau)+C_N\eps\tau^{-1+C\eps}E_{k-1,0}(\tau)^{1/2}E_{k,0}(\tau)^{1/2}\\&\hspace{1em}+ C\eps\tau^{-3/2+C\eps}E_{k,0}(\tau)^{1/2}d\tau\\&\leq\int_t^{2T} C_NB_{k,0}\eps^3\tau^{-2+C_{k,0}\eps}+C_NB_{k,0}\eps^3\tau^{-2+(C+C_{k,0}/2+C_{k-1,0}/2)\eps}\\&\hspace{1em}+ CB_{k,0}^{1/2}\eps^2\tau^{-2+(C+C_{k,0}/2)\eps}\ d\tau\\&\leq CC_NB_{k,0}\eps^3 t^{-1+C_{k,0}\eps}+CC_NB_{k,0}\eps^3 t^{-1+(C+C_{k,0}/2+C_{k-1,0}/2)\eps}\\&\hspace{1em}+CB_{k,0}^{1/2}\eps^2 t^{-1+(C+C_{k,0}/2)\eps}.}
Similarly we can  prove \eqref{ca1} with $B_{k,0}$ replaced by $B_{k,0}/2$, if we assume that $B_{k,0},C_{k,0}$ are large enough and $\eps\ll1$ (say $CC_N\eps<1/8$, $C_{k,0}\geq C_{k-1,0}$, $C\sqrt{B_{k,0}}\leq B_{k,0}/8$).

\subsection{Energy estimate \eqref{ca1} with $k=0$ and $i>0$}
Let $k=0$ and $i>0$ and fix $T_1\leq t\leq 2T$. Also fix $Z^I$ with $|I|= i$.

For $R_2$, we have 
\fm{\norm{R_2}_{L^2(w)}&\leq \norm{[g^{\alpha\beta}u,Z^I]\partial_\alpha\partial_\beta v}_{L^2(w)}+\norm{[\gamma^{\alpha\beta}(u)-g^{\alpha\beta}u,Z^I]\partial_\alpha\partial_\beta v}_{L^2(w)}\\&\leq C\sum_{|J_1|+|J_2|\leq i,\ |J_1|>0}\norm{|Z^{J_1}u||\partial^2 Z^{J_2}v|}_{L^2(w)}\\&\hspace{1em}+C\sum_{|J_1|+|J_2|+|J_3|\leq i,\ |J_3|<i}\norm{Z^{J_1}uZ^{J_2}u\partial^2Z^{J_3}v}_{L^2(w)}.}
The second sum comes from Lemma \ref{l24}. Note that the second sum is controlled by the first sum. In fact, since $|J_1|,|J_2|$ cannot be greater than $i/2$ at the same time, without loss of generality we assume $|J_1|\leq i/2\leq N-2$. Thus $|Z^{J_1}u|\leq C_N\eps t^{-1+C_N\eps}\leq 1$ by \eqref{ptb2} if we choose $\eps\ll 1$. For the first sum, by writing $u=v+u_{app}$, we have the following terms in the sum:
\fm{
&\norm{|Z^{J_1} u_{app}||\partial^{2}Z^{J_2}v|}_{L^2(w)},\hspace{1cm}|J_1|+|J_2|\leq i,|J_1|>0;\\
&\norm{|Z^{J_1} v||\partial^{2}Z^{J_2}v|}_{L^2(w)},\hspace{1cm}|J_1|+|J_2|\leq i,|J_1|>0.
} 
The first term has an upper bound
\fm{C\eps t^{-1+C\eps}E_{1,i-1}(t)^{1/2}.}
By Lemma \ref{l21}, we can see that the second term is controlled by 
\fm{C\norm{|\lra{t-r}^{-1}Z^{J_1} v||\partial ZZ^{J_2}v|}_{L^2(w)},\hspace{1cm}|J_2|<i.}
If $|J_1|\leq N-2$, then by \eqref{ptb1} we have
\fm{|\lra{t-r}^{-1}Z^{J_1}v|\leq \lra{t-r}^{-1}\int_{t-R}^r |\partial_\rho Z^{J_1}(t,\rho\omega)|\ d\rho\leq C\norm{\partial Z^{J_1}v(t)}_{L^\infty(\R^3)}\leq C_N\eps t^{-3/2+C_{N}\eps},}
which implies that 
\fm{C\norm{|\lra{t-r}^{-1}Z^{J_1} v||\partial ZZ^{J_2}v|}_{L^2(w)}\leq CC_N t^{-3/2+C_N\eps}E_{0,i}(t)^{1/2}.}
If $|J_1|\geq N-1$, then $|J_2|\leq 1$. In this case, by \eqref{ptb1}, \eqref{equivnorm} and Lemma \ref{lp1}, we have
\fm{|\partial ZZ^{J_2}v|\leq C_N\eps t^{-3/2+C_N\eps},}
\fm{\norm{\lra{t-r}^{-1} Z^{J_1}v}_{L^2(w)}&\leq C t^{C\eps}\norm{\lra{t-r}^{-1}Z^{J_1}v}_{L^2(\R^3)}\leq C t^{C\eps}\norm{\partial Z^{J_1}v}_{L^2(\R^3)}\\&\leq C t^{C\eps}E_{0,i}(t)^{1/2}.}
Thus, the term above is controlled by \fm{CC_N\eps t^{-3/2+(C_N+C)\eps}E_{0,i}(t)^{1/2}.}

For $R_3$, following the same discussion as above,  we have
\fm{\norm{R_3}_{L^2(w)}&\leq C\sum_{|J|<i}\norm{|u||\partial^2 Z^Jv|}_{L^2(w)}\\&\leq C\eps t^{-1+C\eps}\sum_{|J|<i}\norm{\partial^2 Z^Jv}_{L^2(w)}+C\sum_{|J|<i}\norm{|v||\partial^2 Z^Jv|}_{L^2(w)}\\&\leq C\eps t^{-1+C\eps}E_{1,i-1}(t)^{1/2}+CC_N\eps t^{-3/2+(C_N+C)\eps}E_{0,i}(t)^{1/2}.}

For $R_4$, since $Z^J(\chi(t/T))=O(1)$ for all $J$ by finite speed of propagation,  we have
\fm{\norm{R_4}_{L^2(w)}&\leq C\sum_{|J|\leq i}\norm{g^{\alpha\beta}Z^Jv\partial_\alpha\partial_\beta u_{app}}_{L^2(w)}+C\sum_{|J_1|+|J_2|\leq i,|J_2|>0}\norm{|Z^{J_1}v||\partial^2Z^{J_2}u_{app}|}_{L^2(w)}
\\&\hspace{1em}+C\sum_{|J_1|+|J_2|+|J_3|\leq i}\norm{|Z^{J_1}v|(|Z^{J_2}v|+|Z^{J_2}u_{app}|)|\partial^2Z^{J_3}u_{app}|}_{L^2(w)}
\\&\leq C\eps t^{-1}E_{0,i}(t)^{1/2}+C\eps t^{-1+C\eps}E_{0,i-1}(t)^{1/2}.}The proof is very similar to the proof on estiamte of $R_4$ in the case $i=0$ and $k>0$.

For $R_5$, again we have 
\fm{\norm{R_5}_{L^2(w)}\leq C\eps t^{-3/2+C\eps}.}

For $R_1$,  we have
\fm{|[\square, Z^I]v|&\lesssim\sum_{|J_1|+|J_2|<i} |Z^{J_1}\square Z^{J_2}v|\lesssim\sum_{|J|<i}|Z^J\square v|\\&\lesssim\sum_{|J|<i} |Z^{J}(\wt{g}^{\alpha\beta}(u)\partial_\alpha\partial_\beta v)|+|Z^{J}(\gamma^{\alpha\beta}(u)\partial_\alpha\partial_\beta v)|\\&\lesssim\sum_{|J|<i} |Z^{J}(\chi(t/T)\wt{g}^{\alpha\beta}(u)\partial_\alpha\partial_\beta u_{app})|+|Z^{J}(\gamma^{\alpha\beta}(u)\partial_\alpha\partial_\beta v)|.}
Here all the constants are universal which depend only on $i,N$. The first term is simply $R_4+R_5$ with a lower order $I$. The second term can be controlled in the same way as we control~$R_2,R_3$. In conclusion, 
\fm{E_{0,i}(t)&\leq\int_{t}^{2T} CC_N\eps \tau^{-1}E_{0,i}(\tau)+CC_N\eps \tau^{-1+C\eps}E_{1,i-1}(\tau)^{1/2}E_{0,i}(\tau)^{1/2}\\&\hspace{1em}+C\eps \tau^{-3/2+C\eps}E_{0,i}(\tau)^{1/2}\ d\tau\\&\leq\int_{t}^{2T} CC_NB_{0,i}\eps^3 \tau^{-2+C_{0,i}\eps}+CC_NB_{0,i}\eps^3 \tau^{-2+(C+C_{1,i-1}/2+C_{0,i}/2)\eps}\\&\hspace{1em}+CB_{0,i}^{1/2}\eps^2 \tau^{-2+(C+C_{0,i}/2)\eps}\ d\tau\\&\leq CC_NB_{0,i}\eps^3 t^{-1+C_{0,i}\eps}+CC_NB_{0,i}\eps^3 t^{-1+(C+C_{0,i}/2+C_{1,i-1}/2)\eps}\\&\hspace{1em}+CB_{0,i}^{1/2}\eps^2 t^{-1+(C+C_{0,i}/2)\eps}.}Again, we can choose $B_{0,i},C_{0,i}$ sufficiently large such that \eqref{ca1} holds with $B_{0,i}$ replaced by $B_{0,i}/2$. Note that $B_{1,i-1},C_{1,i-1}$ are already chosen when we consider the case $k=0,i>0$.

\subsection{Energy estimate \eqref{ca1} with $k,i>0$}
Let $k,i>0$ and fix $T_1\leq t\leq 2T$. Also fix $Z^I$ with $|I|= i$. This  case can be viewed as a combination of the case $k=0,i>0$ and the case $i=0,k>0$.

For $R_2$, we have 
\fm{\norm{R_2}_{L^2(w)}&\leq \norm{[g^{\alpha\beta}u,\partial^kZ^I]\partial_\alpha\partial_\beta v}_{L^2(w)}+\norm{[\gamma^{\alpha\beta}(u)-g^{\alpha\beta}u,\partial^kZ^I]\partial_\alpha\partial_\beta v}_{L^2(w)}
\\&\leq C\sum_{k_1+k_2\leq k, |J_1|+|J_2|\leq i, k_1+|J_1|>0}\norm{|\partial^{k_1}Z^{J_1}u||\partial^{2+k_2} Z^{J_2}v|}_{L^2(w)}\\&\hspace{1em}+C\sum_{\tiny \begin{array}{l}k_1+k_2+k_3\leq k\\|J_1|+|J_2|+|J_3|\leq i\\k_3+|J_3|<k+i\end{array}}\norm{\partial^{k_1}Z^{J_1}u\partial^{k_2}Z^{J_2}u\partial^{2+k_3}Z^{J_3}v}_{L^2(w)}.}
The second sum is again easy to handle. For the first sum, we consider the following three cases: $k_1=0$ and $|J_1|>0$; $k_1=1$ and $|J_1|=0$; all the remaining choices of $(k,J_1)$. For the first case, we apply Proposition \ref{mainprop4} and Lemma \ref{l21} to obtain a factor $\lra{t-r}^{-1}$ with one $\partial$ replaced by $Z$; for the second case, we use $|\partial u|\leq C_N\eps t^{-1}$; for the third, we use \eqref{ptb1} directly. The proof here is very similar to the proof in the previous cases. We thus have 
\fm{\norm{R_2}_{L^2(w)}\leq C_N\eps t^{-1}E_{k,i}(t)^{1/2}+C_N\eps t^{-1+C\eps}(E_{k-1,i}(t)^{1/2}+E_{k+1,i-1}(t)^{1/2}).}

For $R_3$, we have 
\fm{\norm{R_3}_{L^2(w)}\leq C\sum_{k_1\leq k,|J|<i}\norm{|u||\partial^{k_1+2}Z^Jv|}_{L^2(w)}.}
We can use Proposition \ref{mainprop4} and Lemma \ref{l21} to obtain \fm{\norm{R_3}_{L^2(w)}\leq C_N\eps t^{-1}E_{k,i}(t)^{1/2}+C_N\eps t^{-1+C\eps}(E_{k-1,i}(t)^{1/2}+E_{k+1,i-1}(t)^{1/2}).}

For $R_4$, we have 
\fm{\norm{R_4}_{L^2(w)}&\leq C\sum_{k_1\leq k,|J|\leq i}\norm{g^{\alpha\beta}\partial^{k_1}Z^Jv\partial_\alpha\partial_\beta u_{app}}_{L^2(w)}\\&\hspace{1em}+C\sum_{k_1+k_2\leq k,|J_1|+|J_2|\leq i,k_2+|J_2|>0}\norm{|\partial^{k_1}Z^{J_1}v||\partial^{k_2+2}Z^{J_2}u_{app}|}_{L^2(w)}
\\&\hspace{1em}+C\sum_{\tiny \begin{array}{l}k_1+k_2+k_3\leq k\\|J_1|+|J_2|+|J_3|\leq i\\k_3+|J_3|<k+i\end{array}}\norm{|\partial^{k_1}Z^{J_1}v|(|\partial^{k_2}Z^{J_2}v|+|\partial^{k_2}Z^{J_2}u_{app}|)|\partial^{k_3+2}Z^{J_3}u_{app}|}_{L^2(w)}
\\&\leq C_N\eps t^{-1}E_{k,i}(t)^{1/2}+C_N\eps t^{-1+C\eps}(E_{k-1,i}(t)^{1/2}+E_{k+1,i-1}(t)^{1/2}).}This can be handled in the same way as we handle $R_4$ in the case $k=0,i>0$ or $i=0,k>0$.

For $R_5$, again we have 
\fm{\norm{R_5}_{L^2(w)}\leq C\eps t^{-3/2+C\eps}.}

For $R_1$, since $[\square,\partial^kZ^I]=\partial^k[\square,Z^I]$, we can conclude that the $L^2(w)$ norm of $R_1$ can be controlled by the bounds of the $L^2(w)$ norms of all other $R_i$.

In conclusion, we have
\fm{E_{k,i}(t)&\leq\int_{t}^{2T} C C_N\eps \tau^{-1}E_{k,i}(\tau)+CC_N\eps \tau^{-1+C\eps}(E_{k-1,i}(\tau)^{1/2}+E_{k+1,i-1}(\tau)^{1/2})E_{k,i}(\tau)^{1/2}\\&\hspace{1em}+C\eps \tau^{-3/2+C\eps}E_{k,i}(\tau)^{1/2}\ d\tau\\&\leq\int_{t}^{2T} C C_NB_{k,i}\eps^3 \tau^{-2+C_{k,i}\eps}+CC_NB_{k,i}\eps^3 \tau^{-2+(C+C_{k+1,i-1}/2+C_{k-1,i}/2+C_{k,i}/2)\eps}\\&\hspace{1em}+CB_{k,i}^{1/2}\eps^2 \tau^{-2+(C+C_{k,i}/2)\eps}\ d\tau\\&\leq C C_NB_{k,i}\eps^3 t^{-1+C_{k,i}\eps}+CC_NB_{k,i}\eps^3 t^{-1+(C+C_{k+1,i-1}/2+C_{k-1,i}/2+C_{k,i}/2)\eps}\\&\hspace{1em}+CB_{k,i}^{1/2}\eps^2 t^{-1+(C+C_{k,i}/2)}.}
Again, we can choose $B_{k,i},C_{k,i}$ sufficiently large such that \eqref{ca1} holds with $B_{k,i}$ replaced by $B_{k,i}/2$. Note that $B_{k+1,i-1},C_{k+1,i-1},B_{k-1,i},C_{k-1,i}$ are already chosen when we consider the case $k,i>0$.

\subsection{Existence for $0\leq t\leq T_{N,R}$}
In the previous subsections, we prove that there exists a solution $v$ to \eqref{eqn} for all $t\geq T_{N,R}$ with \eqref{prop6f1} hold for all $|I|\leq N$ and $t\geq T_{N,R}$. Now we finish the proof of Proposition \ref{prop6} by extending the solution to all $t\geq 0$. At a small time, $u_{app}$ does not approximate $u$ well, but $u_{app}$ and all its derivatives stay bounded for all $(t,x)$ with $0\leq t\leq T_{N,R}$. See Proposition \ref{mainprop4}. So, it is better to use \eqref{qwe} to control $u$ directly instead of using \eqref{eqn}.

Fix $N\geq 6$. By using the pointwise bounds in Proposition \ref{mainprop4} and the support of $u_{app}$, we  have
\fm{\norm{Z^Iu_{app}(t)}_{L^2(\R^3)}\lesssim_{I,N,R}\eps,\hspace{1cm}  0\leq t\leq T_{N,R}.} Thus, it suffices to prove that the solution $u$ to \eqref{qwe} with $u=v	+u_{app}$ for $t\geq T_{N,R}$ exists for $0\leq t\leq T_{N,R}$, with \fm{\norm{\partial Z^Iu(t)}_{L^2(\R^3)}\lesssim_{I,N,R}\eps,\hspace{1cm}0\leq t\leq T_{N,R},\ |I|\leq N.}

If we apply $Z^I$ to \eqref{qwe}, we have
\begin{equation}\label{eqn3}\begin{aligned}
\wt{g}^{\alpha\beta}(u)\partial_\alpha\partial_\beta Z^Iu&=[\square,Z^I]u+[\gamma^{\alpha\beta}(u),Z^I]\partial_\alpha\partial_\beta u+\gamma^{\alpha\beta}(u)[\partial_\alpha\partial_\beta, Z^I]u.
\end{aligned}\end{equation}

We can now set up the continuity argument. Suppose that we have a solution $u$ to \eqref{qwe} for $T_1\leq t\leq T_{N,R}$ for some $0\leq T_1\leq T_{N,R}$, such that 
\begin{equation}\label{ea}
\norm{\partial Z^Iu(t)}_{L^2(\R^3)}\leq B\eps,\hspace{1cm}|I|\leq N,T_1\leq t\leq T_{N,R}.
\end{equation} Here $B=B_N$ depends on $N$. We remark that \eqref{ea} implies \eqref{prop6f1}  for $t\leq T_{N,R}$, where the power is the same but the constant in $\lesssim_I$ now depends on $N$. This is because $1\lesssim_N t^{-1/2+C_I\eps}$ for $t\leq T_{N,R}$, assuming $\eps\ll1$.

By the Klainerman-Sobolev inequality, we conclude that  for $t\geq T_1$
\fm{|\partial Z^I u(t,x)|\leq CB\eps(1+t+r)^{-1}\lra{t-r}^{-1/2},\hspace{1cm} |I|\leq N-2}
and \fm{|Z^I u(t,x)|\leq CB\eps(1+t)^{-1/2},\hspace{1cm} |I|\leq N-2.}
The proof of the second estimate is similar to that of \eqref{ptb2}. Thus, assuming $\eps\ll1$, from \eqref{eqn3}  we have for $|I|\leq N$
\fm{|\wt{g}^{\alpha\beta}(u)\partial_\alpha\partial_\beta Z^Iu|&\leq C\sum_{|J|+|K|\leq |I|,|K|<|I|}|Z^Ju||\partial^2 Z^Ku|\\&\leq C\sum_{|J|+|K|\leq |I|,|K|< |I|}\lra{t-r}^{-1}|Z^Ju||\partial ZZ^Ku|\\&\leq C_N\eps \sum_{|J|\leq |I|}(|\partial Z^Ju|+\lra{t-r}^{-1}|Z^Ju|).}
Here we apply Lemma \ref{l24} in the first inequality and the pointwise bounds in the third one. Note that if $|J|+|K|\leq |I|$ and $|K|<|I|$, then $\min\{|J|,|K|+1\}\leq N/2+1\leq N-2$ when $N\geq 6$.

Now we can use the standard energy estimates, say Proposition 2.1 in Chapter I in Sogge \cite{sogg} or Proposition 6.3.2 in H\"{o}rmander \cite{horm}. We apply the Poincar$\acute{\rm e}$'s lemma, i.e.\ Lemma \ref{lp1}, to $\lra{t-r}^{-1}|Z^Ju|$, so its $L^2(\R^3)$ norm is controlled by the that of $|\partial Z^Ju|$. By setting \fm{E_N(t)=\sum_{|I|\leq N}\norm{\partial Z^Iv(t)}^2_{L^2(\R^3)},} for  small $\eps\ll1$, we have 
\fm{E_N(t)^{1/2}&\leq 2(E_N(T_{N,R})^{1/2}+C_N\eps\int_t^{T_{N,R}}E_N(\tau)^{1/2}\ d\tau)\exp(\int_t^{T_{N,R}}C_N\eps\ d\tau)\\&\leq C_N\eps +C_N B^{1/2}\eps^2.}
Then by choosing $\eps$ small enough and $B$ large enough, both depending on $N$, we can replace $B$ with $B/2$ in \eqref{ea}. We are done.

Finally, we remark that for each $|I|\leq N$ and $\eps\ll1$, we can apply Proposition \ref{prop6} with $N$ replaced by $N^\prime=\max\{6,|I|\}\leq N$. Note that when $\eps<\eps_N\leq \eps_{N^\prime}$ and $T>T_{N,R}\geq T_{N^\prime,R}$, the solution for $N$ and the solution for $N^\prime$ are exactly the same. But the constants in \eqref{prop6f1} now depend on $\max\{6,|I|\}$ instead of $N$. This allows us to remove the dependence of $N$ in the coefficients of \eqref{prop6f1}.

\section{Limit as $T\to\infty$}\label{slim}

Our goal for this section is to prove the following proposition. 

\prop{\label{prop7} Fix $N\geq 6$. Then for the same $\eps_{N}$ in Proposition \ref{prop6} and for  $0<\eps<\eps_{N}$, there is a solution $u$ to \eqref{qwe}  in $C^{N-4}$ for all $t\geq 0$, such that for all  $|I|\leq N-5$
\begin{equation}\label{p7f1}
\norm{\partial Z^I(u-u_{app})(t)}_{L^2(\R^3)}\lesssim_I\eps(1+t)^{-1/2+C_I\eps},\hspace{1cm}t\geq 0.
\end{equation}
Besides, for all $|I|\leq N-5$ and $t\gg_R1$,
\begin{equation}\label{p7f2}
|\partial Z^I(u-u_{app})(t,x)|\lesssim_I\eps t^{-1/2+C_I\eps}\lra{r+t}^{-1}\lra{t-r}^{-1/2},
\end{equation}
\begin{equation}\label{p7f3}
|Z^I(u-u_{app})(t,x)|\lesssim_I\min\{\eps t^{-1+C_I\eps},\eps t^{-3/2+C_I\eps}\lra{r-t}\}.
\end{equation}
}\rm

It should be pointed out that the value of ``$N$'' in the main theorem is equal to $N-4$ for the $N$ in this proposition.

From now on, the constant $C$ is allowed to depend on all the constants in the previous sections (say $C_{k,i},B_{k,i},N$), but it must be independent of $\eps$ and $T$.

\subsection{Existence of the limit}\label{sec7.1} Fix $N\geq 6$ and $T_2>T_1\gg 1$. By Proposition \ref{prop6}, for each $0<\eps<\eps_{N}$, we  get two corresponding solutions $v_1=v^{T_1}$ and $v_2=v^{T_2}$ which exist for all $t\geq 0$. Our goal now is to prove that $v_1-v_2$ tends to $0$ in some Banach space as $T_2>T_1\to\infty$.

Recall that $\eps_{N},T_{N,R}$ are independent of the choice of $T$, as long as $T>T_{N,R}$. In addition, $v_1$ and $v_2$ satisfy \eqref{prop6f1}, \eqref{ca1}, \eqref{ca2}, \eqref{ptb1} and \eqref{ptb2}, as shown in the continuity argument, for $t\geq T_{N,R}$, and they satisfy \eqref{ea} along with the pointwise bounds for $0\leq t\leq T_{N,R}$. All the constants involved in these estimates are independent of $T$.  We define  $u_{1}=v^{T_1}+u_{app}$, $u_{2}=v^{T_2}+u_{app}$ and $\wt{v}=v^{T_2}-v^{T_1}$. Then, for $t\geq T_1$ and $|I|\leq N$, by \eqref{prop6f1}, we have \fm{\norm{\partial Z^I\wt{v}(t)}_{L^2(\R^3)}&\leq\norm{\partial Z^Iv_1(t)}_{L^2(\R^3)}+\norm{\partial Z^Iv_2(t)}_{L^2(\R^3)}\leq C\eps^2T_1^{-1+C\eps}.}
In addition, for $t\leq T_1$ (now $\chi(t/T_1)=\chi(t/T_2)=1$) and for each $|I|\leq N$, we have
\begin{equation}\label{teqn}\begin{aligned}\wt{g}^{\alpha\beta}(u_1)\partial_\alpha\partial_\beta Z^I\wt{v}&=[\square,Z^I]\wt{v}+[\gamma^{\alpha\beta}(u_1),Z^I]\partial_\alpha\partial_\beta \wt{v}+[\gamma^{\alpha\beta}(u_2)-\gamma^{\alpha\beta}(u_1),Z^I]\partial_\alpha\partial_\beta v_2\\&\hspace{1em}+\gamma^{\alpha\beta}(u_1)[\partial_\alpha\partial_\beta,Z^I]\wt{v}+(\gamma^{\alpha\beta}(u_2)-\gamma^{\alpha\beta}(u_1))[\partial_\alpha\partial_\beta,Z^I]v_2\\&\hspace{1em}-Z^I((\gamma^{\alpha\beta}(u_2)-\gamma^{\alpha\beta}(u_1))\partial_\alpha\partial_\beta u_{app})-(\gamma^{\alpha\beta}(u_2)-\gamma^{\alpha\beta}(u_1))\partial_\alpha\partial_\beta Z^Iv_2.\end{aligned}\end{equation}

Define a new energy \fm{\wt{E}_{k,i}(t):=\sum_{l\leq k,|I|\leq i} E_{u_1}(\partial^lZ^I\wt{v})(t).}Here $E_{u_1}$ is defined in \eqref{energydef} with $u$ replaced by $u_1$. For $k+i\leq N-3$ with $|I|=i$, and for $t\geq T_{N,R}$ we have 
\begin{equation}\label{f71}\norm{\wt{g}^{\alpha\beta}(u_1)\partial_\alpha\partial_\beta \partial^kZ^I\wt{v}}_{L^2(w)}\leq C\eps t^{-1}\wt{E}_{k,i}(t)^{1/2}+C\eps t^{-1+C\eps}(\wt{E}_{k-1,i}(t)^{1/2}+\wt{E}_{k+1,i-1}(t)^{1/2})\end{equation}with $\wt{E}_{-1,\cdot}=\wt{E}_{\cdot,-1}=0$. This is a simple application of Lemma \ref{l21}, Lemma \ref{l24} and the estimates for $u_1,v_1,u_2,v_2$.  We skip the detail of the proof here, since it is very similar to the proof of \eqref{fff} on $E_{k,i}$. However, we should always put $L^2(w)$ norm on the terms involving $\wt{v}$ and put $L^\infty$ norm on terms involving $u_1,u_2,v_1,v_2$. The pointwise bounds only holds for $|I|\leq N-2$, as seen in \eqref{ptb1} and \eqref{ptb2}, so we need to assume $k+i\leq N-3$ instead of $k+i\leq N$ above. Besides, there is no term like $R_5$ in the previous section, so we expect $\wt{E}_{k,i}$ to have a better decay than $E_{k,i}$.

Since \eqref{ptb1} and \eqref{ptb2} hold for $v_1$, we can apply energy estimate \eqref{energymain}  for  $E_{u_1}$.  Thus, for all $T_{N,R}\leq t\leq T_1$ and for $k+i\leq N-3$, 
\fm{\wt{E}_{k,i}(t)&\leq C\eps^2 T_1^{-1+C\eps}+B\int_t^{T_1}\eps\tau^{-1}\wt{E}_{k,i}(\tau)\ d\tau\\&\hspace{1em}+C\eps\int_t^{T_1}\tau^{-1+C\eps}(\wt{E}_{k-1,i}(\tau)^{1/2}+\wt{E}_{k+1,i-1}(\tau)^{1/2})\wt{E}_{k,i}(\tau)^{1/2}\ d\tau.}
Using this estimate, we claim that $\wt{E}_{k,i}(t)\leq C\eps^2T_1^{-1+C\eps}$ for all $k+i\leq N-3$.  Here $C$ may depend on $k,i$. To prove this claim, we first induct on $i=0,1,\dots,N$ and then on $k=0,\dots,N-3-i$ for each fixed $i$. If we fix $(k,i)$ and let $V(t)=V_{k,i}(t)$ be the right hand side, then we have
\fm{dV/dt&=-B\eps t^{-1}\wt{E}_{k,i}(t)-C\eps t^{-1+C\eps}(\wt{E}_{k-1,i}(t)^{1/2}+\wt{E}_{k+1,i-1}(t)^{1/2})\wt{E}_{k,i}(t)^{1/2}\\&\geq -B\eps t^{-1}V(t)-C\eps t^{-1+C\eps}(\wt{E}_{k-1,i}(t)^{1/2}+\wt{E}_{k+1,i-1}(t)^{1/2})V(t)^{1/2}.}
Thus,
\fm{\frac{d}{dt}( t^{B\eps/2}\sqrt{V})&=\frac{1}{2}B\eps t^{-1+B\eps/2}\sqrt{V}+ t^{B\eps/2}\frac{dV/dt}{2\sqrt{V}}\\&=\frac{1}{2\sqrt{V}} t^{B\eps/2}(B\eps t^{-1}V+dV/dt)\\&\geq -C\eps t^{-1+(C+B/2)\eps}(\wt{E}_{k-1,i}(t)^{1/2}+\wt{E}_{k+1,i-1}(t)^{1/2})\\&\geq -C\eps^2 t^{-1+C\eps}T_1^{-1+C\eps}.}
The last line holds by induction hypothesis.  We then have
\fm{t^{B\eps/2}\sqrt{V(t)}&\leq T_1^{B\eps/2}\sqrt{V(T_1)}+\int_t^{T_1}C\eps^2\tau^{-1+C\eps}T_1^{-1/2+C\eps}\ d\tau\\&\leq C\eps T_1^{-1/2+C\eps},}
and thus for all $t\geq T_{N,R}$, we have
\fm{\wt{E}_{k,i}(t)\leq V(t)\leq t^{B\eps} V(t)\leq C\eps^2T_1^{-1+C\eps}.}
Here $C$ in different places may denote different values.

For $0\leq t\leq T_{N,R}$, we can also prove that \fm{\norm{\partial Z^I(v_2-v_1)(t)}_{L^2(\R^3)}\leq C_N\eps^2T_1^{-1+C_N\eps}.}The proof is very similar to the proof in Section \ref{smp}. We can use the equation 
\fm{\wt{g}^{\alpha\beta}(u_1)\partial_\alpha\partial_\beta (u_2-u_1)&=-(\gamma^{\alpha\beta}(u_2)-\gamma^{\alpha\beta}(u_1))\partial_\alpha\partial_\beta u_2} and apply the standard energy estimates to establish the continuity argument. Again, we can remove the dependence of $N$ in the constants, using the same argument  in Section \ref{smp}.

By \eqref{equivnorm}, for each $|I|\leq N-3$ and $\eps\ll 1$, we have 
\fm{\sup_{t\geq 0}\norm{\partial Z^I(v_2-v_1)(t)}_{L^2(\R^3)}\leq C\eps T_1^{-1/2+C\eps}\to 0} as $T_2>T_1\to\infty$. By the  Klainerman-Sobolev inequality and \fm{\int_{0}^{\infty} (1+t+\rho)^{-1}\lra{t-\rho}^{-1/2}\ d\rho\lesssim (1+t)^{-1/2},}for all $|I|\leq N-5$, we have
\fm{
\sup_{t\geq 0,\ x\in\R^3}|\partial Z^I (v_2-v_1)(t,x)|&\leq C\eps T_1^{-1/2+C\eps}\to 0\\
\sup_{t\geq 0,\ x\in\R^3}|Z^I (v_2-v_1)(t,x)|&\leq C\eps T_1^{-1/2+C\eps}\to 0,}as $T_2>T_1\to\infty$. Then, there is $v^\infty\in C^{N-4}(\{(t,x):\ t\geq 0\})$, such that $\partial Z^Iv^{T}\to \partial Z^Iv^\infty$ and $Z^Iv^T \to Z^Iv^\infty$ pointwisely for $t\geq 0$ as $T\to\infty$, for each $|I|\leq N-5$. It is clear that the pointwise bounds \eqref{ptb1} and \eqref{ptb2} also hold for $v^\infty$ for $|I|\leq N-5$. By Fatou's lemma, for each $|I|\leq N-5$ we have
\begin{equation}\label{vinfty}\norm{\partial Z^Iv^\infty(t)}_{L^2(\R^3)}\leq \liminf_{T\to\infty}\norm{\partial Z^Iv^T(t)}_{L^2(\R^3)}\leq C_I\eps( 1+t)^{-1/2+C_I\eps}.\end{equation}

Meanwhile, if $N\geq 6$, then by taking $T\to\infty$ in \fm{\wt{g}^{\alpha\beta}(u_{app}+v^T)\partial_\alpha\partial_\beta v^T=-\chi(t/T)\wt{g}^{\alpha\beta}(u_{app}+v^T)\partial_\alpha\partial_\beta u_{app},} we conclude that $u^\infty:=v^\infty+u_{app}$ is a solution to \eqref{qwe} for $t\geq 0$.

\subsection{End of the proof of Theorem \ref{mthm1}} \label{sec7.2}
If $t\gg_R1$ and $t\sim r$, we have $U(t,r,\omega)=U(\eps\ln(t)-\delta,R,\omega)$, so $\partial^kU=O(t^{-k})$; besides, $\partial^k(\psi(r/t))=O(t^{-k})$. Thus, for $t\gg_R1$ and for each $I$, we have 
\fm{|\partial Z^I(u_{app}-\eps r^{-1}U)(t,x)|\chi_{|x|\leq 3t/2}\lesssim_I\eps t^{-2+C_I\eps}}
and
\fm{&\hspace{1.25em}\norm{\partial Z^I(u_{app}-\eps r^{-1}U)(t)}_{L^2(\{x\in\R^3:\ |x|\leq 3t/2\})}\\&=\norm{\partial Z^I((1-\psi(r/t))\eps r^{-1}U)(t)}_{L^2(\{x\in\R^3:\ 5t/4\leq|x|\leq 3t/2\})}\\&\lesssim_I \eps t^{-2+C_I\eps}\cdot |\{x\in\R^3:\ 5t/4\leq|x|\leq 3t/2\}|^{1/2}\\&\lesssim_I \eps t^{-1/2+C_I\eps}.}
These two bounds allows us to get the estimates in the main theorem from \eqref{p7f1}, \eqref{p7f2} and \eqref{p7f3}, since
\fm{u-u_{app}&=(u-\eps r^{-1}U)\chi_{|x|\leq 3t/2}-(u_{app}-\eps r^{-1}U)\chi_{|x|\leq 3t/2}+u\chi_{|x|>3t/2}.}
We also remark that starting from the estimates in the main theorem, we can also derive \eqref{p7f1}, \eqref{p7f2} and \eqref{p7f3}, using the essentially same derivation here.

By \eqref{vinfty}, for $t\gg_R 1$, we have
\fm{|(\partial_t-\partial_r)(u^\infty-u_{app})(t,x)|\lesssim\eps t^{-1/2+C_I\eps}(1+t+r)^{-1}.}
Since $\psi\equiv 0$ unless $t\sim r$, for $t\gg_R 1$ we have
\fm{(\partial_t-\partial_r)u_{app}&=(\partial_t-\partial_r)(\eps r^{-1}\psi(r/t)U)\\&=-\eps r^{-2} \psi U+\eps r^{-1}\psi\mu U_q+\eps^2 r^{-1}t^{-1}\psi U_s\\&\hspace{1em}+\eps r^{-1}t^{-2}(t-r)\psi^\prime U\\&=-2\eps r^{-1}\psi A+O(\eps t^{-2+C\eps}), }
we conclude that  for all $t\gg_R 1$ we have
\begin{equation}
|(\partial_t-\partial_r)u^\infty+\frac{2\eps}{r}A(q(t,r,\omega),\omega)|\lesssim \eps t^{-3/2+C\eps}.
\end{equation}
Note that $A(q(t,r,\omega),\omega)=0$ unless $\psi(r/t)=1$ for $t\gg_R1$, so we do not have $\psi(r/t)$. So we gets the last estimates in the main theorem.

\subsection{Uniqueness} Now we give a brief proof of the uniqueness statement given in the remark of Theorem \ref{mthm1}. It suffices to prove the uniqueness of Proposition \ref{prop7}, assuming $N\geq 11$ and $\eps\ll1$. This is because \eqref{p7f1}, \eqref{p7f2} and \eqref{p7f3} are equivalent to the estimates in the main theorem, even if we replace $5/4$ with a fixed constant $\kappa>1$. We refer to Section \ref{sec7.2} for the proof.

Now, suppose we have two $C^{N-4}$ solutions $u_1,u_2$ constructed in Proposition \ref{prop7}. Fix $T\gg 1$. We can prove that $\norm{\partial Z^I(u_1-u_2)(t)}\lesssim\eps T^{-1/2+C\eps}$ for all $t\geq 0$ and $|I|\leq N-10$. Here the constants are independent of $T$. The proof is essentially the same as that in Section \ref{sec7.1}. Let $T\to\infty$ and we get $u_1\equiv u_2$.

\bibliography{paper}{}
\bibliographystyle{plain}

\end{document}